\documentstyle{article}

\def\dk#1{\mbox{$\lfloor {#1} \rfloor$}}

\def\gk#1{\mbox{$\lceil {#1} \rceil$}}

\def\ss{{\mbox{$/\!/$}}}

\def\bs{{\mbox{$\setminus\!\setminus$}}}

\def\cu{{\mbox{\it cups}}}

\def\ca{{\mbox{\it caps}}}

\def\eql{\equiv_{\cal L}}

\def\eqk{\equiv_{\cal K}}

\def\eqm{\equiv^{\cal J}}

\def\Lo{{\mbox{${\cal L}_\omega$}}}

\def\Lpm{{\mbox{${\cal L}_{\pm\omega}$}}}

\def\Ln{{\mbox{${\cal L}_n$}}}

\def\Lc{{\mbox{${\cal L}_c$}}}

\def\Lnc{{\mbox{${\cal L}_n^{\mbox{\scriptsize\rm cyl}}$}}}

\def\Ko{{\mbox{${\cal K}_\omega$}}}

\def\Kpm{{\mbox{${\cal K}_{\pm\omega}$}}}

\def\Kn{{\mbox{${\cal K}_n$}}}

\def\Kc{{\mbox{${\cal K}_c$}}}

\def\Knc{{\mbox{${\cal K}_n^{\mbox{\scriptsize\rm cyl}}$}}}

\def\M{{\mbox{${\cal J}$}}}

\def\Mo{{\mbox{${\cal J}_\omega$}}}

\def\Mn{{\mbox{${\cal J}_n$}}}

\def\Mc{{\mbox{${\cal J}_c$}}}

\def\Fn{{\mbox{${\cal F}_n$}}}

\def\Bn{{\mbox{${\cal B}_n$}}}

\def\X{{\mbox{${\cal X}$}}}

\def\Mat{{\mbox{${\mbox{\bf Mat}}_{\cal F}$}}}

\def\End{\mbox{\bf End}}

\newcommand{\mj}{\mbox{\bf 1}}

\newcommand{\str}{\rightarrow}

\newcommand{\N}{{\mbox{\bf N}}}

\newcommand{\Z}{{\mbox{\bf Z}}}

\newcommand{\qed}{\hfill $\Box$}

\def\cirk{\raisebox{1pt} {\scriptsize $\;\circ \;$}}

\begin{document}

\title
{\bf Self-Adjunctions and Matrices}

\author{
{\sc Kosta Do\v sen} and {\sc Zoran Petri\' c}
\\[.05cm]
\\Mathematical Institute, SANU
\\Knez Mihailova 35, P.O. Box 367
\\11001 Belgrade, Yugoslavia
\\email: \{kosta, zpetric\}@mi.sanu.ac.yu}

\date{ }
\maketitle
\begin{abstract}
\noindent It is shown that the multiplicative monoids of
Temperley-Lieb algebras generated out of the basis are isomorphic
to monoids of endomorphisms in categories where an endofunctor is
adjoint to itself. Such a self-adjunction is found in a category
whose arrows are matrices, and the functor adjoint to itself is
based on the Kronecker product of matrices. This self-adjunction
underlies the orthogonal group case of Brauer's representation of
the Brauer centralizer algebras. \vspace{0.3cm}



\end{abstract}

\section{Introduction}

As an offshoot of Jones' polynomial approach to knot and link invariants,
Temperley-Lieb algebras have played in the 1990s a prominent
role in knot theory and low-dimensional topology (see \cite{KL}, \cite{L97}
and \cite{PS}).
In this paper we show that the multiplicative monoids of
Temperley-Lieb algebras are closely related to the general notion of
adjunction, one of the fundamental notions of category theory, and of
mathematics in general
(see \cite{ML71}). More precisely, we show that these monoids
are isomorphic to monoids of endomorphisms in categories
involved in one kind of self-adjoint
situation, where an endofunctor is adjoint to itself.

Early work on Temperley-Lieb algebras established the importance
of a self-dual object in a monoidal category for understanding
the categorial underpinnings of the matter (see \cite{Y88}, \cite{FY89},
\cite{JS93} and papers cited therein). The result of the present
paper is in the wake of this earlier work, and shouldn't be
surprising.

We find a self-adjunction in categories whose arrows are
matrices, where the functor adjoint to itself is based on the
Kronecker product of matrices.
This self-adjunction underlies the orthogonal group case
of Brauer's representation of the Brauer algebras,
which can be restricted to the Temperley-Lieb subalgebras
of the Brauer algebras (see \cite{B37}, \cite{W88}, Section 3, and
\cite{J94}, Section 3).
Thereby, building on ideas similar to, but not quite the same as,
those that lead to the representation of braid groups in
Temperley-Lieb
algebras, which is due to Jones,
we obtain a representation of braid groups in matrices.
The question
whether the representation of braid groups in Temperley-Lieb algebras
is faithful is raised in \cite{J00}.
(The review \cite{B98} provides a good survey of questions about the linear
representation of braid groups.)

The representation of monoids of Temperley-Lieb algebras in matrices
provides a faithful, i.e. isomorphic, representation in matrices
of a certain brand
of these algebras; this is the representation mentioned above, which
originates in \cite{B37}.
The faithfulness of this representation is established in
\cite{J94} (Section 3; an elementary self-contained proof of the same
fact may also be found in \cite{DKP}).

Although this is a paper in the theory of Temperley-Lieb
algebras, reading it does not presuppose acquaintance
with works of that particular theory, except for the sake of motivation.
In the latter sections of our paper, where we deal with matters of category
theory, we presuppose some, rather general, acquaintance with a few
notions from that field, which are all explained in \cite{ML71}.

In the paper we proceed roughly as follows. We first
present by generators and relations monoids
for which we will show later that they
are engendered by categories involved in self-adjoint
situations. These categories engender monoids,
whose names will be indexed by $\omega$, when we consider
a total binary operation on all arrows defined with
the help of composition.
Our categories engender monoids of a different kind,
with names indexed by $n$, when we consider just composition,
restricting ourselves to
endomorphisms in the category. We deal first
with the monoids related to the most general notion of self-adjunction,
which we tie to the
label $\cal L$, and next with those related to the more particular notion
of self-adjunction,
tied to the label $\cal K$, which we encounter in connection with
Temperley-Lieb algebras.

Next we prove in detail that our monoids are isomorphic to monoids made
of equivalence classes of diagrams which in knot theory would
be called planar tangles, without crossings, and which we call
{\it friezes}.
In these representations, there are two different notions
of equivalence of friezes: the $\cal L$ notion is based purely
on planar ambient isotopies, whereas the $\cal K$ notion allows
circles to cross lines, which is forbidden in the $\cal L$ notion.
So the mathematical content of the most general notion of
self-adjunction is caught by the notion of planar ambient
isotopy. The diagrammatic representation of the
\Kn\ monoids is not an entirely new result, but,
though some of the facts are well known, we have
not been able to find in the literature a self-contained and complete
treatment of the matter, such as we will try to give.

A theorem connecting self-adjunction to diagrams analogous to
friezes is stated without proof in \cite{FY89} (Theorem 4.1.1,
p. 172). It is not clear, however, whether $\cal L$ or $\cal K$
notions are meant (it seems the self-adjunction is of the $\cal L$
kind, while the diagrams are of the $\cal K$ kind).  Adjunction, as
it occurs in symmetric monoidal closed categories,
is connected to diagrams like ours in \cite{EK66} and \cite{KM71}.
This connection between adjunction and diagrams is made also,
more or less explicitly, in many of the papers mentioned above,
about categorial matters tied to Temperley-Lieb algebras; in
particular, in \cite{Y88} and \cite{FY92}
(see also \cite{KD99}, \S 4.10). The completeness results of
\cite{Y88} and \cite{FY92} are established by showing that
each equivalence move of the Reidemeister kind in the diagrams
corresponds to an equation. Our style is somewhat different,
because we rely on normal forms that are not mentioned
in the other approach.

We consider also a third notion of self-adjunction, $\cal J$
self-ajunction,
in whose diagrammatic representation we don't take account of
circles at all.
This notion is more strict
than $\cal K$ self-adjunction. With the help of
friezes, we show for this third notion that it
is maximal in the sense that we could not extend it with any further
assumption
without trivializing it. This maximality is an essential
ingredient in our proof that we have in matrices an isomorphic
representation of the monoids of Temperley-Lieb algebras.
However, this maximality need not serve only for that
particular goal, which can be reached by other means, as we mentioned
above (see \cite{J94}, Section 3, and \cite{DKP}), or by relying
on \cite{DP02}. Maximality
can serve to establish the isomorphism of other nontrivial
representations of the monoids of Temperley-Lieb algebras.

We deal with matters involving categories towards the end of the
paper. There we introduce our categories of matrices, and exhibit
the $\cal K$ self-adjunction involved in them.
We deal with the orthogonal group case of Brauer's representation
and with our
representation of braid groups in matrices in the last two sections.

As an aside,
we consider in Section 13 monoids interpreted in
friezes with points labeled
by all integers, and not only positive integers, and also monoids
interpreted in cylindrical friezes. This matter, though related to
other matters in the paper, is independent of its main thrust,
and is presented with less details.

In the main body of the paper, however, we strive,
as we said above, to give detailed proofs of
isomorphisms of monoids, and so our style of exposition will be occasionally
rather formal. It will be such at the beginning, and it might help the
reader while going through Sections 2-5 to take a look at Sections 6-8,
and perhaps also at Sections 14-16, to get some motivation.

\vspace{.3cm}

\noindent {\sc Contents}
\begin{tabbing}
1. \hspace{1em} \= {\it Introduction}
\\
2. \> {\it The monoids \Lo\ and \Ko}
\\
3. \> {\it Finite multisets, circular forms and ordinals}
\\
4. \> {\it Normal forms in \Lo}
\\
5. \> {\it Normal forms in \Ko}
\\
6. \> {\it Friezes}
\\
7. \> {\it Generating friezes}
\\
8. \> {\it \Lo\ and \Ko\ are monoids of friezes}
\\
9. \> {\it The monoids \Ln}
\\
10. \> {\it The monoids \Kn}
\\
11. \> {\it The monoid \Mo}
\\
12. \> {\it The maximality of \Mo}
\\
13. \> {\it The monoids \Lpm\ and \Kpm}
\\
14. \> {\it Self-adjunctions}
\\
15. \> {\it Free self-adjunctions}
\\
16. \> {\it \Lc\ and \Lo}
\\
17. \> {\it Self-adjunction in \Mat}
\\
18. \> {\it Representing \Mo\ in \Mat}
\\
19. \> {\it Representing \Kc\ in \Mat}
\\
20. \> {\it The algebras $\End(p^n)$}
\\
21. \> {\it Representing braid groups in $\End(p^n)$}
\\
    \> {\it References}
\end{tabbing}

\section{The monoids \Lo\ and \Ko}

The monoid \Lo\ has for every $k\in \N^+=\N\!-\!\{0\}$
a generator \dk{k}, called a
{\it cup}, and a generator \gk{k}, called a {\it cap}. The {\it terms} of
\Lo\ are defined inductively by stipulating that the generators and \mj\ are
terms, and that if $t$ and $u$ are terms, then $(tu)$ is a term. As usual, we
will omit the outermost parentheses of terms. In the presence of
associativity we will omit all parentheses, since they can be restored as we
please.

The monoid \Lo\ is freely generated from the generators above so that the
following equations hold between terms of \Lo\ for $l\leq k$:
\[
\begin{array}{rl}
(1) & \mj t=t, \quad t\mj=t,
\\[0.1cm]
(2) & t(uv)=(tu)v,
\\[0.2cm]
({\mbox{\it cup}}) & \dk{k}\dk{l} = \dk{l}\dk{k+2},
\\[0.1cm]
({\mbox{\it cap}}) & \gk{l}\gk{k} = \gk{k+2}\gk{l},
\\[0.2cm]
({\mbox{\it cup-cap}}\;1) & \dk{l}\gk{k+2} = \gk{k}\dk{l},
\\[0.1cm]
({\mbox{\it cup-cap}}\;2) & \dk{k+2}\gk{l} = \gk{l}\dk{k},
\\[0.1cm]
({\mbox{\it cup-cap}}\;3) & \dk{k}\gk{k\pm 1} = \mj.
\end{array}
\]

The monoid \Ko\ is defined as the monoid \Lo\ save that we have
the additional equation
\[
\begin{array}{ll}
({\mbox{\it cup-cap}}\; 4) & \dk{k}\gk{k} = \dk{k+1}\gk{k+1},
\end{array}
\]
which, of course, implies
\[ \dk{k}\gk{k} = \dk{l}\gk{l}.\]

\noindent To understand the equations of \Lo\ and \Ko\ it helps to
have in mind their diagrammatic interpretation of Sections 6-8
(see in particular the diagrams corresponding to \dk{k} and \gk{k}
at the beginning of Section 7).

Let $[k]$ be an abbreviation for $\dk{k}\gk{k}$, and let us call such terms
{\it circles}. Then ({\it cup-cap} 4) says that we have only one circle,
which we
designate by $c$. We have the following equations in \Lo\ for $l\leq k$:
\[
\begin{array}{l}
\dk{k}[l] = [l]\dk{k},
\\[0.1cm]
\dk{l}[k+2] = [k]\dk{l}.
\end{array}
\]
For the first equation we have
\[
\begin{array}{rl}
\dk{k}\dk{l}\gk{l} = & \dk{l}\dk{k+2}\gk{l},\quad {\mbox{\rm by ({\it cup})}}
\\[0.1cm]
= & \dk{l}\gk{l}\dk{k},\quad {\mbox{\rm by ({\it cup-cap} 2)}},
\end{array}
\]
and for the second
\[
\begin{array}{rl}
\dk{l}\dk{k+2}\gk{k+2} = & \dk{k}\dk{l}\gk{k+2},\quad {\mbox{\rm by ({\it cup})}}
\\[0.1cm]
= & \dk{k}\gk{k}\dk{l},\quad {\mbox{\rm by ({\it cup-cap} 1)}}.
\end{array}
\]
We derive analogously the following dual equations of \Lo\ for $l\leq k$:
\[
\begin{array}{l}
[l]\gk{k} = \gk{k}[l],
\\[0.1cm]
[k+2]\gk{l} = \gk{l}[k].
\end{array}
\]
So in \Ko\ we have the equations
\[
\begin{array}{c}
\dk{k}c=c\dk{k},
\\
\gk{k}c=c\gk{k},
\end{array}
\]
which yield the equation $tc=ct$ for any term $t$.

\section{Finite multisets, circular forms and ordinals}

Let an $o$-{\it monoid} be a monoid with an arbitrary unary operation $o$, and consider
the free {\it commutative} $o$-monoid $\cal F$ generated by the empty set of generators.
In $\cal F$ the operation $o$ is a one-one function.

The elements of $\cal F$ may be designated by parenthetical words, i.e. well-formed
words in the alphabet $\{ (,)\}$, which will be precisely defined in a moment,
where the empty word stands for the unit of the monoid, concatenation is monoid
multiplication, and $o(a)$ is written simply $(a)$. Parenthetical words are defined
inductively as follows:

\vspace{.2cm}

(0)\quad the empty word is a parenthetical word;

(1)\quad if $a$ is a parenthetical word, then $(a)$ is a parenthetical word;

(2)\quad if $a$ and $b$ are parenthetical words, then $ab$ is a
parenthetical word.

\vspace{.2cm}

We consider next several isomorphic representations of $\cal F$, via finite
multisets, circular forms in the plane and ordinals.

If we take that $()$ is the empty multiset, then the elements of $\cal F$ of
the form $(a)$ may be identified with finite multisets, i.e. the hierarchy
of finite multisets obtained by starting from the empty multiset $\emptyset$ as
the only {\it urelement}. To obtain a more conventional notation for these multisets, just
replace $()$ everywhere by $\emptyset$, replace the remaining left parentheses $($
by left braces $\{$ and the remaining right parentheses $)$ by right braces $\}$,
and put in commas where concatenation occurs.

The elements of $\cal F$ may also be identified with nonintersecting
finite collections
of circles in the plane factored through homeomorphisms of the plane mapping
one collection into another (cf. \cite{KKL}, Section II).
For this interpretation, just replace $(a)$ by
\begin{picture}(10,10)(1,3)
\put(5,5){\circle{10}}
\put(5,5){\makebox(0,0){$a$}}
\end{picture}.
Since we will be interested in particular in this plane interpretation,
we
call the elements of $\cal F$ {\it circular forms}. The {\it empty circular form}
is the unit of $\cal F$. When we need to refer to it we use $e$. We refer
to other circular forms with parenthetical words.

The free commutative $o$-monoid $\cal F$ has another isomorphic representation
in the ordinals contained in the ordinal $\varepsilon_0=\min\{\xi \mid
\omega^\xi=\xi\}$, i.e. in the ordinals lesser than $\varepsilon_0$.
By Cantor's
Normal Form Theorem (see, for example,
\cite{KM}, VII.7, Theorem 2, p. 248, or \cite{L79}, IV.2, Theorem
2.14, p. 127), for every ordinal $\alpha > 0$ in $\varepsilon_0$
there is a unique finite ordinal
$n\geq 1$ and a unique sequence of ordinals
$\alpha_1\geq\ldots\geq\alpha_n$ contained in $\alpha$, i.e. lesser than
$\alpha$, such that $\alpha=\omega^{\alpha_1}+\ldots +\omega^{\alpha_n}$.
The {\it natural sum} $\alpha\:\sharp\:\beta$ of
\[
\begin{array}{c}
\alpha=\omega^{\alpha_1}+\ldots +\omega^{\alpha_n},\quad\alpha_1\geq\ldots\geq\alpha_n,
\\
\beta=\omega^{\beta_1}+\ldots +\omega^{\beta_m},\quad\beta_1\geq\ldots\geq\beta_m,
\end{array}
\]
is defined as $\omega^{\gamma_1}+\ldots +\omega^{\gamma_{n+m}}$ where
$\gamma_1,\ldots,\gamma_{n+m}$ is obtained by permuting the sequence
$\alpha_1,\ldots,\alpha_n,\beta_1,\ldots,\beta_m$ so that $\gamma_1\geq\ldots\geq
\gamma_{n+m}$ (this operation was introduced by Hessenberg; see
\cite{KM}, p. 252, or \cite{L79}, p.130). We also have
$\alpha\:\sharp\: 0 =0\:\sharp\:\alpha=\alpha$. The natural sum $\sharp$ and the
ordinal sum $+$ don't coincide in general: $\sharp$ is commutative, but $+$ is
not (for example, $\omega=1+\omega\neq\omega+1$, but $1\:\sharp\:\omega=
\omega\:\sharp\: 1=\omega+1=\omega^{\omega^0}+\omega^0$). However, if
$\alpha_1\geq\ldots\geq\alpha_n$, then
$\omega^{\alpha_1}+\ldots +\omega^{\alpha_n}=\omega^{\alpha_1}\:\sharp\:
\ldots \:\sharp\:\omega^{\alpha_n}$.

Let $\omega^{\ldots}$ be the unary operation that assigns to every
$\alpha\in\varepsilon_0$ the ordinal $\omega^\alpha\in\varepsilon_0$. Then
it can be shown that the commutative $o$-monoid
$\langle\varepsilon_0,\sharp,0,\omega^{\ldots}\rangle$ is isomorphic to $\cal F$
by the isomorphism $\iota :\varepsilon_0\str {\cal F}$ such that
$\iota(0)$ is the empty word and
\[
\iota(\omega^{\alpha_1}+\ldots+
\omega^{\alpha_n})=\iota(\omega^{\alpha_1}\:\sharp\:\ldots\:\sharp\:\omega^{\alpha_n})
=(\iota(\alpha_1))\ldots(\iota(\alpha_n)).
\]
That the function $\iota^{-1}:{\cal F}\str\varepsilon_0$
defined inductively by
\[
\begin{array}{l}
\iota^{-1}(e)=0,
\\
\iota^{-1}(ab)=\iota^{-1}(a)\:\sharp\:\iota^{-1}(b),
\\
\iota^{-1}((a))=\omega^{\iota^{-1}(a)}
\end{array}
\]
is the inverse of $\iota$ is established by easy inductions relying
on the fact that
\[
\iota(\alpha\:\sharp\:\beta)=\iota(\alpha)\iota(\beta).
\]

It is well known in proof theory that
the ordinal $\varepsilon_0$ and natural sums play an important role in
Gentzen's
proof of the consistency of formal Peano arithmetic PA (see
\cite{G69}, Paper 8, \S4). Induction up to any ordinal lesser
than $\varepsilon_0$ is derivable in PA; induction up to
$\varepsilon_0$, which is not derivable in PA, is not only
sufficient, but also necessary, for proving the consistency of PA
(see \cite{G69}, Paper 9).

From the isomorphism of $\cal F$ with $\langle\varepsilon_0,\sharp,0,
\omega^{\ldots}\rangle$
we obtain immediately a normal form for the elements of $\cal F$. Circular forms
inherit a well-ordering from the ordinals, and we have the following inductive
definition. The empty word is in normal form, and if $a_1,\ldots,
a_n$, $n\geq 1$, are parenthetical words in normal form such that
$a_1\geq\ldots\geq a_n$, then $(a_1)\ldots (a_n)$ is in normal form. We call
this normal form of parenthetical words the {\it Cantor normal form}.

Let a commutative $o$-monoid be called {\it solid} iff it satisfies
\[
({\mbox{\it solid}})\quad o(a)=o(1) a,
\]
where $1$ is the unit of the monoid. The free solid commutative $o$-monoid
${\cal F}'$ generated by the empty set of generators is isomorphic to the structure
$\langle\N ,+,0,\mbox{\scriptsize\ldots}\!+\!1\rangle$
by the isomorphism that assigns to $n$ the
sequence of $n$ pairs $(\;)$.
So ({\it solid}) makes $\langle\varepsilon_0,\sharp,\emptyset,\omega^{\ldots}
\rangle$ collapse into
$\langle\omega,\sharp,\emptyset,\mbox{\scriptsize\ldots}\!+\!1\rangle$.
For $k\in \N$, let
$k\N =\{kn\mid n\in\N\}$ and $k^{\mbox{\scriptsize{\bf N}}}
=\{k^n\mid n\in\N\}$. If $k\geq 1$, then
$\langle\N,+,0,{\mbox{\scriptsize\ldots}}\!+\!1\rangle$ is
isomorphic to $\langle k\N,+,0,\mbox{\scriptsize\ldots}\!+k
\rangle$, which for $k\geq 2$ is isomorphic to
$\langle k^{\mbox{\scriptsize{\bf N}}}, \cdot,1,
{\mbox{\scriptsize\ldots}}\!\cdot k\rangle$.

The equation ({\it solid}) is what a unary function $o:{\cal M}\str{\cal M}$,
for a monoid $\cal M$,
has to satisfy to be in the
image of the Cayley monomorphic representation of $\cal M$ in
${\cal M}^{\cal M}$, which assigns to every $a\in
{\cal M}$ the function $f_a\in
{\cal M}^{\cal M}$ such that $f_a(b)=ab$. In the presence of ({\it solid}),
the function $f_{o(a)}$ will be equal to $o\cirk f_a$.
The equation ({\it solid})
can be replaced by $o(ab)=o(a)b$,
and in commutative $o$-monoids it could, of course, as well be written
$o(a)=ao(1)$.

\section{Normal forms in \Lo}

For $k\in\N^+$ let $c_k^0$ be the term $\mj$ of \Lo. For $\alpha>0$
an ordinal in $\varepsilon_0$ whose Cantor normal form is
$\omega^{\alpha_1}+\ldots+\omega^{\alpha_n}$ let the term $c_k^{\alpha}$ of \Lo\ be
defined inductively as
\[
\dk{k}c_{k+1}^{\alpha_1}\gk{k}\ldots\dk{k}c_{k+1}^{\alpha_n}\gk{k}.
\]
Next, let $a_k^0$ be the term \dk{k}, and let $a_k^{\alpha}$ be the term
$\dk{k}c_{k+1}^{\alpha}$. Similarly, let $b_k^0$ be the term \gk{k}, and let $b_k^{\alpha}$
be the term $c_{k+1}^{\alpha}\gk{k}$.

Consider terms of \Lo\ of the form
\[
b_{j_1}^{\beta_1}\ldots b_{j_m}^{\beta_m}c_{k_1}^{\gamma_1}\ldots
c_{k_l}^{\gamma_l}a_{i_1}^{\alpha_1}\ldots a_{i_n}^{\alpha_n}
\]
where $n,m,l\geq 0$, $n+m+l\geq 1$, $j_1>\ldots>j_m$, $k_1<\ldots<k_l$,
$i_1<\ldots<i_n$, and for every $p\in\{1,\ldots,l\}$ we have $\gamma_p\neq 0$.
If $n$ is 0, the sequence $a_{i_1}^{\alpha_1}\ldots a_{i_n}^{\alpha_n}$ is empty,
and analogously if $m$ or $l$ is 0. Terms of \Lo\ of this form
and the term \mj\
will be said to be in {\it normal form}.

In the definition of normal form we could have required that $k_1>\ldots>k_n$, or,
as a matter of fact, we could have imposed any other order on these particular
indices, with the same effect. We have chosen the order above for the sake of
definiteness. (Putting aside complications involving the terms
$c_k^{\alpha}$ and the ordinals, the idea of our normal form may be
found in \cite{BJ97}, p. 106.)

To reduce terms of \Lo\ to normal form we use an alternative formulation of \Lo,
which is obtained as follows. Now the generators are the $a$ {\it terms}
$a_k^{\alpha}$, the $b$ {\it terms} $b_k^{\alpha}$ and the $c$ {\it terms}
$c_k^{\alpha}$ for $k\in \N^+$ and $\alpha\in\varepsilon_0$.
These terms are now primitive, and not defined. We generate terms with these generators,
\mj\ and multiplication, and we stipulate the following equations for
$l\leq k$:
\[
\begin{array}{ll}
{\makebox[1cm][l]{$(1)$}} & {\makebox[6cm][l]{$\mj t=t, \quad t\mj=t,$}}
\\[.1cm]
(2) & t(uv)=(tu)v,
\\[.1cm]
(aa) & a_k^{\alpha}a_l^{\beta}=a_l^{\beta}a_{k+2}^{\alpha},
\\[.1cm]
(bb) & b_l^{\alpha}b_k^{\beta}=b_{k+2}^{\beta}b_{l}^{\alpha},
\\[.1cm]
(c1) & c_k^0=\mj,
\\[.1cm]
(c2) & c_k^{\alpha}c_k^{\beta}=c_k^{\alpha\sharp\beta},
\\[.1cm]
(cc) & c_k^{\alpha}c_l^{\beta}=c_l^{\beta}c_k^{\alpha},\quad
{\mbox{\rm for }} l< k,
\end{array}
\]
$ab$ {\it equations}:
\[
\begin{array}{ll}
{\makebox[1cm][l]{$(ab\:1)$}} & {\makebox[6cm][l]
{$a_l^{\alpha}b_{k+2}^{\beta}=b_k^{\beta}a_l^{\alpha},$}}
\\[.1cm]
(ab\:2) & a_{k+2}^{\alpha}b_{l}^{\beta}=b_l^{\beta}a_k^{\alpha},
\\[.1cm]
(ab\:3.1) & a_k^{\alpha}b_{k+1}^{\beta}=c_k^{\beta}c_{k+1}^{\alpha},
\\[.1cm]
(ab\:3.2) & a_{k+1}^{\alpha}b_k^{\beta}=c_k^{\alpha}c_{k+1}^{\beta},
\\[.1cm]
(ab\:3.3) & a_k^{\alpha}b_k^{\beta}=c_k^{\omega^{\alpha\sharp\beta}},
\end{array}
\]
$ac$ {\it equations}:
\[
\begin{array}{ll}
{\makebox[1cm][l]{$(ac\:1)$}} &
{\makebox[6cm][l]{$a_k^{\alpha}c_l^{\gamma}=c_l^{\gamma}a_k^{\alpha},$}}
\\[.1cm]
(ac\:2) & a_l^{\alpha}c_{k+2}^{\gamma}=c_k^{\gamma}a_l^{\alpha},
\\[.1cm]
(ac\:3) & a_k^{\alpha}c_{k+1}^{\gamma}=a_k^{\alpha\sharp\gamma},
\end{array}
\]
$bc$ {\it equations}:
\[
\begin{array}{ll}
{\makebox[1cm][l]{$(bc\:1)$}} & {\makebox[6cm][l]
{$c_l^{\gamma}b_k^{\beta}=b_k^{\beta}c_l^{\gamma},$}}
\\[.1cm]
(bc\:2) & c_{k+2}^{\gamma}b_l^{\beta}=b_l^{\beta}c_k^{\gamma},
\\[.1cm]
(bc\:3) & c_{k+1}^{\gamma}b_k^{\beta}=b_k^{\gamma\sharp\beta}.
\end{array}
\]

It is tiresome, but pretty
straightforward, to derive all these equations in the
original formulation of \Lo\ for defined $c_k^{\alpha}$, $a_k^{\alpha}$ and
$b_k^{\alpha}$, while with \dk{k}\ defined as $a_k^0$ and \gk{k}\ defined as
$b_k^0$, we easily derive in the new formulation the equations of the original
formulation of \Lo. We can, moreover, derive in the new formulation the inductive
definitions of $c_k^{\alpha}$, $a_k^{\alpha}$ and $b_k^{\alpha}$.
We can then prove the following lemma for \Lo.

\vspace{.2cm}

\noindent {\sc Normal Form Lemma.}\quad {\it Every term is equal to a term in
normal form.}

\vspace{.1cm}

\noindent {\it Proof.}\quad We will give a reduction procedure that transforms every term
$t$ of \Lo\ into a term $t'$ in normal form such that $t=t'$ in \Lo.
(In logical jargon, we establish that this procedure is strongly
normalizing---namely, that any sequence of reduction steps terminates
in a term in normal form.)

Take a term $t$ in the new alternative formulation of \Lo, and let subterms
of this term of the forms on the left-hand sides of the equations of the alternative
formulation except (2) be called {\it redexes}. A reduction consists in replacing
a redex of $t$ by the term on the right-hand side of the corresponding equation.
Note that the terms on the left-hand sides of these equations cover all
possible cases for terms of the forms $a_k^{\alpha}b_l^{\beta}$,
$a_k^{\alpha}c_k^{\gamma}$
and $c_l^{\gamma}b_k^{\beta}$, and all these cases exclude each
other.

A subterm of $t$ which is an $a$ term will be called an $a$ {\it subterm}
of $t$, and analogously with $b$ and $c$.
For a particular subterm $a_k^{\alpha}$ of $t$ let $\sigma(a_k^{\alpha})$ be the
number of $b$ subterms of $t$ on the right-hand side of $a_k^{\alpha}$ in $t$. Let $n_1$
be the sum of all the numbers $\sigma(a_k^{\alpha})$ for every $a$ subterm
$a_k^{\alpha}$ of $t$. If there are no $a$ subterms of $t$,
then $n_1$ is zero.

For a particular subterm $a_k^{\alpha}$ of $t$ let $\sigma_a(a_k^{\alpha})$
be the number of $a$ subterms $a_l^{\beta}$ of $t$ on the right-hand side
of $a_k^{\alpha}$ in $t$ such that $l\leq k$. Let $\sigma_a$ be the
sum of all the numbers $\sigma_a(a_k^{\alpha})$ for every $a$ subterm
$a_k^{\alpha}$ of $t$.
For a particular subterm $b_k^{\beta}$ of $t$ let $\sigma_b(b_k^{\beta})$
be the number of $b$ subterms $b_l^{\alpha}$ of $t$ on the left-hand side of
$b_k^{\beta}$ in $t$ such that $l\leq k$. Let $\sigma_b$ be the sum
of all the numbers $\sigma_b(b_k^{\beta})$ for every $b$ subterm
$b_k^{\beta}$ of $t$.
For a particular subterm $c_l^{\gamma}$ of $t$ let $\tau(c_l^{\gamma})$
be the number of $a$ subterms on the left-hand side $c_l^{\gamma}$ in $t$
plus the number of $b$ subterms on the right-hand side of $c_l^{\gamma}$
in $t$. Let $\tau$ be the sum of all the numbers $\tau(c_l^{\gamma})$
for every $c$ subterm $c_l^{\gamma}$ of $t$. Let $\nu_c$ be the number of $c$ subterms
of $t$, and let $\nu_{\mbox{\scriptsize\bf 1}}$
be the number of subterms $\mj$ of $t$.
Then let $n_2$ be $\sigma_a+\sigma_b+\tau+2\nu_c+\nu_{\mbox
{\scriptsize\bf 1}}$.

Let $\sigma_c$ be defined as $\sigma_a$ save that $a$ is everywhere
replaced by $c$, and let $n_3$ be $\sigma_c$.
With reductions based on $(c2)$ and $(cc)$, the number
$n_3$ decreases, while $n_1$ and $n_2$ don't increase.
With reductions based on $(1)$, $(c1)$, $(aa)$, $(bb)$ and the $ac$ and $bc$
equations, $n_2$ decreases, $n_1$ doesn't change, and $n_3$
may even increase
in case we apply $(ac\:2)$ or $(bc\:2)$. With reductions based on
the $ab$ equations,
$n_1$ decreases, while $n_2$ and $n_3$ may increase. Then we take as
the {\it complexity measure} of $t$ the ordered triple $(n_1,n_2,n_3)$.
These triples are well-ordered lexicographically, and with every
reduction the complexity measure decreases. So by induction on the complexity measure,
we obtain that every term in the new formulation is equal to a term without redexes,
and it is easy to check that such a term stands for a term in normal form of the
original formulation of \Lo.
\qed

\vspace{.2cm}

\section{Normal forms in \Ko}

Let $c$ stand for $[k]$, where $k \in \N^+$.
Let $c^0$ be the empty sequence, and let $c^{n+1}$ be $c^nc$. Consider
terms of \Lo\ of the form
\[
\gk{j_1}\ldots\gk{j_m}c^l\dk{i_1}\ldots\dk{i_n}
\]
where $n,m,l\geq 0$, $n+m+l\geq 1$, $j_1>\ldots>j_m$ and
$i_1<\ldots<i_n$. Terms of this form and the term $\mj$ will
be said to be in
$\cal K$-{\it normal form}. (We could as well put $c^l$ on the extreme left, or on the
extreme right, or, actually, anywhere, but for the sake of definiteness,
and, by analogy with the normal form of \Lo, we put $c^l$ in the middle.)

We can easily derive from the Normal Form Lemma for \Lo\
the Normal Form Lemma for \Ko, which says that every term is
equal in \Ko\ to a term in $\cal K$-normal form. For that it is
enough to use the uniqueness of $c$ and
$tc=ct$. However, the Normal Form Lemma for \Ko\ has a much simpler
direct proof,
which does not require the introduction of an alternative formulation of \Ko.
This proof is obtained by simplifying the proof of the Normal Form Lemma
for \Lo. The complications of the previous proof were all due to distinguishing
$[k]$ from $[k+1]$ and to the absence of $tc=ct$. In \Ko\ we have in fact assumed
({\it solid}), and the ordinals in $\varepsilon_0$ have collapsed
into natural numbers.

\section{Friezes}

A one-manifold with boundary is a
topological space whose points have
open neighbourhoods homeomorphic to the real intervals $(-1,1)$
or $[0,1)$, the boundary points having the latter kind of neighbourhoods.
For $a>0$ a real number, let $R_a$ be $[0,\infty)\times[0,a]$. Let
$\{(x,a)\mid x\geq 0\}$ be the {\it top} of $R_a$ and $\{(x,0)\mid x\geq 0\}$
the {\it bottom} of $R_a$.

An $\omega$-{\it diagram} $D$ in $R_a$ is a one-manifold with
boundary with denumerably many compact connected components
embedded in $R_a$ such that the intersection of $D$ with the top
of $R_a$ is $t(D)=\{(i,a)\mid i\in\N^+\}$ the intersection of $D$
with the bottom of $R_a$ is $b(D)=\{(i,0)\mid i\in\N^+\}$ and
$t(D)\cup b(D)$  is the set of boundary points of $D$.

It follows from this definition that every $\omega$-diagram has denumerably many
components homeomorphic to $[0,1]$, which are called {\it threads}, and at most
a denumerable number of components homeomorphic to $S^1$, which are called
{\it circular components}. The threads and the circular components make all the
connected components of an $\omega$-diagram. All these components are mutually
disjoint.
Every thread has two {end points} that belong to the boundary $t(D)\cup b(D)$.
When one of these end points is in $t(D)$ and the other in $b(D)$, the thread is
{\it transversal}. A transversal thread is {\it vertical} when the first coordinates
of its end points are equal. A thread that is not transversal is a {\it cup} when
both of its end points are in $t(D)$, and it is a {\it cap} when they are both
in $b(D)$.

A {\it frieze} is an $\omega$-diagram with a finite number of cups,
caps and circular
components. Although many, but not all, of the definitions that follow can be
formulated for all $\omega$-diagrams, and not only for friezes, we
will be interested here only in friezes, and we will formulate our
definitions only with respect to them.
The notion of frieze corresponds to a special kind of {\it tangle} of
knot theory,
in which there are no crossings (see \cite{BZ85}, p. 99,
\cite{M96}, Chapter 9, \cite{K95}, Chapter 12).

For $D_1$ a frieze in $R_a$ and $D_2$ frieze in $R_b$, we say that $D_1$ is
$\cal L$-{\it equivalent} to $D_2$, and write $D_1\cong_{\cal L} D_2$, iff
there is a homeomorphism
$h:R_a\str R_b$ such that $h[D_1]=D_2$ and for every $i\in\N^+$ we have
$h(i,0)=(i,0)$ and $h(i,a)=(i,b)$. It is straightforward to check that
$\cal L$-equivalence between friezes is indeed an equivalence relation.

This definition is equivalent to a definition of $\cal L$-equivalence in terms of ambient
isotopies. The situation is analogous to what one finds in knot theory, where
one can define equivalence of knots either in terms of ambient isotopies or in
a simpler manner, analogous to what we have in the preceding
paragraph. The equivalence of these two
definitions is proved with the help of Alexander's trick
(see \cite{BZ85}, Chapter 1B), an adaptation of which also works in
the case of $\cal L$-equivalence.

For $D_1$ a frieze in $R_a$ and $D_2$ a frieze in $R_b$, we say that $D_1$ is
$\cal K$-{\it equivalent} to $D_2$, and write $D_1\cong_{\cal K} D_2$, iff
there is a homeomorphism
$h:D_1\str D_2$ such that for every $i\in \N^+$ we have $h(i,0)=(i,0)$ and
$h(i,a)=(i,b)$. It is clear that this defines an equivalence relation
on friezes,
which is wider than $\cal L$-equivalence: namely, if $D_1\cong_{\cal L} D_2$, then
$D_1\cong_{\cal K} D_2$, but the converse need not hold. If $D_1$ and $D_2$
are without circular components, then $D_1\cong_{\cal L} D_2$ iff $D_1\cong_{\cal K} D_2$.
The relation of $\cal K$-equivalence takes account only of the number
of circular components,
whereas $\cal L$-equivalence takes also account of whether circular components
are one in another, and, in general, in which region of the diagram they
are located.

If $i$ stands for $(i,a)$ and $-i$ stands for $(i,0)$, we may identify the end
points of each thread in a frieze in $R_a$ by a pair of
integers in $\Z\!-\!\{0\}$. For $M$ an ordered set and for $a,b\in M$ such that
$a<b$, let a {\it segment} $[a,b]$ in $M$
be $\{z\in M\mid a\leq z\leq b\}$. The numbers $a$ and $b$ are the end points
of $[a,b]$. We say that $[a,b]$ {\it encloses} $[c,d]$ iff $a<c$ and $d<b$.
A set of segments is {\it nonoverlapping} iff every two distinct segments
in it are either disjoint or one of these segments encloses the other.

We may then establish a one-to-one correspondence between the set $\Theta$ of threads
of a frieze and a set $S_{\Theta}$ of nonoverlapping segments in $\Z\!-\!\{0\}$.
Every element of $\Z -\!\{0\}$ is an end point of a segment in $S_{\Theta}$.
Since enclosure is irreflexive and transitive, $S_{\Theta}$ is partially ordered
by enclosure. This is a tree-like ordering without root, with a finite number of
branching nodes.
For example, in the frieze

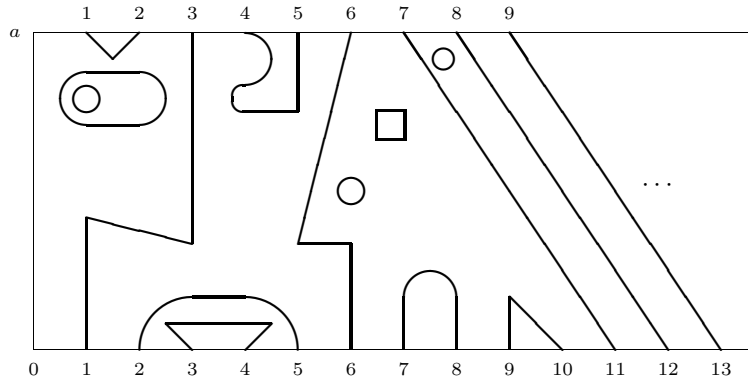
\begin{center}
\begin{picture}(270,160)
{\linethickness{0.05pt}
\put(0,20){\line(1,0){270}}
\put(0,20){\line(0,1){120}}
\put(0,140){\line(1,0){270}}}

\put(-5,140){\makebox(0,0)[r]{\scriptsize $a$}}
\put(0,15){\makebox(0,0)[t]{\scriptsize $0$}}
\put(20,15){\makebox(0,0)[t]{\scriptsize $1$}}
\put(40,15){\makebox(0,0)[t]{\scriptsize $2$}}
\put(60,15){\makebox(0,0)[t]{\scriptsize $3$}}
\put(80,15){\makebox(0,0)[t]{\scriptsize $4$}}
\put(100,15){\makebox(0,0)[t]{\scriptsize $5$}}
\put(120,15){\makebox(0,0)[t]{\scriptsize $6$}}
\put(140,15){\makebox(0,0)[t]{\scriptsize $7$}}
\put(160,15){\makebox(0,0)[t]{\scriptsize $8$}}
\put(180,15){\makebox(0,0)[t]{\scriptsize $9$}}
\put(200,15){\makebox(0,0)[t]{\scriptsize $10$}}
\put(220,15){\makebox(0,0)[t]{\scriptsize $11$}}
\put(240,15){\makebox(0,0)[t]{\scriptsize $12$}}
\put(260,15){\makebox(0,0)[t]{\scriptsize $13$}}

\put(20,145){\makebox(0,0)[b]{\scriptsize $1$}}
\put(40,145){\makebox(0,0)[b]{\scriptsize $2$}}
\put(60,145){\makebox(0,0)[b]{\scriptsize $3$}}
\put(80,145){\makebox(0,0)[b]{\scriptsize $4$}}
\put(100,145){\makebox(0,0)[b]{\scriptsize $5$}}
\put(120,145){\makebox(0,0)[b]{\scriptsize $6$}}
\put(140,145){\makebox(0,0)[b]{\scriptsize $7$}}
\put(160,145){\makebox(0,0)[b]{\scriptsize $8$}}
\put(180,145){\makebox(0,0)[b]{\scriptsize $9$}}

\thicklines
\put(20,20){\line(0,1){50}}
\put(20,70){\line(4,-1){40}}
\put(60,60){\line(0,1){80}}

\put(20,140){\line(1,-1){10}}
\put(40,140){\line(-1,-1){10}}

\put(30,115){\oval(40,20)}

\put(20,115){\circle{10}}

\put(70,20){\oval(60,40)[t]}

\put(60,20){\line(-1,1){10}}
\put(50,30){\line(1,0){40}}
\put(80,20){\line(1,1){10}}

\put(80,130){\oval(20,20)[r]}
\put(80,115){\oval(10,10)[l]}
\put(80,110){\line(1,0){20}}
\put(100,140){\line(0,-1){30}}

\put(120,20){\line(0,1){40}}
\put(120,60){\line(-1,0){20}}
\put(100,60){\line(1,4){20}}

\put(120,80){\circle{10}}

\put(130,100){\framebox(10,10){}}

\put(140,20){\line(0,1){20}}
\put(160,20){\line(0,1){20}}
\put(150,40){\oval(20,20)[t]}

\put(180,20){\line(0,1){20}}
\put(200,20){\line(-1,1){20}}

\put(220,20){\line(-2,3){80}}
\put(240,20){\line(-2,3){80}}
\put(260,20){\line(-2,3){80}}

\put(155,130){\circle{8}}

\put(230,80){$\cdots$}

\end{picture}
\end{center}

\noindent the set $\Theta$ of threads corresponds to the
following tree in $S_{\Theta}$:
{\scriptsize
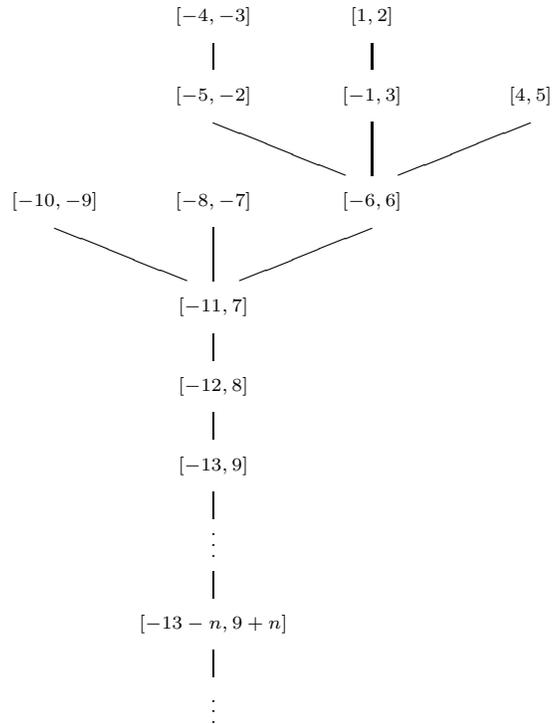
\begin{center}
\begin{picture}(180,290)

\put(60,10){\makebox(0,0){$\vdots$}}

\put(60,20){\line(0,1){10}}

\put(60,40){\makebox(0,0){$[-13-n,9+n]$}}

\put(60,50){\line(0,1){10}}

\put(60,73){\makebox(0,0){$\vdots$}}

\put(60,80){\line(0,1){10}}

\put(60,100){\makebox(0,0){$[-13,9]$}}

\put(60,110){\line(0,1){10}}

\put(60,130){\makebox(0,0){$[-12,8]$}}

\put(60,140){\line(0,1){10}}

\put(60,160){\makebox(0,0){$[-11,7]$}}

\put(60,170){\line(0,1){20}}

\put(50,170){\line(-5,2){50}}

\put(70,170){\line(5,2){50}}

\put(0,200){\makebox(0,0){$[-10,-9]$}}

\put(60,200){\makebox(0,0){$[-8,-7]$}}

\put(120,200){\makebox(0,0){$[-6,6]$}}

\put(120,210){\line(0,1){20}}

\put(110,210){\line(-5,2){50}}

\put(130,210){\line(5,2){50}}

\put(60,240){\makebox(0,0){$[-5,-2]$}}

\put(120,240){\makebox(0,0){$[-1,3]$}}

\put(180,240){\makebox(0,0){$[4,5]$}}

\put(60,250){\line(0,1){10}}

\put(120,250){\line(0,1){10}}

\put(60,270){\makebox(0,0){$[-4,-3]$}}

\put(120,270){\makebox(0,0){$[1,2]$}}

\end{picture}
\end{center}}
The branching points of this tree are $[-6,6]$ and $[-11,7]$. This tree-like ordering
of $S_{\Theta}$ induces an isomorphic ordering of $\Theta$.

If from a frieze $D$ in $R_a$ we omit all the threads, we obtain a
disjoint family
of connected sets in $R_a$, which are called the {\it regions}
of $D$. Every circular component of $D$ is included in a unique region of $D$. The closure
of a region of $D$ has a border that includes a nonempty set of threads.
In the
tree-like ordering, this set must have a lowest thread, and all the
other threads in the set, if any, are its immediate successors. Every
thread is the lowest thread for some region.
In our example, in the region in which one finds as circular components
a circle and a square,
the lowest thread is the one corresponding to $[-11,7]$, and its
immediate successors correspond to $[-10,-9]$, $[-8,-7]$ and $[-6,6]$. Assigning
to every region of a frieze the corresponding lowest thread in the border
establishes a one-to-one correspondence between regions and threads.

The collection (possibly empty) of circular components in a
single region of a frieze
corresponds to a circular form (see Section 3),
which can then be coded by an ordinal in $\varepsilon_0$. In every
frieze we can assign to every thread the ordinal that corresponds to the
collection of circular components in the region for which this is the lowest
thread.
This describes all the circular components of a frieze. (In an $\omega$-diagram
that is not a frieze it is possible that one collection of circular components, which is
in a region without lowest thread, is not covered.)

Then it is easy to establish the following.

\vspace{.2cm}

\noindent {\sc Remark 1$\cal L$.}\quad {\it The friezes} $D_1$ {\it and} $D_2$ {\it are}
$\cal L$-{\it equivalent iff}

\vspace{.1cm}

\noindent $(i)$ {\it the end points of the threads in} $D_1$
{\it are identified with the same pairs of
integers as the end points of the threads in} $D_2$,

\vspace{.1cm}

\noindent $(ii)$ {\it the same ordinals in} $\varepsilon_0$
{\it are assigned to the threads of} $D_1$ {\it and} $D_2$
{\it that are identified with the same pairs of integers.}

\vspace{.2cm}

\noindent This means that the $\cal L$-equivalence class of
a frieze may be identified with a function
$f:S_{\Theta}\str\varepsilon_0$, where the domain $S_{\Theta}$
of $f$ is a set of nonoverlapping segments in $\Z\!-\!\{0\}$.

\vspace{.2cm}

\noindent {\sc Remark 1$\cal K$.}\quad {\it The friezes} $D_1$ {\it and} $D_2$ {\it are}
$\cal K$-{\it equivalent iff}

\vspace{.1cm}

\noindent $(i)$ {\it the end points of the threads in} $D_1$
{\it are identified with the same pairs of integers as the end points
of the threads in} $D_2$,

\vspace{.1cm}

\noindent $(ii)$ $D_1$ {\it and} $D_2$ {\it have the same number of circular
components.}

\vspace{.2cm}

\noindent This means that the $\cal K$-equivalence class of a frieze may be identified with
a pair $(S_{\Theta},l)$ where $S_{\Theta}$ is a set of nonoverlapping
segments in $\Z\!-\!\{0\}$, and $l$ is a natural number, which is the
number of circular components.

The set of $\cal L$-equivalence classes of friezes is endowed with the structure of
a monoid in the following manner. Let the {\it unit frieze} $I$ be
$\{(i,y)\mid i\in\N^+\;{\mbox{\it and }} y\in[0,1]\}$ in $R_1$. So
$I$ has no circular components and all of its threads are vertical threads.
We draw $I$ as follows:
{\scriptsize
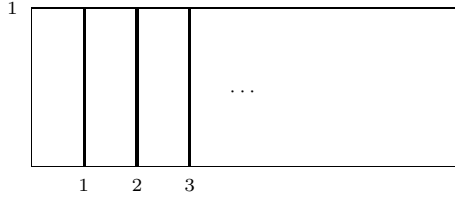
\begin{center}
\begin{picture}(160,100)
{\linethickness{0.05pt}
\put(0,20){\line(1,0){160}}
\put(0,20){\line(0,1){60}}
\put(0,80){\line(1,0){160}}}

\put(-5,80){\makebox(0,0)[r]{$1$}}
\put(20,15){\makebox(0,0)[t]{$1$}}
\put(40,15){\makebox(0,0)[t]{$2$}}
\put(60,15){\makebox(0,0)[t]{$3$}}

\thicklines
\put(20,20){\line(0,1){60}}
\put(40,20){\line(0,1){60}}
\put(60,20){\line(0,1){60}}

\put(80,50){\makebox(0,0)[t]{$\cdots$}}

\end{picture}
\end{center}}

For two friezes $D_1$ in $R_a$ and $D_2$ in $R_b$ let the {\it composition} of
$D_1$ and $D_2$ be defined as follows:
\[
D_2\cirk D_1=\{(x,y+b)\mid(x,y)\in D_1\}\cup D_2.
\]
It is easy to see that $D_2\cirk D_1$ is a frieze in $R_{a+b}$.

For $1\leq i \leq 4$, let $D_i$ be a frieze in $R_{a_i}$ and suppose
$D_1\cong_{\cal L}D_3$ with the homeomorphism $h_1:R_{a_1}\str R_{a_3}$
and
$D_2\cong_{\cal L}D_4$ with the homeomorphism $h_2:R_{a_2}\str R_{a_4}$.
Then $D_2\cirk D_1\cong_{\cal L}D_4\cirk D_3$ with the homeomorphism
$h:R_{a_1+a_2}\str R_{a_2+a_4}$ defined as follows. For $p^1$ the first and $p^2$
the second projection, let
\[h(x,y)=\left\{ \begin{array}{ll}
(p^1(h_1(x,y-a_2)),p^2(h_1(x,y-a_2))+a_4) & {\mbox{\rm if }} y>a_2
\\
h_2(x,y) & {\mbox{\rm if }} y\leq a_2
\end{array}
\right.
\]
So the composition $\cirk$ defines an operation on $\cal L$-equivalence
classes of friezes.

We can then establish that
\[
\begin{array}{ll}
(1) & I\cirk D\cong_{\cal L} D, \quad D\cirk I\cong_{\cal L} D,
\\
(2) & D_3\cirk(D_2\cirk D_1)\cong_{\cal L} (D_3\cirk D_2)\cirk D_1.
\end{array}
\]
The equivalences of $(1)$ follow from the fact that the threads of $I\cirk D$,
$D\cirk I$ and $D$ are identified with the same pairs of integers, because
all the threads of $I$ are vertical transversal threads, and from the fact
that $I$ has no circular component. Then we apply Remark 1$\cal L$. For
the equivalence $(2)$, it is clear that $D_3\cirk(D_2\cirk D_1)$ is
actually identical to $(D_3\cirk D_2)\cirk D_1$.
So the set of $\cal L$-equivalence classes of friezes has the structure of a monoid,
and the monoid structure of the set of $\cal K$-equivalence classes of friezes
is defined quite analogously. We will show for these monoids that they
are isomorphic to \Lo\ and \Ko\ respectively.

\section{Generating friezes}

For $k\in\N^+$ let the {\it cup frieze} ${\mbox{\rm V}}_k$ be the frieze in
$R_1$
without circular components, with a single semicircular cup with the end points
$(k,1)$ and $(k+1,1)$; all the other threads are straight line segments
connecting $(i,0)$ and $(i,1)$ for $i<k$ and $(i,0)$ and $(i+2,1)$ for
$i\geq k$. This frieze looks as follows:

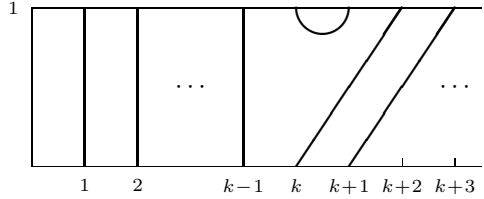
\begin{center}
\begin{picture}(170,100)
{\linethickness{0.05pt}
\put(0,20){\line(1,0){170}}
\put(0,20){\line(0,1){60}}
\put(0,80){\line(1,0){170}}
\put(140,20){\line(0,1){3}}
\put(160,20){\line(0,1){3}}}

\put(-5,80){\makebox(0,0)[r]{\scriptsize $1$}}
\put(20,15){\makebox(0,0)[t]{\scriptsize $1$}}
\put(40,15){\makebox(0,0)[t]{\scriptsize $2$}}
\put(80,15){\makebox(0,0)[t]{\scriptsize $k\!-\!1$}}
\put(100,15){\makebox(0,0)[t]{\scriptsize $k$}}
\put(120,15){\makebox(0,0)[t]{\scriptsize $k\!+\!1$}}
\put(140,15){\makebox(0,0)[t]{\scriptsize $k\!+\!2$}}
\put(160,15){\makebox(0,0)[t]{\scriptsize $k\!+\!3$}}

\thicklines
\put(20,20){\line(0,1){60}}
\put(40,20){\line(0,1){60}}
\put(80,20){\line(0,1){60}}

\put(100,20){\line(2,3){40}}
\put(120,20){\line(2,3){40}}

\put(60,50){\makebox(0,0){$\cdots$}}
\put(160,50){\makebox(0,0){$\cdots$}}

\put(110,80){\oval(20,20)[b]}
\end{picture}
\end{center}

For $k\in \N^+$ let the {\it cap frieze} $\Lambda_k$ be
the frieze in $R_1$
that is defined analogously to ${\mbox{\rm V}}_k$ and looks as follows:

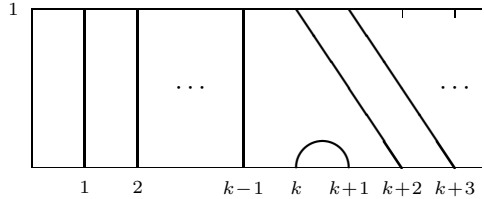
\begin{center}
\begin{picture}(170,100)
{\linethickness{0.05pt}
\put(0,20){\line(1,0){170}}
\put(0,20){\line(0,1){60}}
\put(0,80){\line(1,0){170}}
\put(140,80){\line(0,-1){3}}
\put(160,80){\line(0,-1){3}}}

\put(-5,80){\makebox(0,0)[r]{\scriptsize $1$}}
\put(20,15){\makebox(0,0)[t]{\scriptsize $1$}}
\put(40,15){\makebox(0,0)[t]{\scriptsize $2$}}
\put(80,15){\makebox(0,0)[t]{\scriptsize $k\!-\!1$}}
\put(100,15){\makebox(0,0)[t]{\scriptsize $k$}}
\put(120,15){\makebox(0,0)[t]{\scriptsize $k\!+\!1$}}
\put(140,15){\makebox(0,0)[t]{\scriptsize $k\!+\!2$}}
\put(160,15){\makebox(0,0)[t]{\scriptsize $k\!+\!3$}}

\thicklines
\put(20,20){\line(0,1){60}}
\put(40,20){\line(0,1){60}}
\put(80,20){\line(0,1){60}}

\put(100,80){\line(2,-3){40}}
\put(120,80){\line(2,-3){40}}

\put(60,50){\makebox(0,0){$\cdots$}}
\put(160,50){\makebox(0,0){$\cdots$}}

\put(110,20){\oval(20,20)[t]}
\end{picture}
\end{center}

Let a frieze without cups and caps be called a {\it circular frieze}.
Note that according to this definition the unit frieze $I$ is a
circular frieze.
For circular friezes we can prove the following lemma.

\vspace{.2cm}

\noindent {\sc Generating Circles Lemma.}\quad {\it Every circular frieze is}
$\cal L$-{\it equivalent to a frieze generated from
the unit frieze I and the
cup and cap friezes with
the operation of composition} $\cirk$.

\vspace{.1cm}

\noindent {\it Proof.}\quad If there are no circular components in our
circular frieze, then, by Remark 1$\cal L$, this frieze is
$\cal L$-equivalent to the unit frieze $I$. Suppose then
that there are circular components in our circular frieze,
and take a circular component in this frieze that is not within
another circular component. For example, let that be the right outer
circle in the following frieze

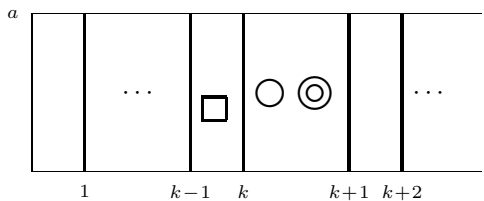
\begin{center}
\begin{picture}(170,100)
{\linethickness{0.05pt}
\put(0,20){\line(1,0){170}}
\put(0,20){\line(0,1){60}}
\put(0,80){\line(1,0){170}}}

\put(-5,80){\makebox(0,0)[r]{\scriptsize $a$}}
\put(20,15){\makebox(0,0)[t]{\scriptsize $1$}}
\put(60,15){\makebox(0,0)[t]{\scriptsize $k\!-\!1$}}
\put(80,15){\makebox(0,0)[t]{\scriptsize $k$}}
\put(120,15){\makebox(0,0)[t]{\scriptsize $k\!+\!1$}}
\put(140,15){\makebox(0,0)[t]{\scriptsize $k\!+\!2$}}

\thicklines
\put(20,20){\line(0,1){60}}
\put(60,20){\line(0,1){60}}
\put(80,20){\line(0,1){60}}
\put(120,20){\line(0,1){60}}
\put(140,20){\line(0,1){60}}

\put(65,40){\framebox(8,8){}}

\put(90,50){\circle{10}}
\put(107,50){\circle{6}}
\put(107,50){\circle{12}}

\put(40,50){\makebox(0,0){$\cdots$}}
\put(150,50){\makebox(0,0){$\cdots$}}

\end{picture}
\end{center}

\noindent We replace this by

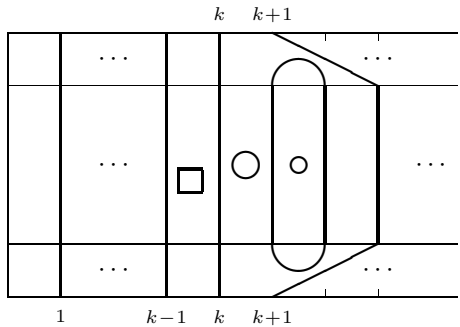
\begin{center}
\begin{picture}(170,140)
{\linethickness{0.05pt}
\put(0,20){\line(1,0){170}}
\put(0,40){\line(1,0){170}}
\put(0,20){\line(0,1){100}}
\put(0,100){\line(1,0){170}}
\put(0,120){\line(1,0){170}}
\put(120,20){\line(0,1){3}}
\put(140,20){\line(0,1){3}}
\put(120,120){\line(0,-1){3}}
\put(140,120){\line(0,-1){3}}}

\put(20,15){\makebox(0,0)[t]{\scriptsize $1$}}
\put(60,15){\makebox(0,0)[t]{\scriptsize $k\!-\!1$}}
\put(80,15){\makebox(0,0)[t]{\scriptsize $k$}}
\put(100,15){\makebox(0,0)[t]{\scriptsize $k\!+\!1$}}
\put(80,125){\makebox(0,0)[b]{\scriptsize $k$}}
\put(100,125){\makebox(0,0)[b]{\scriptsize $k\!+\!1$}}

\put(40,30){\makebox(0,0){$\cdots$}}
\put(40,70){\makebox(0,0){$\cdots$}}
\put(40,110){\makebox(0,0){$\cdots$}}

\put(140,30){\makebox(0,0){$\cdots$}}
\put(140,110){\makebox(0,0){$\cdots$}}
\put(160,70){\makebox(0,0){$\cdots$}}

\thicklines
\put(20,20){\line(0,1){100}}
\put(60,20){\line(0,1){100}}
\put(80,20){\line(0,1){100}}

\put(100,40){\line(0,1){60}}
\put(120,40){\line(0,1){60}}
\put(140,40){\line(0,1){60}}

\put(100,20){\line(2,1){40}}
\put(100,120){\line(2,-1){40}}

\put(65,60){\framebox(8,8){}}

\put(90,70){\circle{10}}
\put(110,70){\circle{6}}

\put(110,100){\oval(20,20)[t]}
\put(110,40){\oval(20,20)[b]}

\end{picture}
\end{center}

\noindent which is $\cal L$-equivalent to the original frieze.
In the frieze in the middle there
are less circular components than in the original frieze, and
the lemma follows
by induction. (By judicious choices, we can ensure that the
composition of cup
and cap friezes we obtain at the end corresponds to a term of
\Lo\ in normal form.)
\qed

\vspace{.2cm}

Note that the unit frieze $I$ is $\cal L$-equivalent to
${\mbox{\rm V}}_k \cirk \Lambda_{k+1}$
(or to ${\mbox{\rm V}}_{k+1} \cirk \Lambda_k$),
for any $k \in \N^+$, so
that, strictly
speaking, the mentioning of $I$ is superfluous in the preceding and in
the following lemma.

\vspace{.2cm}

\noindent {\sc Generating Lemma.}\quad {\it Every frieze is} $\cal L$-{\it equivalent to a frieze generated
from the unit frieze} $I$ {\it and the cup and cap friezes
with the operation of composition} $\cirk$.

\vspace{.1cm}

\noindent {\it Proof.}\quad We proceed by induction on the sum of the numbers of cups and
caps in the given frieze. The basis of the induction is covered by the
Generating Circles Lemma. If our frieze has cups, it must have a cup whose end
points are $(i,a)$ and $(i+1,a)$. If, for example, we have

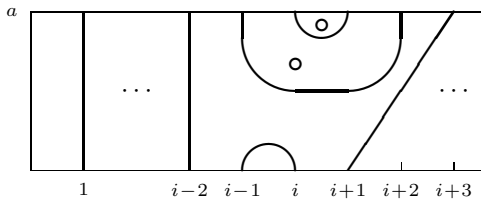
\begin{center}
\begin{picture}(170,100)
{\linethickness{0.05pt}
\put(0,20){\line(1,0){170}}
\put(0,80){\line(1,0){170}}
\put(0,20){\line(0,1){60}}
\put(160,20){\line(0,1){3}}
\put(140,20){\line(0,1){3}}}

\put(-5,80){\makebox(0,0)[r]{\scriptsize $a$}}
\put(20,15){\makebox(0,0)[t]{\scriptsize $1$}}
\put(60,15){\makebox(0,0)[t]{\scriptsize $i\!-\!2$}}
\put(80,15){\makebox(0,0)[t]{\scriptsize $i\!-\!1$}}
\put(100,15){\makebox(0,0)[t]{\scriptsize $i$}}
\put(120,15){\makebox(0,0)[t]{\scriptsize $i\!+\!1$}}
\put(140,15){\makebox(0,0)[t]{\scriptsize $i\!+\!2$}}
\put(160,15){\makebox(0,0)[t]{\scriptsize $i\!+\!3$}}

\put(40,50){\makebox(0,0){$\cdots$}}
\put(160,50){\makebox(0,0){$\cdots$}}

\thicklines
\put(20,20){\line(0,1){60}}
\put(60,20){\line(0,1){60}}

\put(120,20){\line(2,3){40}}

\put(100,60){\circle{4}}
\put(110,75){\circle{4}}

\put(110,80){\oval(20,20)[b]}
\put(90,20){\oval(20,20)[t]}
\put(110,80){\oval(60,60)[b]}

\end{picture}
\end{center}

\noindent we replace this by

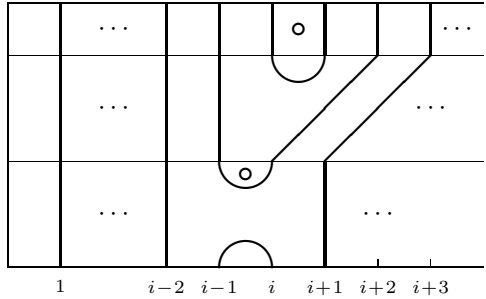
\begin{center}
\begin{picture}(180,140)
{\linethickness{0.05pt}
\put(0,20){\line(1,0){180}}
\put(0,60){\line(1,0){180}}
\put(0,100){\line(1,0){180}}
\put(0,120){\line(1,0){180}}
\put(0,20){\line(0,1){100}}
\put(160,20){\line(0,1){3}}
\put(140,20){\line(0,1){3}}}

\put(20,15){\makebox(0,0)[t]{\scriptsize $1$}}
\put(60,15){\makebox(0,0)[t]{\scriptsize $i\!-\!2$}}
\put(80,15){\makebox(0,0)[t]{\scriptsize $i\!-\!1$}}
\put(100,15){\makebox(0,0)[t]{\scriptsize $i$}}
\put(120,15){\makebox(0,0)[t]{\scriptsize $i\!+\!1$}}
\put(140,15){\makebox(0,0)[t]{\scriptsize $i\!+\!2$}}
\put(160,15){\makebox(0,0)[t]{\scriptsize $i\!+\!3$}}

\put(40,40){\makebox(0,0){$\cdots$}}
\put(40,80){\makebox(0,0){$\cdots$}}
\put(40,110){\makebox(0,0){$\cdots$}}
\put(140,40){\makebox(0,0){$\cdots$}}
\put(160,80){\makebox(0,0){$\cdots$}}
\put(170,110){\makebox(0,0){$\cdots$}}

\thicklines
\put(20,20){\line(0,1){100}}
\put(60,20){\line(0,1){100}}
\put(80,60){\line(0,1){60}}
\put(100,100){\line(0,1){20}}
\put(120,100){\line(0,1){20}}
\put(140,100){\line(0,1){20}}
\put(160,100){\line(0,1){20}}
\put(120,20){\line(0,1){40}}

\put(100,60){\line(1,1){40}}
\put(120,60){\line(1,1){40}}

\put(90,55){\circle{4}}
\put(110,110){\circle{4}}

\put(110,100){\oval(20,20)[b]}
\put(90,60){\oval(20,20)[b]}
\put(90,20){\oval(20,20)[t]}

\end{picture}
\end{center}

\noindent which is $\cal L$-equivalent to the original frieze.
In the lowest frieze there are less cups,
and the same number of caps. We apply to this frieze the induction hypothesis,
and we apply the Generating Circles Lemma to the highest frieze. We proceed
analogously with caps. (Again, by judicious choices, we can ensure that the composition
of cup and cap friezes we obtain at the end corresponds to a term of \Lo\ in normal
form.)
\qed

\vspace{.2cm}

\noindent Since $\cal L$-equivalence implies $\cal K$-equivalence, we have the
Generating Circles Lem\-ma
and the Generating Lemma also for $\cal L$-equivalence replaced by $\cal K$-equivalence.

It follows from the Generating Lemma that there are only denumerably
many $\cal L$-equivalence
classes of friezes, and the same holds a fortiori for $\cal K$-equivalence
classes. If we had allowed infinitely many cups or caps in friezes, then
we would have a continuum of different $\cal L$ or $\cal K$-equivalence
classes of friezes (which is clear from the fact that we can code 0-1
sequences with such friezes). The corresponding monoids could not then
be finitely generated, as \Lo\ and \Ko\ are. With infinitely many circular
components
we would have a continuum of different $\cal L$-equivalence classes,
but not so
for $\cal K$-equivalence classes (see Section 11).

\section{\Lo\ and \Ko\ are monoids of friezes}

Let $\cal F$ be the set of friezes. We define as follows a map $\delta$
from the terms of \Lo\ into $\cal F$:
\[
\begin{array}{l}
\delta(\dk{k})={\mbox{\rm V}}_k,
\\
\delta(\gk{k})=\Lambda_k,
\\
\delta(\mj)=I,
\\
\delta(tu)=\delta(t)\cirk\delta(u).
\end{array}
\]
We can then prove the following.

\vspace{.2cm}

\noindent {\sc Soundness Lemma.}\quad {\it If} $t=u$ {\it in} \Lo, {\it then}
$\delta(t)\cong_{\cal L}\delta(u)$.

\vspace{.1cm}

\noindent {\it Proof.}\quad We already verified in Section 6
that we have replacement of equivalents, and that the equations
$(1)$ and $(2)$
of the axiomatization of \Lo\ are satisfied for $I$ and $\cirk$. It just remains
to verify the remaining equations, which is quite straightforward.
\qed

\vspace{.2cm}

\noindent We have an analogous Soundness Lemma for \Ko\ and $\cong_{\cal K}$, involving the additional
checking of ({\it cup-cap} 4).

Let $[{\cal F}]_{\cal L}$ be the set of
$\cal L$-equivalence classes $[D]_{\cal L}= \{D': D\cong_{\cal L}D'\}$
for all friezes $D$
(and analogously with $\cal L$ replaced by $\cal K$).
This set is a monoid whose
unit is $[I]_{\cal L}$ and whose multiplication is defined
by taking that $[D_1]_{\cal L}[D_2]_{\cal L}$
is $[D_1\cirk D_2]_{\cal L}$. The Soundness Lemma guarantees that there
is a homomorphism, defined via $\delta$,
from \Lo\ to the monoid $[{\cal F}]_{\cal L}$, and
the Generating Lemma guarantees
that this homomorphism is onto. We have the same with $\cal L$
replaced by $\cal K$.
It remains to establish that these homomorphisms from \Lo\
onto $[{\cal F}]_{\cal L}$
and from \Ko\ onto $[{\cal F}]_{\cal K}$ are also one-one.

We can prove the following lemmata.

\vspace{.2cm}

\noindent {\sc Auxiliary Lemma.}\quad {\it If} $t$ {\it and} $u$ {\it are terms of}
\Lo\ {\it in normal form and}
$\delta(t)\cong_{\cal L}\delta(u)$, {\it then} $t$ {\it and} $u$ {\it are the same term.}

\vspace{.1cm}

\noindent {\it Proof.}\quad Let $t$ and $u$ be the following two terms:
\[
\begin{array}{l}
b_{j_1}^{\beta_1}\ldots b_{j_m}^{\beta_m}c_{k_1}^{\gamma_1}\ldots
c_{k_l}^{\gamma_l}a_{i_1}^{\alpha_1}\ldots a_{i_n}^{\alpha_n},
\\[.2cm]
b_{j_1'}^{\beta_1'}\ldots b_{j_{m'}'}^{\beta_{m'}'}c_{k_1'}^{\gamma_1'}\ldots
c_{k_{l'}'}^{\gamma_{l'}'}a_{i_1'}^{\alpha_1'}
\ldots a_{i_{n'}'}^{\alpha_{n'}'}.
\end{array}
\]
If $a_{i_1}^{\alpha_1}\ldots a_{i_n}^{\alpha_n}$ is different from
$a_{i_1'}^{\alpha_1'}\ldots a_{i_{n'}'}^{\alpha_{n'}'}$, then either $n<n'$
or $n'<n$ or ($n=n'$ and for some $p\in \{1,\ldots, n\}$ either $i_p\neq i_p'$
or $\alpha_p\neq \alpha_p'$). Since $i_1<\ldots< i_n$, each index $i_p$
corresponds to the left end point of a cup. So if $n<n'$ or $n'<n$ or
$i_p\neq i_p'$, then $\delta(t)$ and $\delta(u)$ don't have the same left
end points of cups, and hence, they cannot be $\cal L$-equivalent by
Remark 1$\cal L$$(i)$. If, on the other hand,
$\delta(t)$ and $\delta(u)$ have
cups identified with the same pairs of integers,
then for some $p\in \{1,\ldots, n\}$ we have
$\alpha_p\neq \alpha_p'$, and, since different ordinals are assigned
to threads of $\delta(t)$ and $\delta(u)$ identified with
the same pairs of integers,
by Remark 1$\cal L$$(ii)$, the friezes
$\delta(t)$ and $\delta(u)$ cannot be $\cal L$-equivalent.
We reason analogously with $a$ replaced by $b$ and $c$.
\qed

\vspace{.2cm}

\noindent {\sc Completeness Lemma.}\quad {\it If} $\delta(t)\cong_{\cal L}\delta(u)$,
{\it then} $t=u$ {\it in \Lo.}

\vspace{.1cm}

\noindent {\it Proof.}\quad By the Normal Form Lemma of Section 4,
for every term $t$ and every term
$u$ of \Lo\ there are terms $t'$ and $u'$ in normal form such that $t=t'$ and
$u=u'$ in \Lo.
By the Soundness Lemma, we  obtain $\delta(t)\cong_{\cal L}\delta(t')$ and
$\delta(u)\cong_{\cal L}\delta(u')$, and if $\delta(t)\cong_{\cal L}\delta(u)$,
it follows that $\delta(t')\cong_{\cal L}\delta(u')$. Then, by the Auxiliary Lemma,
the terms $t'$ and $u'$ are the same term, and hence $t=u$ in \Lo.
\qed

\vspace{.2cm}

The Auxiliary Lemma and the Completeness Lemma are easily obtained when
$\cal L$ is replaced by $\cal K$. So we may conclude that our
homomorphisms from \Lo\ onto $[{\cal F}]_{\cal L}$ and from \Ko\ onto $[{\cal F}]_{\cal K}$
are one-one, and hence \Lo\ is isomorphic to  $[{\cal F}]_{\cal L}$ and \Ko\ is
isomorphic to $[{\cal F}]_{\cal K}$.

We may also conclude that for every term $t$ of \Lo\ there is a {\it unique}
term $t'$ in normal form such that $t=t'$ in \Lo. If $t=t'$ and $t=t''$ in \Lo,
then $t'=t''$ in \Lo, and hence, by the Soundness Lemma,
$\delta(t')\cong_{\cal L}\delta(t'')$. If $t'$ and $t''$ are in normal form, by the Auxiliary
Lemma we obtain that $t'$ and $t''$ are the same term. We conclude analogously that
the $\cal K$-normal form is unique in the same sense with respect to \Ko.

\section{The monoids \Ln}

The monoid \Ln\ has for every $i\in\{1,\ldots,n-1\}$ a generator $h_i$,
called a {\it diapsis} (plural {\it diapsides}), and also for every ordinal
$\alpha\in\varepsilon_0$ and every $k\in\{1,\ldots,n+1\}$
a generator $c_k^{\alpha}$, called a $c$-{\it term}. The number $n$ here could in
principle be any natural number, but the interesting monoids \Ln\ have
$n\geq 2$. When $n$ is 0 or 1, we have no diapsides. The diapsis $h_i$
corresponds
to the term $\gk{i}\dk{i}$ of \Lo. The terms of \Ln\ are obtained
from these generators and \mj\ by closing under multiplication.

We assume the following equations for \Ln:
\[
\begin{array}{ll}
(1) & \mj t=t, \quad t\mj=t,
\\
(2) & t(uv)=(tu)v,
\\[.2cm]
(c1) & \mj=c_k^0,
\\[.1cm]
(c2) & c_k^{\alpha}c_k^{\beta}=c_k^{\alpha\sharp\beta},
\\[.1cm]
(cc) & c_k^{\alpha}c_l^{\beta}=c_l^{\beta}c_k^{\alpha},\quad
{\mbox{\rm for }} k\neq l,
\\[.2cm]
(h1) & h_i h_{j+2}=h_{j+2} h_i,\quad {\mbox{\rm for }} i \leq j,
\\
(h2) & h_i h_{i\pm 1} h_i=h_i,
\\[.2cm]
(hc1') & h_ic_k^{\alpha}=c_k^{\alpha}h_i,\quad
{\mbox{\rm for }} k\neq i+1,
\\[.1cm]
(hc2') & h_ic_{i+1}^{\alpha}h_i=c_i^{\omega^{\alpha}}h_i,
\\[.1cm]
(hc3) & c_i^{\alpha}h_i=c_{i+2}^{\alpha}h_i.
\end{array}
\]
With the help of $(c2)$ we can derive $(cc)$ for $k=l$ too.

With $h_i$ defined as $\gk{i}\dk{i}$ and $c_k^{\alpha}$ defined as in
Section 4,
we can check easily that all the equations above hold in \Lo. We can make
this checking also with friezes. So \Ln\ is a submonoid of \Lo.

An $n$-{\it frieze} is a frieze such that for every $k\geq n+1$
we have a vertical
thread identified with $[-k,k]$ and for every $k\geq n+2$ the
ordinal of circular components assigned to the thread $[-k,k]$ is $0$.
Each $n$-frieze without circular components may be conceived
up to $\cal L$-equivalence or
$\cal K$-equivalence, which here coincide, as an element of the free
(noncommutative) $o$-monoid generated by the empty set of generators
(cf. Section 3).
This is because the threads of each $n$-frieze without
circular components are
identified with a rooted subtree of $S_{\Theta}$ (see
Section 6),
whose root is $[-(n+1),n+1]$, and this rooted tree may be coded
by a parenthetical
word.

If ${\cal F}_n$ is the set of $n$-friezes, let $[{\cal F}_n]_{\cal L}$ be the set of
$\cal L$-equivalence classes of these friezes, and analogously for $[{\cal F}_n]_{\cal K}$.
The set $[{\cal F}_n]_{\cal L}$ has the structure of a monoid defined as
for $[{\cal F}]_{\cal L}$.

Then it can be shown that the monoid $[{\cal F}_n]_{\cal L}$ is
isomorphic to \Ln\
with the help of a map $\delta:\Ln\str {\cal F}_n$ that maps a diapsis
$h_k$ into the {\it diapsidal n-frieze} $H_k$, which
is the {\it n}-frieze in $R_b$, for some $b>1$, without circular components,
with a single semicircular cup with
the end points $(k,b)$ and $(k+1,b)$, and a single
semicircular cap
with the end points $(k,0)$ and $(k+1,0)$; all the other threads are
vertical
threads orthogonal to the $x$ axis. A diapsidal {\it n}-frieze $H_k$
looks as follows:

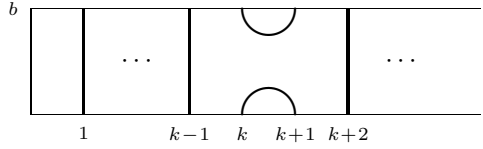
\begin{center}
\begin{picture}(170,80)
{\linethickness{0.05pt}
\put(0,20){\line(1,0){170}}
\put(0,60){\line(1,0){170}}
\put(0,20){\line(0,1){40}}}

\put(-5,60){\makebox(0,0)[r]{\scriptsize $b$}}
\put(20,15){\makebox(0,0)[t]{\scriptsize $1$}}
\put(60,15){\makebox(0,0)[t]{\scriptsize $k\!-\!1$}}
\put(80,15){\makebox(0,0)[t]{\scriptsize $k$}}
\put(100,15){\makebox(0,0)[t]{\scriptsize $k\!+\!1$}}
\put(120,15){\makebox(0,0)[t]{\scriptsize $k\!+\!2$}}

\put(40,40){\makebox(0,0){$\cdots$}}
\put(140,40){\makebox(0,0){$\cdots$}}

\thicklines
\put(20,20){\line(0,1){40}}
\put(60,20){\line(0,1){40}}
\put(120,20){\line(0,1){40}}

\put(90,60){\oval(20,20)[b]}
\put(90,20){\oval(20,20)[t]}

\end{picture}
\end{center}

\noindent The $c$-term $c_k^{\alpha}$ is mapped by $\delta$ into the frieze

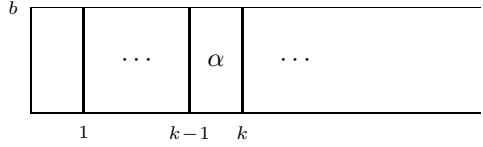
\begin{center}
\begin{picture}(170,80)
{\linethickness{0.05pt}
\put(0,20){\line(1,0){170}}
\put(0,60){\line(1,0){170}}
\put(0,20){\line(0,1){40}}}

\put(-5,60){\makebox(0,0)[r]{\scriptsize $b$}}
\put(20,15){\makebox(0,0)[t]{\scriptsize $1$}}
\put(60,15){\makebox(0,0)[t]{\scriptsize $k\!-\!1$}}
\put(80,15){\makebox(0,0)[t]{\scriptsize $k$}}
\put(70,40){\makebox(0,0){$\alpha$}}

\put(40,40){\makebox(0,0){$\cdots$}}
\put(100,40){\makebox(0,0){$\cdots$}}

\thicklines
\put(20,20){\line(0,1){40}}
\put(60,20){\line(0,1){40}}
\put(80,20){\line(0,1){40}}

\end{picture}
\end{center}

\noindent where $\alpha$ stands for an arbitrary
circular form corresponding to $\alpha$.
We also have $\delta(\mj)=I$ and $\delta(tu)=\delta(t)\cirk\delta(u)$,
as before. It is clear that the unit frieze $I$ is an
{\it n}-frieze for every $n \in \N$, and that the composition
of two {\it n}-friezes is an {\it n}-frieze.

We will not go into the details of the proof that we have an isomorphism
here, because we
don't have much use for \Ln\ in this work. A great part of this proof is analogous
to what we had for \Lo, or to what we will have for \Kn\ in the next section.
The essential part of the proof is the definition of unique normal form for
elements of \Ln. Here is how such a normal form would look like.

For $1\leq j\leq i\leq n-1$ and $\alpha,\beta\in\varepsilon_0$,
let the {\it block} $h_{[i,j]}^{\alpha,\beta}$ be defined as
\[
c_{i+1}^{\alpha}h_ih_{i-1}\ldots h_{j+1}h_jc_{j+1}^{\beta}.
\]
A term of \Ln\ in {\it normal form} will be \mj, or it looks as follows:
\[
c_{k_1}^{\gamma_1}\ldots c_{k_l}^{\gamma_l}h_{[b_1,c_1]}^{\alpha_1,\beta_1}
\ldots h_{[b_n,c_n]}^{\alpha_n,\beta_n},
\]
where $n,l\geq 0$, $k_1<\ldots<k_l$, $b_1<\ldots<b_n$, and
$c_1<\ldots<c_n$. All the $c$-terms on the left-hand side are
such that they could
be permuted with all the blocks, and pass to the right-hand side;
i.e. they would not
be ``captured'' by a block. We must also make a choice for the
indices $k_p$ of
these $c$-terms to ensure uniqueness, and $\gamma_p$ should not be 0.

\section{The monoids \Kn}

One way to define the monoid \Kn\ is to have the same generators
as for \Ln, and
the following equations, which we add to those of \Ln:
\[
\begin{array}{l}
c_k^{\omega^{\alpha}}=c_k^{\alpha+1},
\\[.2cm]
c_k^{\alpha}=c_{k+1}^{\alpha}.
\end{array}
\]
The first equation has the effect of collapsing the ordinals in
$\varepsilon_0$ into natural numbers (as the equation ({\it solid}) of
Section 3),
while the second equation has the effect of making superfluous the lower index
of $c$-terms.

An alternative, and simpler, axiomatization of \Kn\ is obtained as follows.
The monoid \Kn\ has for
every $i \in \{1,\ldots, n-1\}$ a generator $h_i$,
called again a {\it diapsis}, and also the generator $c$, called the
{\it circle}.
The terms of \Kn\
are obtained from these generators and \mj\ by closing under multiplication.
We assume the following equations for \Kn:
\[
\begin{array}{ll}
(1) & \mj t=t, \quad t \mj = t,
\\
(2) & t(uv)=(tu)v,
\\[0.1cm]
(h1) & h_i h_{j+2}=h_{j+2} h_i,\quad {\mbox{\rm for }} i \leq j,
\\
(h2) & h_i h_{i\pm 1}h_i =  h_i,
\\[0.1cm]
(hc1) & h_i c = c h_i,
\\
(hc2) & h_i h_i = c h_i.
\end{array}
\]

The equations $(h1)$, $(h2)$ and $(hc2)$,
which may be derived from Jones' paper \cite{J83} (p. 13),
and which appear in the form above in many works of Kauffman (see
\cite{KL}, \cite{KA97}, Section 6, and references therein), are usually
tied to the presentation of Temperley-Lieb algebras. They may,
however, be found in Brauer algebras too (see \cite{W88}, p. 180-181).

With $h_i$ defined as $\gk{i}\dk{i}$ and $c$ defined as $\dk{i}\gk{i}$ we can
check easily that \Kn\ is a submonoid of \Ko.

For $1\leq j \leq i \leq n-1$, let the {\it block} $h_{[i,j]}$ be defined as
$h_i h_{i-1}  \ldots  h_{j+1} h_j$. The block $h_{[i,i]}$, which is
defined as
$h_i$, will be called {\it singular}.
(One could conceive \dk{i}\ as the infinite
block $\ldots h_{i+2}h_{i+1}h_i$, whereas \gk{i}\ would be
$h_i h_{i+1}h_{i+2}\ldots$)
Let $c^1$ be $c$, and let $c^{l+1}$ be $c^lc$.

A term is in {\it Jones normal form} iff it is either of the form
$c^l h_{[b_1,a_1]} \ldots h_{[b_k,a_k]}$ for $l,k \geq 0$, $l+k\geq 1$,
$a_1< \ldots <a_k$ and $b_1< \ldots <b_k$,
or it is the term {\bf 1}
(see \cite{J83}, \S 4.1.4, p. 14).
As before, if $l=0$, then $c^l$ is the empty sequence, and if $k=0$, then
$h_{[b_1,a_1]} \ldots h_{[b_k,a_k]}$ is empty.

That every term of ${\cal K}_{n}$ is equal to a term in Jones normal form
will be demonstrated with the help of an alternative formulation of
\Kn,
called the {\it block formulation}, which is obtained as follows. Besides
the circle $c$, we take as generators the blocks $h_{[i,j]}$
instead of the diapsides, we generate terms with these generators,
{\bf 1} and multiplication, and to the equations
$(1)$ and $(2)$ we add the equations
\[
\begin{array}{ll}
(h {\rm I}) & \mbox{for }j \geq k+2,\quad
h_{[i,j]}h_{[k,l]} = h_{[k,l]}h_{[i,j]},
\\[0.2cm]
(h {\rm II}) & \mbox{for } i \geq l \; \mbox{and }
|k-j|=1, \quad h_{[i,j]}h_{[k,l]} = h_{[i,l]},
\\[0.2cm]
(hc {\rm I}) &  h_{[i,j]} c = c h_{[i,j]},
\\[0.2cm]
(hc {\rm II}) & h_{[i,j]}h_{[j,l]}= c h_{[i,l]}.
\end{array}
\]

We verify first that with $h_i$ defined as the singular block $h_{[i,i]}$
the equations
$(h1),(h2), (hc1)$ and $(hc2)$ are instances of the new equations of
the block formulation. The equation
$(h1)$ is $(h \rm I)$ for $i=j$ and $k=l$, the equation $(h2)$
is $(h \rm II)$ for
$i=j=l$ and $k=i+1$, or $i=k=l$ and $j=i-1$, the equation $(hc1)$
is $(hc \rm I)$
for $i=j$, and the equation $(hc2)$ is $(hc \rm II)$ for $i=j=l$.
We also have to verify
that in the new axiomatization we can deduce the definition of $h_{[i,j]}$
with diapsides replaced by singular blocks; namely, we have to verify
$$h_{[i,j]}=h_{[i,i]}h_{[i-1,i-1]} \ldots h_{[j+1,j+1]}h_{[j,j]},$$
which readily follows from $(h \rm II)$ for $j=k+1$. To
finish showing that the block formulation of ${\cal K}_{n}$
is equivalent to
the old formulation, we have to verify that with blocks
defined via diapsides we can deduce
$(h {\rm I}), (h {\rm II}), (hc {\rm I})$ and $(hc \rm II)$
from the old equations,
which is a straightforward exercise.

We can deduce the following equations in ${\cal K}_{n}$ for
$j+2 \le k$:
\begin{tabbing}
\hspace{3em} \= $(h \rm III.1)$ \hspace{2em} \=
$h_{[i,j]}h_{[k,l]} = h_{[k-2,l]}h_{[i,j+2]} \; \;
\mbox {if} \; \;
i \ge k \; \; \mbox {and} \; \; j \ge l,$
\\[0.2cm]
\> $(h \rm III.2)$ \> $h_{[i,j]}h_{[k,l]} = h_{[i,l]}h_{[k,j+2]}
\; \; \; \; \; \; \mbox {if} \; \;
i < k \; \; \mbox {and} \; \; j \ge l,$
\\[0.2cm]
\> $(h \rm III.3)$ \> $h_{[i,j]}h_{[k,l]} = h_{[k-2,j]}h_{[i,l]} \;
\; \; \; \; \; \mbox {if} \; \;
i \ge k \; \; \mbox {and} \; \; j < l,$
\end{tabbing}

\noindent which is also pretty straightforward. Then we
prove the following lemma. (A lemma with the same content is established
in a different
manner in \cite{J83}, pp. 13-14, and \cite{GHJ}, pp. 87-89.)

\vspace{.2cm}

\noindent {\sc Normal Form Lemma.}\quad {\it Every term of} ${\cal K}_{n}$
{\it is equal in} ${\cal K}_{n}$ {\it to a term in Jones normal form.}

\vspace{.1cm}

\noindent {\it Proof.} We will give a reduction procedure that transforms
every term into a term in Jones normal form, every reduction step
being justified by an equation of ${\cal K}_{n}$. (In logical jargon, we
establish that this procedure is strongly normalizing---namely, that
any sequence of reduction steps terminates in a term in normal form.)

Take a term in the block formulation of ${\cal K}_{n}$, and let
subterms of this term of the forms
\[
\begin{array}{l}
h_{[i,j]} h_{[k,l]}, \; \; \mbox {for} \; \; i \ge k \; \; \mbox {or} \; \;
j \ge l,
\\[0.1cm] h_{[i,j]}c,
\\[0.1cm] {\bf 1}t, t{\bf 1}
\end{array}
\]
be called {\it redexes}. A reduction of the first sort
consists in replacing
a redex of the first form by the corresponding term on the right-hand
side of one
of the equations $(h {\rm I}),(h {\rm II}),(hc \rm II),
(h {\rm III.1}),(h \rm III.2)$
and
$(h \rm III.3)$. (Note that the terms on the left-hand sides of these
equations cover
all possible redexes of the first form, and the conditions of these
equations
exclude each other.) A reduction of the second sort consists in replacing
a redex of the second form by
the right-hand side of $(hc \rm I)$, and, finally, a reduction of the
third sort consists in replacing
a redex of one the forms in the third line by $t$, according to the
equations $(1)$.

Let the {\it weight} of a block $h_{[i,j]}$ be $i-j+2$. For any
subterm  $h_{[i,j]}$ of
a term $t$ in the block formulation of ${\cal K}_{n}$, let $\rho (h_{[i,j]})$
be the number of subterms $h_{[k,l]}$ of $t$ on the right-hand side
of $h_{[i,j]}$ such that $i \ge k$ or $j \ge l$. The subterms $h_{[k,l]}$
are not necessarily immediately on the right-hand side of $h_{[i,j]}$ as in
redexes of the first form:
they may also be separated by other terms. For any subterm $c$ of a term $t$
in the block formulation ${\cal K}_{n}$, let $\tau (c)$ be the number of
blocks on the left-hand side of this $c$.

The {\it complexity measure} of a term $t$ in the block formulation is
$\mu (t)= (n_1,n_2)$ where $n_1 \ge 0$ is the
sum of the weights of all the blocks in $t$,
and $n_2 \ge 0$
is the sum of all the numbers $\rho (h_{[i,j]})$ for all blocks
$h_{[i,j]}$ in $t$ plus the sum of all the
numbers $\tau (c)$ for all circles $c$ in $t$ and plus the number of
occurrences of {\bf 1} in $t$. The ordered pairs $(n_1,n_2)$
are well-ordered lexicographically.

Then we check that if $t'$ is obtained from $t$ by a reduction,
then $\mu (t')$
is strictly smaller than $\mu (t)$. With reductions of the first sort
we have that
if they are based on $(h \rm I)$, then $n_2$ diminishes while
$n_1$ doesn't change, and if they are based on
the remaining equations, then $n_1$ diminishes.
With reductions of the
second and third sort, $n_2$ diminishes while $n_1$ doesn't change.

So, by induction on the complexity measure, we obtain that every term
is equal to a term without redexes, and it is easy to see that a term
is without redexes iff it is in Jones normal form.
\qed

\vspace{.2cm}

Note that for a term  $c^l h_{[b_1,a_1]} \ldots h_{[b_k,a_k]}$
in Jones normal form the number $a_i$ is strictly smaller than all indices
of diapsides on the right-hand side of $h_{a_i}$, and $b_i$
is strictly greater than all indices
of diapsides on the left-hand side of $h_{b_i}$. The following remark is an
immediate consequence of that.

\vspace{.2cm}

\noindent {\sc Remark 2.}\quad
{\it If in a term in Jones normal form a diapsis} $h_i$ {\it occurs more
than once, then in between any two occurrences of} $h_i$
{\it we have an occurrence of $h_{i+1}$ and an occurrence of} $h_{i-1}$.

\vspace{.2cm}

A normal form dual to Jones' is obtained with blocks $h_{[i,j]}$ where
$i \le j$,
which are defined as $h_i h_{i+1} \ldots h_{j-1} h_j$. Then in
$c^l h_{[a_1,b_1]} \ldots h_{[a_k,b_k]}$
we require that $a_1> \ldots >a_k$ and $b_1> \ldots >b_k$.
The length of this new normal form will be the same as the length of Jones'.
As a matter of fact, we could take as a term in normal form many other terms
of the same reduced length as terms in the Jones normal form. For all these
alternative normal forms we can establish the property of Remark 2.

Consider the monoid $[\Fn]_{\cal K}$, i.e. the monoid of $n$-friezes
modulo $\cal K$-equiva\-lence,
which we introduced in the previous section. We will show that \Kn\ is isomorphic
to $[\Fn]_{\cal K}$.
It is clear that the following holds.

\vspace{.2cm}

\noindent {\sc Remark 3.}\quad {\it In every} $n$-{\it frieze the number of
cups is equal to the number of caps.}

\vspace{.2cm}

If the end points of a thread of an {\it n}-frieze are $(i,x)$
and $(j,y)$, where
$x,y \in \{0,a \}$, let us say that this thread {\it covers} a pair of
natural
numbers $(m,l)$, where $1 \le m<l \le n$, iff $min \{ i,j \} \le m$ and
$l \le max \{ i,j \}$.
Then we can establish the following.

\vspace{.2cm}

\noindent {\sc Remark 4.}\quad
{\it In every} $n$-{\it frieze, every pair} $(m,m+1)$, {\it where}
$1 \le m<n$,
{\it is covered by an even number of threads.}

\vspace{.1cm}

\noindent {\it Proof.}\quad For an {\it n}-frieze $D$ the cardinality of the set
$P= \{ (i,x) \in t(D) \cup b(D)\;|\; i \le m \}$ is $2m$.
Every thread of $D$ that covers $(m,m+1)$ has a single end point in $P$,
and other threads of $D$ have 0 or 2 end points in $P$. Since $P$ is even,
the remark follows.
\qed

\vspace{.2cm}

If the end points of a thread of a frieze are $(i,x)$
and $(j,y)$, let
the {\it span of the thread} be $|i-j|$. Let the {\it span} $\sigma (D)$
{\it of an n-frieze} $D$ be the sum of the spans of all the
threads of $D$. The span of an $n$-frieze must be finite, because
every $n$-frieze
has a finite number of nonvertical threads.
Remark 4 entails that the span of an
{\it n}-frieze is an even number greater than or equal to 0; the
span of $I$ is 0. It is clear that $\cal K$-equivalent {\it n}-friezes
have the
same span. It is also easy to see that the following holds.

\vspace{.2cm}

\noindent {\sc Remark 5.}
{\it If an n-frieze has cups, then it must have at least one cup
whose span is} $1$. {\it The same holds for caps.}

\vspace{.2cm}

\noindent It is clear that the same holds for friezes in general, and we have
used that fact in the proof of the Generating Lemma of Section 7.

We ascertained in the preceding section that
the unit frieze $I$ of Section 6
is an {\it n}-frieze for every $n\in \N$. We have also
defined there what is the diapsidal
{\it n}-frieze $H_i$ for $i \in \{1, \cdots ,n-1\}$.
The {\it circular} {\it n}-frieze $C$ is the {\it n}-frieze that
differs from the unit
frieze $I$ by having a single circular component,
which, for the sake of definiteness,
we choose to be a circle of radius $1/4$,
with centre $(1/2,1/2)$. We have also mentioned that
the composition of two {\it n}-friezes is an {\it n}-frieze.
Then we can prove the following lemma. (Different, and more sketchy, proofs
of this lemma may be found in \cite{PS}, Chapter VIII, Section 26, and
\cite{KA97}, Section 6; in \cite{BJ97}, Proposition 4.1.3, one may find a
proof of something more general, and somewhat more complicated.)

\vspace{.2cm}

\noindent {\sc Generating Lemma.}\quad
{\it Every n-frieze is} $\cal K$-{\it equivalent to an n-frieze generated
from I, C and the diapsidal n-friezes}
$H_i$, {\it for} $i \in \{1, \cdots ,n-1\}$,
{\it with the operation of composition} $\cirk$.

\vspace{.1cm}

\noindent {\it Proof.}\quad Take an arbitrary {\it n}-frieze $D$ in $R_a$, and
let $D$ have $l \ge 1$ circular components.
Let $C^1$ be $C$, and let $C^{l+1}$
be $C^l \cirk C$. It is clear that $D$ is $\cal K$-equivalent
to an {\it n}-frieze
$C^l \cirk D_1$ in  $R_{a+l}$ with $D_1$ an
{\it n}-frieze in $R_a$ without circular components.
If $l=0$, then $D_1$ is $D$.

Then we proceed by induction on  $\sigma (D_1)$. If $\sigma (D_1)=0$,
then $D_1 \cong_{\cal K} I$, and the lemma holds. Suppose
$\sigma (D_1)>0$. Then there must be a cup or a cap in $D_1$,
and by Remarks 3 and 5, there must be at least one cup of $D_1$
whose span is 1.
The end points of such cups are of the form $(i,a)$ and $(i+1,a)$.
Then select among
all these cups that one where $i$ is the greatest number; let that number
be $j$, and let us call the cup we have selected $\upsilon_j$.
By Remark 4,
there must be at least one other thread of $D_1$,
different from $\upsilon_j$,
that also covers $(j,j+1)$. We have to consider four cases, which exclude
each other:
\begin{itemize}
\item[$(1)$] {\rm $(j,j+1)$ is covered by a cup $\xi$ of $D_1$ different
from $\upsilon_j$, whose end points are $(p,a)$ and $(q,a)$ with $p<q$;}
\item[$(2.1)$] $(j,j+1)$ {\rm is not covered by a cup of $D_1$ different
from $\upsilon_j$, but it is covered by a transversal thread
$\xi$ of $D_1$ whose end points are $(p,a)$ and $(q,0)$ with $p<q$;}
\item[$(2.2)$] {\rm same as case (2.1) save that the end points of
$\xi$ are $(p,0)$ and $(q,a)$ with $p<q$;}
\item[$(3)$] $(j,j+1)$ {\rm is covered neither by a cup of $D_1$ different
from $\upsilon_j$, nor by a transversal thread
of $D_1$, but it is covered by a cap $\xi$ of $D_1$, whose end points
are $(p,0)$ and $(q,0)$ with $p<q$.}
\end{itemize}

In cases (1) and (2.1) we select among the threads $\xi$ mentioned
the one where $p$ is maximal, and in cases
(2.2) and (3) we select the $\xi$ where $p$ is minimal.
(We obtain the same result if in cases (1) and (2.2) we take $q$ minimal,
while in (2.1) and (3) we take $q$ maximal.)

We build out of the $n$-frieze $D_1$ a new $n$-frieze $D_2$ in
$R_a$ by replacing the thread $\upsilon_j$ and the selected thread
$\xi$, whose end points are $(p,x)$ and $(q,y)$, where $x,y \in \{ 0,a \}$,
with two new threads: one whose end points are $(p,x)$ and $(j,a)$,
and the other
whose end points are $(j+1,a)$ and $(q,y)$. We can easily check that
$D_1 \cong_{\cal K} D_2 \cirk H_j$. This is clear from the following picture:

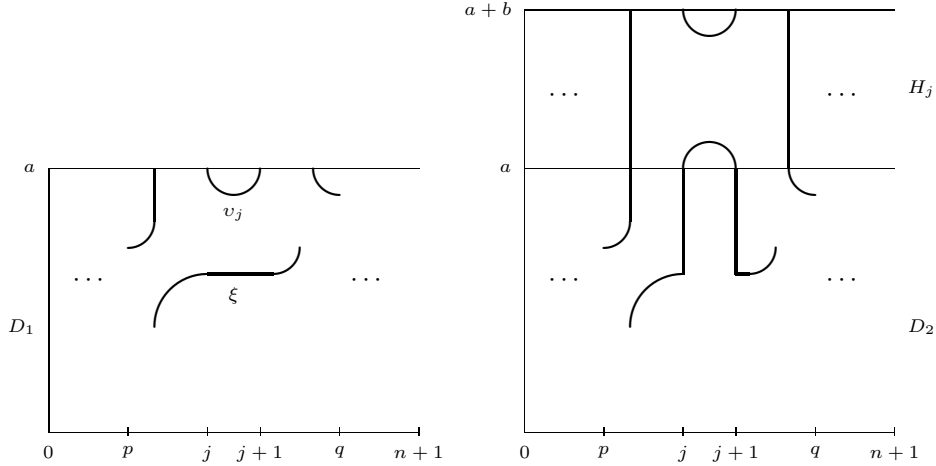
\begin{center}
\begin{picture}(320,200)
{\linethickness{0.05pt}
\put(0,20){\line(1,0){140}}
\put(0,20){\line(0,1){100}}
\put(0,120){\line(1,0){140}}

\put(180,20){\line(1,0){140}}
\put(180,20){\line(0,1){160}}
\put(180,180){\line(1,0){140}}
\put(30,19){\line(0,1){3}}
\put(60,19){\line(0,1){3}}
\put(80,19){\line(0,1){3}}
\put(110,19){\line(0,1){3}}
\put(140,19){\line(0,1){3}}
\put(210,19){\line(0,1){3}}
\put(240,19){\line(0,1){3}}
\put(260,19){\line(0,1){3}}
\put(290,19){\line(0,1){3}}
\put(320,19){\line(0,1){3}}
}

\put(-5,120){\makebox(0,0)[r]{\scriptsize $a$}}
\put(0,15){\makebox(0,0)[t]{\scriptsize $0$}}
\put(30,15){\makebox(0,0)[t]{\scriptsize $p$}}
\put(60,15){\makebox(0,0)[t]{\scriptsize $j$}}
\put(80,15){\makebox(0,0)[t]{\scriptsize $j+1$}}
\put(110,15){\makebox(0,0)[t]{\scriptsize $q$}}
\put(140,15){\makebox(0,0)[t]{\scriptsize $n+1$}}

\put(180,15){\makebox(0,0)[t]{\scriptsize $0$}}
\put(175,120){\makebox(0,0)[r]{\scriptsize $a$}}
\put(175,180){\makebox(0,0)[r]{\scriptsize $a+b$}}
\put(210,15){\makebox(0,0)[t]{\scriptsize $p$}}
\put(240,15){\makebox(0,0)[t]{\scriptsize $j$}}
\put(260,15){\makebox(0,0)[t]{\scriptsize $j+1$}}
\put(290,15){\makebox(0,0)[t]{\scriptsize $q$}}
\put(320,15){\makebox(0,0)[t]{\scriptsize $n+1$}}
\put(70,105){\makebox(0,0)[t]{\scriptsize $\upsilon_j$}}
\put(70,75){\makebox(0,0)[t]{\scriptsize $\xi$}}

\thicklines
\put(30,120){\oval(20,60)[br]}
\put(70,120){\oval(20,20)[b]}
\put(110,120){\oval(20,20)[bl]}
\put(80,60){\oval(80,40)[tl]}
\put(80,90){\oval(30,20)[br]}

\put(15,80){\makebox(0,0)[t]{$\cdots$}}
\put(120,80){\makebox(0,0)[t]{$\cdots$}}

\put(-5,60){\makebox(0,0)[r]{\scriptsize $D_1$}}

\put(210,120){\oval(20,60)[br]}
\put(290,120){\oval(20,20)[bl]}
\put(240,60){\oval(40,40)[tl]}
\put(260,90){\oval(30,20)[br]}
\put(240,80){\line(0,1){40}}
\put(260,80){\line(0,1){40}}

\put(220,120){\line(0,1){60}}
\put(280,120){\line(0,1){60}}

\put(250,120){\oval(20,20)[t]}
\put(250,180){\oval(20,20)[b]}

\put(195,80){\makebox(0,0)[t]{$\cdots$}}
\put(300,80){\makebox(0,0)[t]{$\cdots$}}
\put(195,150){\makebox(0,0)[t]{$\cdots$}}
\put(300,150){\makebox(0,0)[t]{$\cdots$}}

\put(325,60){\makebox(0,0)[l]{\scriptsize $D_2$}}
\put(325,150){\makebox(0,0)[l]{\scriptsize $H_j$}}

{\linethickness{0.05pt}
\put(180,120){\line(1,0){140}}}

\end{picture}
\end{center}

Neither of the new threads of $D_2$ that have replaced $\upsilon_j$
and $\xi$ covers $(j,j+1)$, and $\sigma (D_1)= \sigma (D_2)+2$.
So, by the induction hypothesis, $D_2$ is $\cal K$-equivalent to an
{\it n}-frieze $D_3$ generated from $I$, $C$ and $H_i$ with $\cirk$,
and since
$D_1 \cong_{\cal K} D_3 \cirk H_j$, this proves the lemma.
\qed

\vspace{.2cm}

Note that we need not require in this proof that in $\upsilon_j$
the number $j$ should be the greatest number $i$ for cups with end points
$(i,a)$ and $(i+1,a)$. The proof would go through without making this choice.
But with this choice we will end up with a composition of $n$-friezes
that
corresponds exactly to a term of \Kn\ in Jones normal form.

Let ${\cal D}_{n}$ be the set of {\it n}-friezes. We define as follows a map
$\delta$ from the terms of \Kn\ into ${\cal D}_{n}$:

\[\begin{array}{lcl}
\delta (h_i) & = & H_{i},\\[0.1cm]
\delta (c) & = & C,\\[0.1cm]
\delta ({\bf 1}) & = & I,\\[0.1cm]
\delta (tu) & = & \delta (t) \cirk \delta (u).
\end{array}\]
We can then prove the following.

\vspace{.2cm}

\noindent {\sc Soundness Lemma.}\quad
{\it If} $t=u$ {\it in}  ${\cal K}_{n}$, {\it then}
$\delta (t) \cong_{\cal K} \delta (u)$.

\vspace{.1cm}

\noindent {\it Proof.}\quad We verify first as in Section 6
that we have
replacement of equivalents, and that the equations (1) and (2)
of the axiomatization of ${\cal K}_{n}$ are satisfied for $I$ and
$\cirk$. Then it remains to verify $(h1), (h2), (hc1)$ and $(hc2)$,
which is quite straightforward.
\qed

\vspace{.2cm}

We want to show that the homomorphism
from \Kn\ to $[\Fn]_{\cal K}$ defined via $\delta$, whose existence is guaranteed
by the Soundness Lemma, is an
isomorphism. The Generating Lemma guarantees that this homomorphism
is onto, and it remains to establish that it is one-one.

Let a transversal thread in an {\it n}-frieze be called {\it falling} iff
its end points are $(i,a)$ and $(j,0)$ with $i<j$. If the end points of
a thread of an {\it n}-frieze $D$ are
$(i,a)$ and $(j,x)$ with $x \in \{ 0,a \}$ and $i<j$, then we say that
$(i,a)$ is a {\it top slope point} of $D$. If the end points of a
thread of $D$ are $(i,x)$ and $(j,0)$ with $x \in \{ 0,a \}$ and $i<j$,
then we say that $(j,0)$ is a {\it bottom slope point} of $D$. Each cup
has a single top slope point, each cap has a single bottom slope point,
and each falling transversal thread has one top slope point and one
bottom slope point. Other transversal threads have no slope points.
So we can ascertain the following.

\vspace{.2cm}

\noindent {\sc Remark 6.}\quad
{\it In every} $n$-{\it frieze the number of top slope points is equal
to the number  of bottom slope points.}

\vspace{.2cm}

Remember that by Remark 3 the number of cups is equal to the
number of caps.

Let $(a_1,a), \ldots , (a_k,a)$ be the sequence of all top slope points
of an  {\it n}-frieze $D$, ordered so that $a_1< \ldots <a_k$, and let
$(b_1+1,0), \ldots ,(b_k+1,0)$ be the sequence of all bottom slope
points of $D$,
ordered so that $b_1< \ldots <b_k$ (as we just saw with Remark 6,
these sequences must be of equal length). Then let $T_D$ be the sequence
of natural numbers $a_1, \ldots , a_k$ and $B_D$ the sequence of natural
numbers  $b_1, \ldots , b_k$.

\vspace{.2cm}

\noindent {\sc Remark 7.}\quad
{\it The sequence} $T_{\delta (h_{[i,j]})}$ {\it has a single member} $j$
{\it and the sequence} $B_{\delta (h_{[i,j]})}$ {\it has a single member} $i$.

\vspace{.2cm}

This is clear from the $n$-frieze $\delta (h_{[i,j]})$,
which is $\cal K$-equivalent to an $n$-frieze of the following form:
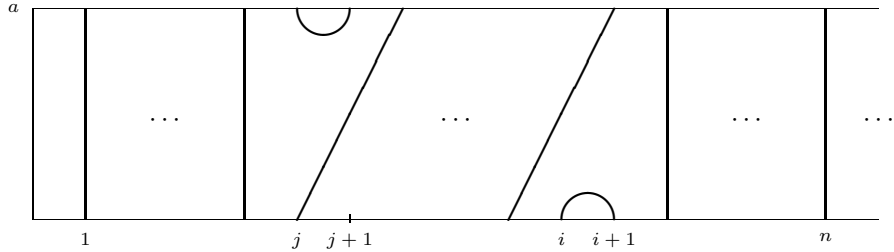
\begin{center}
\begin{picture}(320,120)
{\linethickness{0.05pt}
\put(0,20){\line(1,0){320}}
\put(0,20){\line(0,1){80}}
\put(0,100){\line(1,0){320}}}
\put(120,19){\line(0,1){3}}

\put(-5,100){\makebox(0,0)[r]{\scriptsize $a$}}
\put(20,15){\makebox(0,0)[t]{\scriptsize $1$}}
\put(100,15){\makebox(0,0)[t]{\scriptsize $j$}}
\put(120,15){\makebox(0,0)[t]{\scriptsize $j+1$}}
\put(200,15){\makebox(0,0)[t]{\scriptsize $i$}}
\put(220,15){\makebox(0,0)[t]{\scriptsize $i+1$}}
\put(300,15){\makebox(0,0)[t]{\scriptsize $n$}}

\thicklines
\put(20,20){\line(0,1){80}}
\put(80,20){\line(0,1){80}}
\put(240,20){\line(0,1){80}}
\put(300,20){\line(0,1){80}}

\put(100,20){\line(1,2){40}}
\put(180,20){\line(1,2){40}}

\put(110,100){\oval(20,20)[b]}
\put(210,20){\oval(20,20)[t]}

\put(50,60){\makebox(0,0)[t]{$\cdots$}}
\put(160,60){\makebox(0,0)[t]{$\cdots$}}
\put(270,60){\makebox(0,0)[t]{$\cdots$}}
\put(320,60){\makebox(0,0)[t]{$\cdots$}}

\end{picture}
\end{center}

\noindent provided $1<j<i<n-1$ (in other cases we simplify this
picture by omitting some
transversal threads).

\vspace{.2cm}

\noindent {\sc Remark 8.}\quad
{\it If} $D_1 \cong_{\cal K} D_2$, {\it then} $T_{D_1} = T_{D_2}$ {\it and}
$B_{D_1} = B_{D_2}.$

\vspace{.2cm}

This follows from Remark 1$\cal K$$(i)$ of Section 6.

Then we can prove the following lemmata.

\vspace{.2cm}

\noindent {\sc Key Lemma.} \quad
{\it If} $t$ {\it is the term}
$h_{[b_1,a_1]} \ldots h_{[b_k,a_k]}$ {\it with}
$a_1< \ldots <a_k$ {\it and} $b_1< \ldots <b_k$,
{\it then} $T_{\delta (t)}$ {\it is} $a_1, \ldots , a_k$ {\it and}
$B_{\delta (t)}$ {\it is} $b_1, \ldots , b_k$.

\vspace{.1cm}

\noindent {\it Proof.}\quad We proceed by induction on $k$. If $k=1$,
we use Remark 7. If $k>1$, then, by the induction hypothesis, the lemma
has been established for the term
$h_{[b_1,a_1]} \ldots h_{[b_{k-1},a_{k-1}]}$,
which we call $t'$. So $T_{\delta (t')}$ is $a_1, \ldots , a_{k-1}$.

Since in $\delta (h_{[b_k,a_k]})$ every point in the top with the first
coordinate $i<a_k$ is the end point of a vertical transversal thread,
and since $\delta (t)= \delta (t') \cirk \delta (h_{[b_k,a_k]})$, the beginning
of the sequence $T_{\delta (t)}$ must be $a_1, \ldots , a_{k-1}$. To this
sequence we have to add $a_k$ because $\delta (t)$ inherits the cup of
$\delta (h_{[b_k,a_k]})$. This shows immediately that $a_k+1$ is not in
$T_{\delta (t)}$. It remains to show that for no $i \ge a_k+2$ we can have in
$\delta (t)$ a top slope point with the first coordinate $i$.

If $i>b_k+1$, then every point in the top with the first coordinate $i$ is
the end point of a vertical transversal thread in both $\delta (h_{[b_k,a_k]})$
and $\delta (t')$. So $i$ is not in $T_{\delta (t)}$. It remains to consider
$i$ for $a_k+2 \le i \le b_k+1$. Every point in the top with
this first coordinate $i$ is the end point of a transversal thread in
$\delta (h_{[b_k,a_k]})$ whose other end point is $(i-2,0)$. If $i$ were to
be added to  $T_{\delta (t)}$, the number $i-2$ would be in $T_{\delta (t')}$,
but this contradicts the fact that $T_{\delta (t')}$ ends with $a_{k-1}$. So
$T_{\delta (t)}$ is $a_1, \ldots , a_k$.

To show that  $B_{\delta (t)}$ is $b_1, \ldots , b_k$ we reason analogously by
applying the induction hypothesis to $h_{[b_2,a_2]} \ldots h_{[b_k,a_k]}$.
\qed

\vspace{.2cm}

\noindent {\sc Auxiliary Lemma}\quad
{\it If} $t$ {\it and} $u$ {\it are terms of} ${\cal K}_{n}$
{\it in Jones normal form and}
$\delta (t) \cong_{\cal K} \delta (u)$, {\it then} $t$ {\it and} $u$
{\it are the same term.}

\vspace{.1cm}

\noindent {\it Proof.}\quad Let $t$ be $c^l h_{[b_1,a_1]} \ldots h_{[b_k,a_k]}$ and
let $u$ be $c^j h_{[d_1, c_1]} \ldots h_{[d_m, c_m]}$. If $l \neq j$, then
$\delta (t)$ is not $\cal K$-equivalent to $\delta (u)$ by
Remark 1$\cal K$$(ii)$ of Section 6, because
$\delta (t)$ and $\delta (u)$ have different numbers of circular components.
If $a_1, \ldots , a_k$ is different from  $c_1, \ldots , c_m$, or
$b_1, \ldots , b_k$
is different from $d_1, \ldots , d_m$, then $\delta (t)$ is not
$\cal K$-equivalent to
$\delta (u)$ by the Key Lemma and Remark 8.
\qed

\vspace{.2cm}

\noindent {\sc Completeness Lemma.}\quad
{\it If} $\delta (t) \cong_{\cal K} \delta (u)$, {\it then} $t=u$ {\it in}
${\cal K}_{n}$.

\vspace{.2cm}

This last lemma is proved analogously to the Completeness Lemma
of Section 8 by using the Normal Form Lemma, the Soundness Lemma and the
Auxiliary Lemma of the present section.
With this lemma we have established that \Kn\ is isomorphic
to $[\Fn]_{\cal K}$.

By reasoning as at the end of Section 8, we
can now conclude that for every term $t$ of \Kn\ there is a
unique term $t'$ in Jones normal form such that $t=t'$ in
\Kn.

\section{The monoid \Mo}

Let \Mo\ be the monoid defined as \Lo\ save that for every $k\in\N^+$ we
require also
\[
\dk{k}\gk{k}=\mj,
\]
i.e. $[k]=\mj$. It is clear that all the equations of
\Ko\ are satisfied in \Mo,
but not conversely. In \Mo\ circles are irrelevant.

The monoid \Mn\ is obtained by extending \Kn\ with $c_k^1=\mj$, or
$c=\mj$. Alternatively,
we may omit $c$-terms, or the generator $c$,
and assume only the equations $(1)$,
$(2)$, $(h1)$ and $(h2)$ of Sections 8 and 9, together with the
idempotency of $h_i$, namely, $h_ih_i=h_i$. (These axioms may be
found in \cite{J83}, p. 13.)
The monoids \Mn\ are submonoids of \Mo.

Let a \M-{\it frieze} be an $\omega$-diagram with a finite number
of cups and caps
and denumerably many circular components. (Instead
of ``denumerably many circular
components'' we could put ``$\kappa$ circular
components for a fixed infinite
cardinal $\kappa$''; for the sake of definiteness,
we chose $\kappa$ to be the
least infinite cardinal $\omega$.) We define $\cal K$-equivalence of
\M-friezes as for
friezes, and we transpose other definitions of Section 6
to \M-friezes in the same manner. It is clear that the following holds.

\vspace{.2cm}

\noindent {\sc Remark 1$\cal J$.}\quad {\it The} \M-{\it friezes} $D_1$ {\it and}
$D_2$ {\it are} $\cal K$-{\it equivalent
iff the end points of the threads in} $D_1$ {\it are identified with the same
pairs of integers as the end points of the threads in} $D_2$.

\vspace{.2cm}

\noindent So we need not pay attention any more to circular components.

The unit \M-frieze is defined as the unit frieze $I$ save that we
assume that it has denumerably many circular components, which are
located in some arbitrary regions. With composition of \M-friezes
defined as before, the set of $\cal K$-equivalence classes of
\M-friezes makes a monoid.

By adapting the argument in Sections 7 and 8, we can
show that this monoid is isomorphic to \Mo. We don't need any more the Generating
Circles Lemma, since circular \M-friezes are $\cal K$-equivalent to the unit
\M-frieze. The cup and cap friezes ${\mbox{\rm V}}_k$ and $\Lambda_k$ have
now denumerably many circular components, which are located in some arbitrary
regions.

A \M-$n$-frieze is defined as an $n$-frieze save that it has
denumerably many
circular components. Then we can show by adapting the argument in the
preceding section
that \Mn\ is isomorphic to the monoid
$[{\cal F}_{{\cal J}{\mbox {\scriptsize -}} n}]_{\cal K}$
of $\cal K$-equivalence classes of \M-$n$-friezes.

An alternative proof that the map from \Mn\ to
$[{\cal F}_{{\cal J}{\mbox {\scriptsize -}} n}]_{\cal K}$,
defined analogously to what we had
in the preceding section,
is one-one may be obtained as follows. One can establish that
the cardinality of
$[{\cal F}_{{\cal J}{\mbox {\scriptsize -}} n}]_{\cal K}$ is
the $n$-th Catalan number $(2n)!/(n!(n+1)!)$
(see the comment after the definition of $n$-friezes in Section 9;
see also \cite{KA97}, Section 6.1, and references therein).
Independently, one establishes as in \cite{J83}
(p. 14) that the number of terms of \Mn\ in Jones normal form is also the
$n$-th Catalan number. So, by the Normal Form Lemma of the preceding section, the cardinality of
\Mn\ is at most the $n$-th Catalan number. Since, by the Generating
Lemma of that section, it is known that the map above is onto,
it follows that it is
one-one. This argument is on the lines of the argument in \cite{B11}
(Note C, pp. 464-465), which establishes that the standard presentation
of symmetric groups is complete with respect to permutations. It can also
be adapted to give an alternative proof of
the Completeness Lemma of the preceding section, which is not based on
the Key Lemma and the Auxiliary Lemma of that section.

\section{The maximality of \Mo}

We will now show that \Mo\ is maximal in the following sense. Let $t$ and $u$
be terms of \Lo\ such that $t=u$ does not hold in \Mo. If \X\ is defined
as \Mo\
save that we require also $t=u$, then for
every $k\in \N^+$ we have $\gk{k}\dk{k}=\mj$ in \X. With the same
assumptions, for some $n\in \N$ we have that
\X\ is isomorphic to the monoid $\Z/n$, i.e. the additive
commutative monoid \Z\ with equality modulo $n$.

For $t$ a term of \Lo, and for $\delta(t)$ the corresponding \M-frieze,
defined analogously to what we had
in Section 8, let $\cu(t)\in\N$ be the number of cups
in $\delta(t)$, and $\ca(t)\in\N$ the number of caps in $\delta(t)$.
For $t_1$ and $t_2$ terms of \Lo, let the {\it balance} $\beta(t_1,t_2)\in\N$
of the pair $(t_1,t_2)$ be defined by
\[
\beta(t_1,t_2)=|\cu(t_1)-\cu(t_2)+\ca(t_2)-\ca(t_1)|.
\]

Let \X\ be defined as above. We will show that \X\ is isomorphic to
$\Z/\beta(t,u)$. In order to prove that we need first the following lemma.

\vspace{.2cm}

\noindent {\sc Balance Lemma.}\quad
{\it If} $t_1=t_2$ {\it holds in} \X, {\it then for some} $n\in\N$
{\it we have that}
$\beta(t_1,t_2)=n\beta(t,u)$.

\vspace{.1cm}

\noindent {\it Proof.}\quad We proceed by induction on the length of the derivation of
$t_1=t_2$ in \X. If $t_1=t_2$ holds in \Mo,
then $\beta(t_1,t_2)=0=0\cdot\beta(t,u)$, and if $t_1=t_2$ is $t=u$,
then $\beta(t_1,t_2)=\beta(t,u)$. It is easy to see that
$\beta(t_1,t_2)=\beta(t_2,t_1)$. Next, if for some $n_1,n_2\in\N$ we have
$\beta(t_1,t_2)=n_1\beta(t,u)$ and $\beta(t_2,t_3)=n_2\beta(t,u)$,
then for some $z_1,z_2 \in\Z$ such that $|z_1|=n_1$ and $|z_2|=n_2$
\[
\begin{array}{l}
\cu(t_1)-\cu(t_2)+\ca(t_2)-\ca(t_1)=z_1 \beta(t,u),
\\
\cu(t_2)-\cu(t_3)+\ca(t_3)-\ca(t_2)=z_2 \beta(t,u).
\end{array}
\]
Then
\[
\cu(t_1)-\cu(t_3)+\ca(t_3)-\ca(t_1)=(z_1+ z_2) \beta(t,u),
\]
and hence
$\beta(t_1,t_3)=|z_1+z_2| \beta(t,u)$.

We also have that $\beta(t_1,t_2)=\beta(\dk{k}t_1,\dk{k}t_2)$.
To show that, we have the following cases for $\delta(t_i)$, $i\in\{1,2\}$:
\begin{description}
\item[$(1.1)$] $(k,0)$ and $(k+1,0)$ are the end points of a single cap;
\item[$(1.2)$] $(k,0)$ and $(k+1,0)$ are the end points of two caps;
\item[$(1.3)$] one of $(k,0)$ and $(k+1,0)$ is the end point of a
cap, and the other is the end point of a transversal thread;
\item[$(2)$] $(k,0)$ and $(k+1,0)$ are the end points of two transversal threads.
\end{description}

Here we illustrate $\delta(\dk{k}t_i)$ in these various cases:


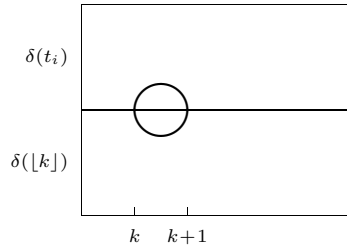
\begin{center}
\begin{picture}(300,110)
{\linethickness{0.05pt}
\put(80,20){\line(1,0){100}}
\put(80,60){\line(1,0){100}}
\put(80,100){\line(1,0){100}}
\put(80,20){\line(0,1){80}}
\put(100,20){\line(0,1){3}}
\put(120,20){\line(0,1){3}}}

\put(75,40){\makebox(0,0)[r]{\scriptsize $\delta(\dk{k})$}}
\put(75,80){\makebox(0,0)[r]{\scriptsize $\delta(t_i)$}}
\put(0,100){\makebox(0,0)[l]{(1.1)}}
\put(100,15){\makebox(0,0)[t]{\scriptsize $k$}}
\put(120,15){\makebox(0,0)[t]{\scriptsize $k\!+\!1$}}

\thicklines
\put(110,60){\oval(20,20)[b]}
\put(110,60){\oval(20,20)[t]}

\end{picture}
\end{center}


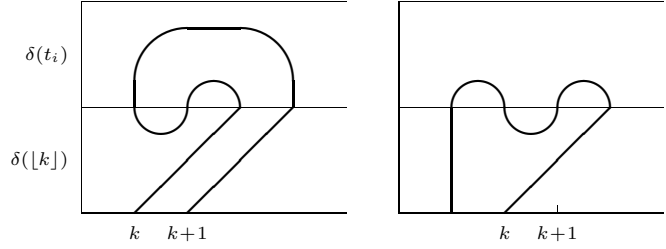
\begin{center}
\begin{picture}(300,110)
{\linethickness{0.05pt}
\put(80,20){\line(1,0){100}}
\put(80,60){\line(1,0){100}}
\put(80,100){\line(1,0){100}}
\put(80,20){\line(0,1){80}}}

{\linethickness{0.05pt}
\put(200,20){\line(1,0){100}}
\put(200,60){\line(1,0){100}}
\put(200,100){\line(1,0){100}}
\put(200,20){\line(0,1){80}}
\put(260,20){\line(0,1){3}}}

\put(75,40){\makebox(0,0)[r]{\scriptsize $\delta(\dk{k})$}}
\put(75,80){\makebox(0,0)[r]{\scriptsize $\delta(t_i)$}}
\put(0,100){\makebox(0,0)[l]{(1.2)}}
\put(100,15){\makebox(0,0)[t]{\scriptsize $k$}}
\put(120,15){\makebox(0,0)[t]{\scriptsize $k\!+\!1$}}

\put(240,15){\makebox(0,0)[t]{\scriptsize $k$}}
\put(260,15){\makebox(0,0)[t]{\scriptsize $k\!+\!1$}}

\thicklines
\put(110,60){\oval(20,20)[b]}
\put(130,60){\oval(20,20)[t]}
\put(130,60){\oval(60,60)[t]}
\put(100,20){\line(1,1){40}}
\put(120,20){\line(1,1){40}}

\put(230,60){\oval(20,20)[t]}
\put(250,60){\oval(20,20)[b]}
\put(270,60){\oval(20,20)[t]}
\put(220,20){\line(0,1){40}}
\put(240,20){\line(1,1){40}}

\end{picture}
\end{center}


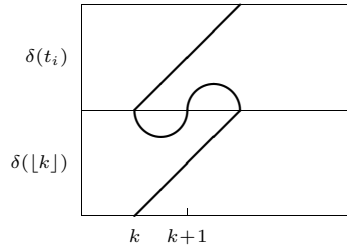
\begin{center}
\begin{picture}(300,110)
{\linethickness{0.05pt}
\put(80,20){\line(1,0){100}}
\put(80,60){\line(1,0){100}}
\put(80,100){\line(1,0){100}}
\put(80,20){\line(0,1){80}}
\put(120,20){\line(0,1){3}}}

\put(75,40){\makebox(0,0)[r]{\scriptsize $\delta(\dk{k})$}}
\put(75,80){\makebox(0,0)[r]{\scriptsize $\delta(t_i)$}}
\put(0,100){\makebox(0,0)[l]{(1.3)}}
\put(100,15){\makebox(0,0)[t]{\scriptsize $k$}}
\put(120,15){\makebox(0,0)[t]{\scriptsize $k\!+\!1$}}

\thicklines
\put(110,60){\oval(20,20)[b]}
\put(130,60){\oval(20,20)[t]}
\put(100,20){\line(1,1){40}}
\put(100,60){\line(1,1){40}}

\end{picture}
\end{center}


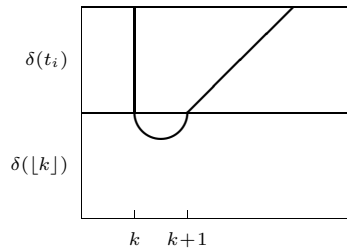
\begin{center}
\begin{picture}(300,110)
{\linethickness{0.05pt}
\put(80,20){\line(1,0){100}}
\put(80,60){\line(1,0){100}}
\put(80,100){\line(1,0){100}}
\put(80,20){\line(0,1){80}}
\put(100,20){\line(0,1){3}}
\put(120,20){\line(0,1){3}}}

\put(75,40){\makebox(0,0)[r]{\scriptsize $\delta(\dk{k})$}}
\put(75,80){\makebox(0,0)[r]{\scriptsize $\delta(t_i)$}}
\put(0,100){\makebox(0,0)[l]{(2)}}
\put(100,15){\makebox(0,0)[t]{\scriptsize $k$}}
\put(120,15){\makebox(0,0)[t]{\scriptsize $k\!+\!1$}}

\thicklines
\put(110,60){\oval(20,20)[b]}
\put(100,60){\line(0,1){40}}
\put(120,60){\line(1,1){40}}

\end{picture}
\end{center}

In cases $(1.1)$, $(1.2)$ and $(1.3)$ we have
\[
\begin{array}{ll}
(1) & \cu(t_i)=\cu(\dk{k}t_i),
\\
    & \ca(t_i)=\ca(\dk{k}t_i)+1;
\end{array}
\]
and in case $(2)$ we have
\[
\begin{array}{ll}
(2) & \cu(t_i)=\cu(\dk{k}t_i)-1,
\\
    & \ca(t_i)=\ca(\dk{k}t_i).
\end{array}
\]
When for both $\delta(t_1)$ and $\delta(t_2)$ we have $(1)$, it is clear that
$\beta(t_1,t_2)=\beta(\dk{k}t_1,\dk{k}t_2)$, and the same if for both we have
$(2)$. If for one of $\delta(t_i)$ we have $(1)$ and for the other $(2)$,
we have again this equality of balances.

We show analogously that $\beta(t_1,t_2)=\beta(t_1\gk{k},t_2\gk{k})$, and
since we have trivially that
\[
\begin{array}{rl}
\beta(t_1,t_2) & =\beta(\gk{k}t_1,\gk{k}t_2)
\\
               & =\beta(t_1\dk{k},t_2\dk{k})
\\
               & =\beta(\mj t_1,\mj t_2)
\\
               & =\beta(t_1\mj,t_2\mj),
\end{array}
\]
we can conclude that for every term $s$ of \Lo\
\[
\begin{array}{rl}
\beta(t_1,t_2) & =\beta(st_1,st_2)
\\
               & =\beta(t_1s,t_2s).
\end{array}
\]
From that the lemma follows.
\qed

\vspace{.2cm}

We have seen in Section 6 that every
$\cal K$-equivalence class of a
frieze may be identified
with a pair $(S_{\Theta},l)$ where $S_{\Theta}$ is a set of
nonoverlapping segments
in $\Z\!-\!\{0\}$, and $l$ is the number of circular components. For
\M-friezes $l$ is always $\omega$, and hence every \M-frieze is identified up to
$\cal K$-equivalence with $S_{\Theta}$. To identify $S_{\Theta}$ of a \M-frieze
it is enough to identify the rooted tree that makes the branching part of the tree
of $S_{\Theta}$, which is the part of $S_{\Theta}$ from the leaves down to the
lowest node after which no node is branching. (This tree may consist of
a single node.) We call this rooted subtree of $S_{\Theta}$ the
{\it crown} of $S_{\Theta}$. In the example in Section 6,
the lowest node after which no node is branching is $[-11,7]$ and
the crown is the tree above $[-11,7]$, whose root is this node.

The transversal thread in the \M-frieze that corresponds to the root of the
crown
will be called the {\it crown thread}. Every \M-frieze has a crown thread.
If the end points of a crown thread are $(k,0)$ and $(l,a)$, we call $(k,l)$
the {\it crown pair}. So in our example the crown pair is $(11,7)$. All
threads on the right-hand side of the crown thread are transversal threads identified
with $[-(k+n),l+n]$, for $n\geq 1$ and $(k,l)$ the crown pair.

Suppose $t=u$ does not hold in \Mo, and let \X\ be as before \Mo\ plus $t=u$.
Let $\gk{j_1}\ldots\gk{j_m}\dk{i_1}\ldots\dk{i_n}$ be the normal form of $t$,
and $\gk{k_1}\ldots\gk{k_p}\dk{l_1}\ldots\dk{l_q}$ the normal form of $u$. We
show that there is a term $v$ (built out of $t$ and $u$) such that
$v=\mj$ in \X, but not in \Mo.

Let $s_1$ be $\dk{k_p+1} \ldots \dk{k_1+1}$,
let $s_2$ be $\gk{l_q+1} \ldots \gk{l_1+1}$,
let $s'_1$ be $\dk{j_m+1}$ $\ldots \dk{j_1+1}$,
and, finally, let $s'_2$ be $\gk{i_n+1} \ldots \gk{i_1+1}$.
It is clear that by ({\it cup-cap} 3) we have
$s_1us_2=\mj$ and $s'_1ts'_2=\mj$ in \Mo,
while $s_1ts_2=\mj$ and $s'_1us'_2=\mj$ hold in \X.

($i$) If $m>p$, then $s_1ts_2=\mj$ cannot hold in \Mo, because the normal
form of $s_1ts_2$ has at least one cap, and analogously if $n>q$,
because then the normal form of $s_1ts_2$ has at least one cup.
If $m<p$ or $n<q$, then $s'_1us'_2=\mj$ cannot hold in \Mo.

($ii$) If $m=p$ and $n=q$, then we proceed by induction on $m+n$. If $m+n=1$,
then we have in \X\ either $\gk{j_1}=\gk{k_1}$ for $j_1\neq k_1$, or
$\dk{i_1}=\dk{l_1}$ for $i_1\neq l_1$. If $j_1<k_1$, then we have
$\gk{j_1}\dk{k_1-1}=\mj$ in \X, but not in \Mo. We proceed analogously
in the other cases of $j_1\neq k_1$ and $i_1\neq l_1$.

Suppose now $m+n>1$ and $m\geq 1$. Let a {\it cap-block} $\gk{r,\ldots,r-k}$
be $\gk{r}\gk{r-1}\ldots\gk{r-k}$ for $r\in \N^+$ and $k\in \N$. Then the
sequence of caps $\gk{j_1}\ldots\gk{j_m}$ can be written in terms of
cap-blocks as
\[
\gk{r_1,\ldots,r_1-k_1}\ldots \gk{r_h,\ldots,r_h-k_h}
\]
such that $1\leq h\leq m$, $r_1=j_1$, $r_h-k_h=j_m$ and
$r_i-k_i-r_{i+1}\geq 2$. Let $\gk{r'_1,\ldots,r'_1-k'_1}$ be the
leftmost cap-block of $\gk{k_1}\ldots\gk{k_m}$, as
$\gk{r_1,\ldots,r_1-k_1}$ is the leftmost cap-block of
$\gk{j_1}\ldots\gk{j_m}$. We have
\[
\dk{r_1+k_1+1}\gk{r_1,\ldots,r_1-k_1}=\gk{r_1,\ldots,r_1-k_1+1}.
\]
If $r_1+k_1=r'_1+k'_1$, then
\[
\dk{r_1+k_1+1}\gk{r'_1,\ldots,r'_1-k'_1}=\gk{r'_1,\ldots,r'_1-k'_1+1},
\]
and from $t=u$ in \X, we obtain $\dk{r_1+k_1+1}t=\dk{r_1+k_1+1}u$ in \X,
but not in \Mo, since the difference in the normal forms of $t$ and $u$
persists. We can then apply the induction hypothesis, since the new $m$
has decreased. If $r_1+k_1 > r'_1+k'_1$, then for some $r'$
\[
\dk{r_1+k_1+1}\gk{k_1}\ldots\gk{k_m}=\gk{k_1}\ldots\gk{k_m}\dk{r'}
\]
(as can be ascertained from the corresponding $\cal J$-friezes), and
with $\dk{r_1+k_1+1}t$ $=\dk{r_1+k_1+1}u$ we are in case ($i$). All the
other cases, where $r_1+k_1 < r'_1+k'_1$, and where $n\geq 1$, are dealt
with analogously.

We can verify that $\beta(t,u)=\beta(v,\mj)$. If $v$ is $s_1ts_2$, then
\[
\beta(t,u)=\beta(s_1ts_2,s_1us_2),
\]
as we have seen in the proof of the Balance Lemma, and the right-hand side
is equal to $\beta(v,\mj)$. We reason analogously for the other
possible forms of $v$.

If $(k,l)$ is the crown pair of $\delta(v)$, then $|k-l|=2\beta(v,\mj)$.
The numbers $k$ and $l$ cannot both be 1; otherwise, $v=\mj$ would hold in \Mo.
We have the following cases.

\vspace{.1cm}

(1)\quad Suppose $k=l>1$. If for some $i$ we have $v=h_i$ in \Mo, where $h_i$
abbreviates $\gk{i}\dk{i}$, then in \X\ we have $h_i=\mj$, and hence
\[
\begin{array}{rl}
h_i = & h_i h_{i\pm1} h_i,\quad {\mbox{\rm by ($h2$)}}
\\
   = & h_{i\pm 1},\quad {\mbox{\rm by (1)}}.
\end{array}
\]
So for every $k\in\N^+$ we have $h_k=\mj$ in \X.

If for every $i$ we don't have $v=h_i$ in \Mo, then $\dk{k-1}v\gk{k-1}=\mj$
in \X, and in the crown pair $(k',l')$ of $\dk{k-1}v\gk{k-1}$ we have $k'<k$ and
$l'<l$. This is clear from the following picture:

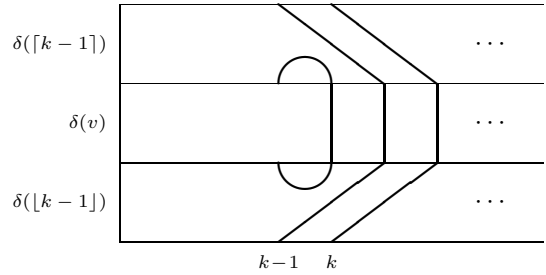
\begin{center}
\begin{picture}(200,130)
{\linethickness{0.05pt}
\put(40,20){\line(1,0){160}}
\put(40,50){\line(1,0){160}}
\put(40,80){\line(1,0){160}}
\put(40,110){\line(1,0){160}}
\put(40,20){\line(0,1){90}}}

\put(35,35){\makebox(0,0)[r]{\scriptsize $\delta(\dk{k-1})$}}
\put(35,65){\makebox(0,0)[r]{\scriptsize $\delta(v)$}}
\put(35,95){\makebox(0,0)[r]{\scriptsize $\delta(\gk{k-1})$}}
\put(100,15){\makebox(0,0)[t]{\scriptsize $k\!-\!1$}}
\put(120,15){\makebox(0,0)[t]{\scriptsize $k$}}

\put(180,35){\makebox(0,0){$\cdots$}}
\put(180,65){\makebox(0,0){$\cdots$}}
\put(180,95){\makebox(0,0){$\cdots$}}

\thicklines
\put(110,50){\oval(20,20)[b]}
\put(110,80){\oval(20,20)[t]}
\put(120,50){\line(0,1){30}}
\put(140,50){\line(0,1){30}}
\put(160,50){\line(0,1){30}}
\put(100,20){\line(4,3){40}}
\put(120,20){\line(4,3){40}}
\put(100,110){\line(4,-3){40}}
\put(120,110){\line(4,-3){40}}

\end{picture}
\end{center}

We would have $\dk{k-1}v\gk{k-1}=\mj$ in \Mo\ only if we had
$v=h_{k-2}$ in \Mo, as can be seen from the following picture:

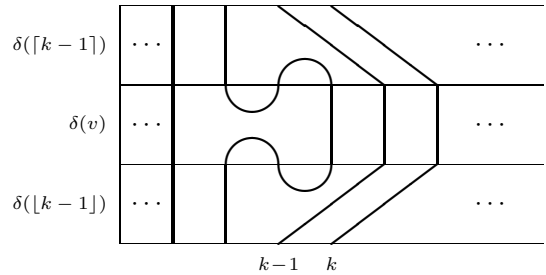
\begin{center}
\begin{picture}(200,130)
{\linethickness{0.05pt}
\put(40,20){\line(1,0){160}}
\put(40,50){\line(1,0){160}}
\put(40,80){\line(1,0){160}}
\put(40,110){\line(1,0){160}}
\put(40,20){\line(0,1){90}}}

\put(35,35){\makebox(0,0)[r]{\scriptsize $\delta(\dk{k-1})$}}
\put(35,65){\makebox(0,0)[r]{\scriptsize $\delta(v)$}}
\put(35,95){\makebox(0,0)[r]{\scriptsize $\delta(\gk{k-1})$}}
\put(100,15){\makebox(0,0)[t]{\scriptsize $k\!-\!1$}}
\put(120,15){\makebox(0,0)[t]{\scriptsize $k$}}

\put(180,35){\makebox(0,0){$\cdots$}}
\put(180,65){\makebox(0,0){$\cdots$}}
\put(180,95){\makebox(0,0){$\cdots$}}

\put(50,35){\makebox(0,0){$\cdots$}}
\put(50,65){\makebox(0,0){$\cdots$}}
\put(50,95){\makebox(0,0){$\cdots$}}

\thicklines
\put(110,50){\oval(20,20)[b]}
\put(110,80){\oval(20,20)[t]}
\put(90,50){\oval(20,20)[t]}
\put(90,80){\oval(20,20)[b]}
\put(120,50){\line(0,1){30}}
\put(140,50){\line(0,1){30}}
\put(160,50){\line(0,1){30}}
\put(60,20){\line(0,1){90}}
\put(80,20){\line(0,1){30}}
\put(80,80){\line(0,1){30}}
\put(100,20){\line(4,3){40}}
\put(120,20){\line(4,3){40}}
\put(100,110){\line(4,-3){40}}
\put(120,110){\line(4,-3){40}}

\end{picture}
\end{center}

\noindent There are several straightforward cases to consider in order to prove this assertion.
So $\dk{k-1}v\gk{k-1}=\mj$ does not hold in \Mo, and by induction we obtain
$h_i=\mj$ in \X\ for some $i$. So, as above, for every $k\in\N^+$ we have
$h_k=\mj$ in \X.

\vspace{.1cm}

(2)\quad If $k\neq l$, and $\min(k,l)=m_1>1$, then
$\dk{m_1-1}v\gk{m_1-1}=\mj$ in \X, and in the crown pair $(k',l')$ of
$\dk{m_1-1}v\gk{m_1-1}$ we have $k'<k$ and $l'<l$. This is clear from
the following picture:

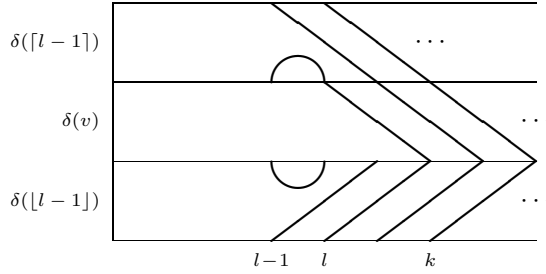
\begin{center}
\begin{picture}(200,130)
{\linethickness{0.05pt}
\put(40,20){\line(1,0){160}}
\put(40,50){\line(1,0){160}}
\put(40,80){\line(1,0){160}}
\put(40,110){\line(1,0){160}}
\put(40,20){\line(0,1){90}}}

\put(35,35){\makebox(0,0)[r]{\scriptsize $\delta(\dk{l-1})$}}
\put(35,65){\makebox(0,0)[r]{\scriptsize $\delta(v)$}}
\put(35,95){\makebox(0,0)[r]{\scriptsize $\delta(\gk{l-1})$}}
\put(100,15){\makebox(0,0)[t]{\scriptsize $l\!-\!1$}}
\put(120,15){\makebox(0,0)[t]{\scriptsize $l$}}
\put(160,15){\makebox(0,0)[t]{\scriptsize $k$}}

\put(160,95){\makebox(0,0){$\cdots$}}
\put(200,65){\makebox(0,0){$\cdots$}}
\put(200,35){\makebox(0,0){$\cdots$}}

\thicklines
\put(110,50){\oval(20,20)[b]}
\put(110,80){\oval(20,20)[t]}
\put(100,20){\line(4,3){40}}
\put(120,20){\line(4,3){40}}
\put(140,20){\line(4,3){40}}
\put(160,20){\line(4,3){40}}
\put(100,110){\line(4,-3){80}}
\put(120,110){\line(4,-3){80}}
\put(120,80){\line(4,-3){40}}

\end{picture}
\end{center}

We cannot have $\dk{m_1-1}v\gk{m_1-1}=\mj$ in \Mo\ because we must
have $|k'-l'|=|k-l|$, as it is clear from the picture. We continue
in the same manner until we reach a term
$\dk{m_p-1}\ldots\dk{m_1-1}v\gk{m_1-1}\ldots\gk{m_p-1}$, which we
abbreviate by $v'$, such that $v'=\mj$ in \X, and the crown pair of
$\delta(v')$ is either $(1,l-k+1)$ or $(k-l+1,1)$. Suppose
$l>k$. Then in $\delta(v')$ there must be a cup $[i,i+1]$, and we have
\[
h_i=v'h_i=v'=\mj.
\]
That $v'h_i=v'$ is clear from the picture

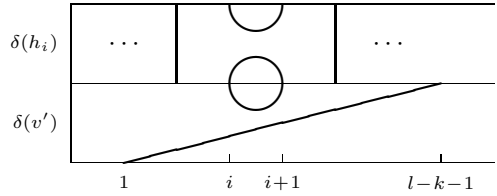
\begin{center}
\begin{picture}(200,100)
{\linethickness{0.05pt}
\put(40,20){\line(1,0){160}}
\put(40,50){\line(1,0){160}}
\put(40,80){\line(1,0){160}}
\put(40,20){\line(0,1){60}}
\put(100,20){\line(0,1){3}}
\put(120,20){\line(0,1){3}}
\put(180,20){\line(0,1){3}}}

\put(35,35){\makebox(0,0)[r]{\scriptsize $\delta(v')$}}
\put(35,65){\makebox(0,0)[r]{\scriptsize $\delta(h_i)$}}
\put(60,15){\makebox(0,0)[t]{\scriptsize $1$}}
\put(100,15){\makebox(0,0)[t]{\scriptsize $i$}}
\put(120,15){\makebox(0,0)[t]{\scriptsize $i\!+\!1$}}
\put(180,15){\makebox(0,0)[t]{\scriptsize $l\!-\!k\!-\!1$}}

\put(60,65){\makebox(0,0){$\cdots$}}
\put(160,65){\makebox(0,0){$\cdots$}}

\thicklines
\put(110,50){\oval(20,20)[b]}
\put(110,50){\oval(20,20)[t]}
\put(110,80){\oval(20,20)[b]}
\put(80,50){\line(0,1){30}}
\put(140,50){\line(0,1){30}}
\put(60,20){\line(4,1){120}}

\end{picture}
\end{center}

We proceed analogously when $k>l$. As before, from $h_i=\mj$ we derive
that for every $k\in\N^+$ we have $h_k=\mj$ in \X.

So, both in case (1) and in case (2), for every $k\in\N^+$ we have
$h_k=\mj$ in \X. Since
\[
\dk{i}\gk{i+1}\dk{i+1}=\dk{i+1},\quad {\mbox{\rm by ({\it cup-cap} 3)}},
\]
with $h_{i+1}=\mj$, we obtain $\dk{i}=\dk{i+1}$ in \X\ for every $i\in\N^+$.
Analogously, we obtain $\gk{i}=\gk{i+1}$ in \X\ for every $i\in\N^+$.

If $\dk{1}^0$ and $\gk{1}^0$ are $\mj$, while $\dk{1}^{k+1}$
is $\dk{1}^k\dk{1}$
and $\gk{1}^{k+1}$ is $\gk{1}^k\gk{1}$, then every element of \X\
is either of the form $\dk{1}^k$ or of the form
$\gk{1}^k$ for some $k\in\N$. To see that,
start from the normal form of an element, identify all
cups with $\dk{1}$, all
caps with $\gk{1}$, and then use $\gk{1}\dk{1}=\mj$.

So all the elements of \X\ are the following
\[
\ldots,\gk{1}^3,\gk{1}^2,\gk{1}^1,\mj,\dk{1}^1,\dk{1}^2,\dk{1}^3,\ldots
\]
If $\beta(t,u)=0$, then the Balance Lemma guarantees that all the elements of
\X\ above are mutually distinct. Composition in \X\ then behaves as addition of
integers, where \mj\ is zero, and so \X\ is isomorphic to $\Z=\Z/0$.

In case $\beta(t,u)=n>0$, the elements of \X\ above are not all mutually distinct.
We are then in case (2), and from $v'=\mj$ we can infer $\dk{1}^n=\mj$, and
also $\gk{1}^n=\mj$, in \X. Hence every element of \X\ is equal to one
of the following
\[
\mj,\dk{1}^1,\ldots \dk{1}^{n-1}\; ,
\]
and by the Balance Lemma these are all mutually distinct. Composition in \X\ then
behaves as addition of integers modulo $n$, where \mj\ is zero, and so
\X\ is isomorphic to \Z/$n$. So, in any
case, we have that \X\ is isomorphic to $\Z/\beta(t,u)$.

\section{The monoids \Lpm\ and \Kpm}

The monoid \Lpm\ is defined as \Lo\ save that for every $k\in\Z$
there is a generator $\dk{k}$ and a generator $\gk{k}$, and there are two
additional generators \ss\ and \bs.
The equations of \Lpm\ are those of \Lo\ plus

\[
\begin{array}{l}
\dk{k}\ss=\ss\dk{k+1},
\\[0.1cm]
\gk{k}\ss=\ss\gk{k+1},
\\[0.1cm]
\bs\ss=\ss\bs=\mj.
\end{array}
\]

We derive easily
\[
\begin{array}{l}
\dk{k+1}\bs=\bs\dk{k},
\\[0.1cm]
\gk{k+1}\bs=\bs\gk{k}.
\end{array}
\]

The monoid \Kpm\ has in addition the equation ({\it cup-cap} 4).

Let $\ss^0$ be the empty sequence, while $\ss^{k+1}$ is $\ss^k\ss$,
and analogously
for $\bs^k$. We can define the normal form for terms of \Lpm\
as the normal form
for \Lo\ of Section 4 prefixed with either $\ss^k$
or $\bs^k$ for $k\in\N$. It is clear that combined with the
Normal Form Lemma of
Section 4 the equations above enable us
to reduce every term
to a term in normal form equal to the original term.

Let now $R_a$ be $(-\infty,\infty)\times[0,a]$. An $\pm\omega$-{\it diagram}
in $R_a$ is defined as an $\omega$-diagram save that $\N^+$ is replaced by
\Z. A $\pm$ {\it frieze} is a $\pm\omega$-diagram with a finite number of cups,
caps and circular components.

The cups and caps $\dk{k}$ and $\gk{k}$ are mapped to $\pm$ friezes
analogously to what we had before, while \ss\ and \bs\ are mapped into the
$\cal L$-equivalence classes of the following $\pm$ friezes:

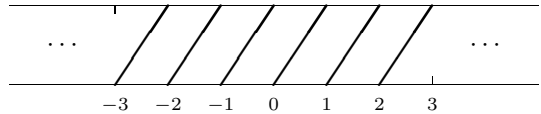
\begin{center}
\begin{picture}(200,70)
{\linethickness{0.05pt}
\put(0,20){\line(1,0){200}}
\put(0,50){\line(1,0){200}}
\put(40,50){\line(0,-1){3}}
\put(160,20){\line(0,1){3}}}

\put(40,15){\makebox(0,0)[t]{\scriptsize $-3$}}
\put(60,15){\makebox(0,0)[t]{\scriptsize $-2$}}
\put(80,15){\makebox(0,0)[t]{\scriptsize $-1$}}
\put(100,15){\makebox(0,0)[t]{\scriptsize $0$}}
\put(120,15){\makebox(0,0)[t]{\scriptsize $1$}}
\put(140,15){\makebox(0,0)[t]{\scriptsize $2$}}
\put(160,15){\makebox(0,0)[t]{\scriptsize $3$}}

\put(20,35){\makebox(0,0){$\cdots$}}
\put(180,35){\makebox(0,0){$\cdots$}}

\thicklines
\put(40,20){\line(2,3){20}}
\put(60,20){\line(2,3){20}}
\put(80,20){\line(2,3){20}}
\put(100,20){\line(2,3){20}}
\put(120,20){\line(2,3){20}}
\put(140,20){\line(2,3){20}}

\end{picture}
\end{center}


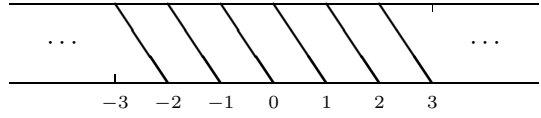
\begin{center}
\begin{picture}(200,70)
{\linethickness{0.05pt}
\put(0,20){\line(1,0){200}}
\put(0,50){\line(1,0){200}}
\put(40,20){\line(0,1){3}}
\put(160,50){\line(0,-1){3}}}

\put(40,15){\makebox(0,0)[t]{\scriptsize $-3$}}
\put(60,15){\makebox(0,0)[t]{\scriptsize $-2$}}
\put(80,15){\makebox(0,0)[t]{\scriptsize $-1$}}
\put(100,15){\makebox(0,0)[t]{\scriptsize $0$}}
\put(120,15){\makebox(0,0)[t]{\scriptsize $1$}}
\put(140,15){\makebox(0,0)[t]{\scriptsize $2$}}
\put(160,15){\makebox(0,0)[t]{\scriptsize $3$}}

\put(20,35){\makebox(0,0){$\cdots$}}
\put(180,35){\makebox(0,0){$\cdots$}}

\thicklines
\put(40,50){\line(2,-3){20}}
\put(60,50){\line(2,-3){20}}
\put(80,50){\line(2,-3){20}}
\put(100,50){\line(2,-3){20}}
\put(120,50){\line(2,-3){20}}
\put(140,50){\line(2,-3){20}}

\end{picture}
\end{center}

One could
conceive \ss\ss, i.e. \ss\ multiplied with itself, as the
infinite block $\ldots h_2 h_1 h_0 h_{-1} h_{-2}\ldots$,
and \bs\bs\ as the
infinite block $\ldots h_{-2} h_{-1} h_0 h_1 h_2\ldots$

That the monoid \Lpm\ is isomorphic to a monoid made of $\cal L$-equivalence
classes of $\pm$ friezes is shown analogously to what we had for \Lo, and the
same can be shown with $\cal L$ replaced by $\cal K$.
From these isomorphisms it follows that \Lo\ can be embedded in \Lpm, and analogously
with $\cal K$. This is because we can identify every frieze with a $\pm$ frieze
such that for every $z\leq 0$ we have a vertical thread with the end points
$(z,0)$ and $(z,a)$.

The monoid \Lnc\ is defined as the monoid \Ln\ save that we have also the diapsides
$h_0$ and $h_n$. To the equations of \Ln\ we add the equations
\[
\begin{array}{l}
h_0=h_n,
\\[0.1cm]
c_1^{\alpha}=c_{n+1}^{\alpha},
\\[0.1cm]
h_k\ss=\ss h_{k+1},\quad {\mbox{\rm for }} k\in\{0,\ldots,n-1\}
\\[0.1cm]
\bs\ss=\ss\bs=\mj.
\end{array}
\]
The monoid \Knc\ is obtained from \Lnc\ as \Kn\ is obtained from \Ln.

Note that while \Ln\ was a submonoid of \Lo, because of $h_0=h_n$ we don't have
that \Lnc\ is a submonoid of \Lpm. For the same reason \Lnc\ is not a submonoid
of ${\cal L}_{n+1}^{\mbox{\scriptsize\rm cyl}}$, while \Ln\ was isomorphic to
two submonoids of ${\cal L}_{n+1}$ (we map $h_i$ either
to $h_i$ or to $h_{i+1}$).
The same holds when $\cal L$ is replaced by $\cal K$.

The monoid \Lnc\ may be shown isomorphic to a monoid made of equivalence
classes of cylindric friezes, which are roughly defined as follows. Instead
of diagrams in $R_a$ we now have diagrams in cylinders where the top
and bottom are copies of a circle with $n$ points labelled counterclockwise with
the numbers from 1 to $n$. We interpret $h_i$ for $i\in\{1,\ldots,n-1\}$ by

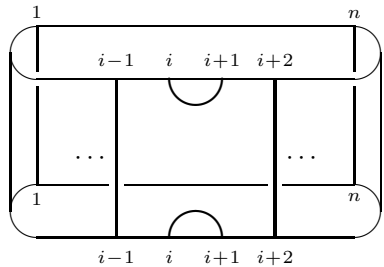
\begin{center}
\begin{picture}(140,120)

\put(10,20){\line(1,0){120}}
\put(10,40){\line(1,0){27}}
\put(43,40){\line(1,0){54}}
\put(103,40){\line(1,0){27}}
\put(10,30){\oval(20,20)[l]}
\put(130,30){\oval(20,20)[r]}

\put(10,80){\line(1,0){120}}
\put(10,100){\line(1,0){120}}
\put(10,90){\oval(20,20)[l]}
\put(130,90){\oval(20,20)[r]}
\put(0,30){\line(0,1){60}}
\put(140,30){\line(0,1){60}}

\put(10,37){\makebox(0,0)[t]{\scriptsize $1$}}
\put(10,103){\makebox(0,0)[b]{\scriptsize $1$}}
\put(40,15){\makebox(0,0)[t]{\scriptsize $i\!-\!1$}}
\put(60,15){\makebox(0,0)[t]{\scriptsize $i$}}
\put(80,15){\makebox(0,0)[t]{\scriptsize $i\!+\!1$}}
\put(100,15){\makebox(0,0)[t]{\scriptsize $i\!+\!2$}}
\put(130,37){\makebox(0,0)[t]{\scriptsize $n$}}
\put(130,103){\makebox(0,0)[b]{\scriptsize $n$}}
\put(40,85){\makebox(0,0)[b]{\scriptsize $i\!-\!1$}}
\put(60,85){\makebox(0,0)[b]{\scriptsize $i$}}
\put(80,85){\makebox(0,0)[b]{\scriptsize $i\!+\!1$}}
\put(100,85){\makebox(0,0)[b]{\scriptsize $i\!+\!2$}}

\put(30,50){\makebox(0,0){$\cdots$}}
\put(110,50){\makebox(0,0){$\cdots$}}

\thicklines
\put(70,80){\oval(20,20)[b]}
\put(70,20){\oval(20,20)[t]}
\put(40,20){\line(0,1){60}}
\put(100,20){\line(0,1){60}}
\put(10,40){\line(0,1){37}}
\put(10,83){\line(0,1){17}}
\put(130,40){\line(0,1){37}}
\put(130,83){\line(0,1){17}}

\end{picture}
\end{center}

\noindent while $h_0$ and $h_n$ are interpreted by

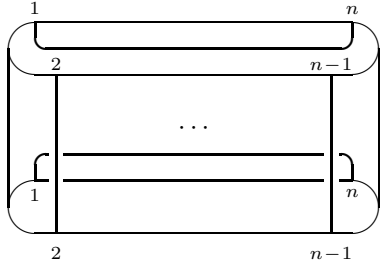
\begin{center}
\begin{picture}(140,120)

\put(10,20){\line(1,0){120}}
\put(10,40){\line(1,0){5}}
\put(21,40){\line(1,0){98}}
\put(125,40){\line(1,0){5}}
\put(10,30){\oval(20,20)[l]}
\put(130,30){\oval(20,20)[r]}

\put(10,80){\line(1,0){120}}
\put(10,100){\line(1,0){120}}
\put(10,90){\oval(20,20)[l]}
\put(130,90){\oval(20,20)[r]}
\put(0,30){\line(0,1){60}}
\put(140,30){\line(0,1){60}}

\put(10,37){\makebox(0,0)[t]{\scriptsize $1$}}
\put(10,103){\makebox(0,0)[b]{\scriptsize $1$}}
\put(130,37){\makebox(0,0)[t]{\scriptsize $n$}}
\put(130,103){\makebox(0,0)[b]{\scriptsize $n$}}
\put(18,15){\makebox(0,0)[t]{\scriptsize $2$}}
\put(18,82){\makebox(0,0)[b]{\scriptsize $2$}}
\put(122,15){\makebox(0,0)[t]{\scriptsize $n\!-\!1$}}
\put(122,82){\makebox(0,0)[b]{\scriptsize $n\!-\!1$}}

\put(70,60){\makebox(0,0){$\cdots$}}

\thicklines
\put(15,40){\oval(10,20)[tl]}
\put(15,100){\oval(10,20)[bl]}
\put(125,40){\oval(10,20)[tr]}
\put(125,100){\oval(10,20)[br]}
\put(18,20){\line(0,1){60}}
\put(122,20){\line(0,1){60}}
\put(21,50){\line(1,0){98}}
\put(15,90){\line(1,0){110}}

\end{picture}
\end{center}

\noindent We interpret \ss\ and \bs\ by
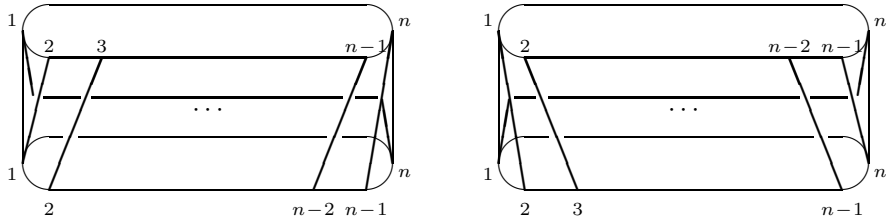
\begin{center}
\begin{picture}(320,110)

\put(10,20){\line(1,0){120}}
\put(10,40){\line(1,0){5}}
\put(21,40){\line(1,0){94}}
\put(121,40){\line(1,0){9}}
\put(10,30){\oval(20,20)[l]}
\put(130,30){\oval(20,20)[r]}

\put(10,70){\line(1,0){120}}
\put(10,90){\line(1,0){120}}
\put(10,80){\oval(20,20)[l]}
\put(130,80){\oval(20,20)[r]}
\put(0,30){\line(0,1){50}}
\put(140,30){\line(0,1){50}}

\put(-2,28){\makebox(0,0)[tr]{\scriptsize $1$}}
\put(-2,82){\makebox(0,0)[br]{\scriptsize $1$}}
\put(10,15){\makebox(0,0)[t]{\scriptsize $2$}}
\put(10,72){\makebox(0,0)[b]{\scriptsize $2$}}
\put(30,72){\makebox(0,0)[b]{\scriptsize $3$}}
\put(110,15){\makebox(0,0)[t]{\scriptsize $n\!-\!2$}}
\put(130,15){\makebox(0,0)[t]{\scriptsize $n\!-\!1$}}
\put(130,72){\makebox(0,0)[b]{\scriptsize $n\!-\!1$}}
\put(142,28){\makebox(0,0)[tl]{\scriptsize $n$}}
\put(142,82){\makebox(0,0)[bl]{\scriptsize $n$}}

\put(70,50){\makebox(0,0){$\cdots$}}

{\thicklines
\put(0,80){\line(1,-6){4}}
\put(140,30){\line(-1,6){4}}
\put(0,30){\line(1,4){10}}
\put(10,20){\line(2,5){20}}
\put(110,20){\line(2,5){20}}
\put(130,20){\line(1,6){10}}
\put(8,55){\line(1,0){14}}
\put(26,55){\line(1,0){95}}
\put(126,55){\line(1,0){8}}}

\put(190,20){\line(1,0){120}}
\put(190,40){\line(1,0){9}}
\put(205,40){\line(1,0){94}}
\put(305,40){\line(1,0){5}}
\put(190,30){\oval(20,20)[l]}
\put(310,30){\oval(20,20)[r]}

\put(190,70){\line(1,0){120}}
\put(190,90){\line(1,0){120}}
\put(190,80){\oval(20,20)[l]}
\put(310,80){\oval(20,20)[r]}
\put(180,30){\line(0,1){50}}
\put(320,30){\line(0,1){50}}

\put(178,28){\makebox(0,0)[tr]{\scriptsize $1$}}
\put(178,82){\makebox(0,0)[br]{\scriptsize $1$}}
\put(190,15){\makebox(0,0)[t]{\scriptsize $2$}}
\put(190,72){\makebox(0,0)[b]{\scriptsize $2$}}
\put(210,15){\makebox(0,0)[t]{\scriptsize $3$}}
\put(290,72){\makebox(0,0)[b]{\scriptsize $n\!-\!2$}}
\put(310,15){\makebox(0,0)[t]{\scriptsize $n\!-\!1$}}
\put(310,72){\makebox(0,0)[b]{\scriptsize $n\!-\!1$}}
\put(322,28){\makebox(0,0)[tl]{\scriptsize $n$}}
\put(322,82){\makebox(0,0)[bl]{\scriptsize $n$}}

\put(250,50){\makebox(0,0){$\cdots$}}

\thicklines
\put(180,30){\line(1,6){4}}
\put(320,80){\line(-1,-6){4}}
\put(320,30){\line(-1,4){10}}
\put(310,20){\line(-2,5){20}}
\put(210,20){\line(-2,5){20}}
\put(190,20){\line(-1,6){10}}
\put(186,55){\line(1,0){8}}
\put(199,55){\line(1,0){95}}
\put(298,55){\line(1,0){14}}

\end{picture}
\end{center}

Cylindric friezes are special three-dimensional tangles, whereas with
friezes and
$\pm$ friezes we had only two-dimensional tangles. In these special
three-dimensional
tangles we have only ``cyclic braidings'' or ``torsions'' like
those obtained from the last two pictures, and further ``cyclic braidings''
obtained by composing these.
Diagrams like our cylindric friezes were considered in \cite{J94}.

\section{Self-adjunctions}

A {\it self-adjunction}, which we will also call $\cal L$-{\it adjunction},
is an adjunction in which an endofunctor is adjoint to itself, which means
that it is both left and right adjoint to itself
(for the general notion of adjunction see \cite{ML71}, Chapter IV).
More precisely, an
$\cal L$-adjunction is $\langle{\cal A}, \cirk,\mj,F,\varphi,\gamma\rangle$
where
$\langle{\cal A},\cirk,\mj\rangle$ is a category, which means that for
$f:a\str b$, $g:b\str c$ and $h:c\str d$ arrows of $\cal A$ we have
the equations
\[
\begin{array}{ll}
{\makebox[1cm][l]{$({\mbox{\it cat }} 1)$}} &
{\makebox[6cm][l]{$f\cirk\mj_a=f, \quad \mj_b\cirk f=f,$}}
\\
({\mbox{\it cat }} 2) & h\cirk(g\cirk f)=(h\cirk g)\cirk f;
\end{array}
\]
$F$ is a functor from $\cal A$ to $\cal A$, which means that we have
the equations
\[
\begin{array}{ll}
{\makebox[1cm][l]{$({\mbox{\it fun }} 1)$}} &
{\makebox[6cm][l]{$F\mj_a=\mj_{Fa},$}}
\\
({\mbox{\it fun }} 2) & F(g\cirk f)=Fg\cirk Ff;
\end{array}
\]
$\varphi$ is a natural transformation (the counit of the adjunction) with
components $\varphi_a:FFa\str a$, and $\gamma$ is a natural transformation
(the unit of the adjunction) with components $\gamma_a:a\str FFa$,
which means that we have the equations
\[
\begin{array}{ll}
{\makebox[1cm][l]{$({\mbox{\it nat }} \varphi)$}} &
{\makebox[6cm][l]{$f\cirk\varphi_a=\varphi_b\cirk FFf,$}}
\\
({\mbox{\it nat }} \gamma) & FFf\cirk\gamma_a=\gamma_b\cirk f;
\end{array}
\]
and, finally, we have the {\it triangular equations}
\[
\begin{array}{ll}
{\makebox[1cm][l]{$(\varphi\gamma)$}} &
{\makebox[6cm][l]{$F{\varphi_a}\cirk\gamma_{Fa}
=\varphi_{Fa}\cirk F\gamma_a=\mj_{Fa}.$}}
\end{array}
\]
We will call the equations from ({\it cat} 1)  to ($\varphi\gamma$) we have displayed
above the $\cal L$-{\it equations}.

A $\cal K$-{\it adjunction} is an $\cal L$-adjunction that satisfies the additional
equation
\[
\begin{array}{ll}
{\makebox[1cm][l]{$(\varphi\gamma{\cal K})$}} &
{\makebox[6cm][l]{$F(\varphi_a\cirk\gamma_a)=\varphi_{Fa}\cirk\gamma_{Fa}.$}}
\end{array}
\]
The $\cal L$-equations plus this equation will be called $\cal K$-{\it equations}.

Let $\kappa_a$ be an abbreviation for $\varphi_a\cirk\gamma_a:a\str a$. Then
in every $\cal L$-adjunction we have that
$f\cirk\kappa_a=\kappa_b\cirk f$, and
($\varphi\gamma{\cal K}$) is expressed by $F\kappa_a=\kappa_{Fa}$.
We will see that the arrows $\kappa_a$ are in general not equal
to the identity
arrows $\mj_a$ in arbitrary $\cal K$-adjunctions, but they have
some properties
of identity arrows: they commute with other arrows, and they are preserved
by the functor $F$.

A $\cal J$-{\it adjunction} is an $\cal L$-adjunction that satisfies the additional
equation
\[
\begin{array}{ll}
(\varphi\gamma{\cal J}) & \varphi_a\cirk\gamma_a=\mj_a,
\end{array}
\]
i.e. $\kappa_a=\mj_a$. The $\cal L$-equations plus this equation will be
called $\cal J$-{\it equations}. Every $\cal J$-adjunction is a $\cal K$-adjunction,
but not vice versa, as we will see later.

\section{Free self-adjunctions}

The {\it free} $\cal L$-{\it adjunction} generated by an arbitrary object, which
we will denote by $0$, is defined as follows. The category of this self-adjunction,
which we will call \Lc, has the objects $0$, $F0$, $FF0$, etc., which may be
identified with the natural numbers 0, 1, 2, etc.

An {\it arrow-term} of \Lc\ will be a word $f$ that has a {\it type}
$(n,m)$, where $n,m\in\N$. That $f$ is of type $(n,m)$ is expressed by
$f:n\str m$. Now we define the arrow-terms of \Lc\ inductively. We stipulate
first for every $n\in\N$ that $\mj_n:n\str n$, $\varphi_n:n+2\str n$ and
$\gamma_n:n\str n+2$ are arrow-terms of \Lc. Next, if $f:m\str n$ is
an arrow-term of \Lc, then $Ff:m+1\str n+1$ is an arrow-term of \Lc, and if
$f:m\str n$ and $g:n\str k$ are arrow-terms of \Lc, then $(g\cirk f):m\str k$
is an arrow-term of \Lc. As usual, we don't write parentheses
in $(g\cirk f)$ when they are not essential.

On these arrow-terms we impose the $\cal L$-equations, where $a$ and $b$
are replaced by $m$ and $n$,
and $f$, $g$ and $h$ stand for arrow-terms of \Lc. Formally, we take the smallest equivalence
relation $\equiv$ on the arrow-terms of \Lc\ satisfying, first,
congruence conditions with respect to $F$ and $\cirk$, namely,
\[
\begin{array}{l}
{\mbox{\rm if }} f\equiv g,\quad {\mbox{\rm then }} Ff\equiv Fg,
\\
{\mbox{\rm if }} f_1\equiv f_2\;{\mbox{\rm and }}g_1\equiv g_2,
\quad {\mbox{\rm then }} g_1\cirk f_1\equiv g_2\cirk f_2,
\end{array}
\]
provided $g_1\cirk f_1$ and $g_2\cirk f_2$ are defined, and, second,
the conditions obtained from the $\cal L$-equations by replacing
the equality sign by $\equiv$. Then we take the equivalence classes of arrow-terms
as arrows, with the obvious source and target, all arrow-terms in the same
class having the same type. On these equivalence classes we define
\mj, $\varphi$, $\gamma$, $F$ and $\cirk$ in the obvious way. This defines
the category \Lc, in which we have clearly an $\cal L$-adjunction.

The category \Lc\ satisfies the following universal property. If $\lambda$
maps the object $0$ into an arbitrary object of the category $\cal A$ of an
arbitrary $\cal L$-adjunction, then there is  a unique functor $\Lambda$ of
$\cal L$-adjunctions (defined in the obvious way, so that the $\cal L$-adjunction structure
is preserved) such that $\Lambda$ maps $0$ into $\lambda(0)$.
This property characterizes \Lc\ up to isomorphism with a functor
of $\cal L$-adjunctions. This justifies calling {\it free} the $\cal L$-adjunction
of \Lc.

The category \Kc\ of the {\it free} $\cal K$-{\it adjunction}
and the
category \Mc\ of the {\it free} $\cal J$-{\it adjunction}, both
generated by $0$, are defined as \Lc\ save that we replace everywhere
$\cal L$ by $\cal K$ and $\cal J$ respectively. So the equations
$(\varphi\gamma{\cal K})$ and $(\varphi\gamma{\cal J})$ come into play.
The categories \Kc\ and \Mc\ satisfy universal properties analogous to
the one above.

Let $\kappa_m$ be $\varphi_m\cirk\gamma_m$ as in the preceding section,
let $\kappa_m^0$ be $\mj_m$, and let $\kappa_m^{l+1}$ be $\kappa_m^l\cirk\kappa_m$.
It is easy to show by induction on the length of derivation that if
$f=g$ in \Mc, then for some $k,l\geq 0$ we have $f\cirk\kappa_0^k=g\cirk\kappa_0^l$
in \Kc. (The converse implication holds trivially.)
We will establish in the next section that
if $f=g$ in \Mc\ but not in \Kc, then $k$ must be different from $l$.

\section{\Lc\ and \Lo}

Let $F^0$ be the empty sequence, and let $F^{k+1}$ be $F^kF$. On the arrows
of \Lc, we define
a total binary operation $\ast$ based on
composition of
arrows in the following manner. For $f:m\str n$ and $g:k\str l$,
\[
g\ast f=_{\mbox{\scriptsize\it def}}\left\{
\begin{array}{ll}
g\cirk F^{k-n}f & {\mbox{\rm if }}n\leq k
\\
F^{n-k}g\cirk f & {\mbox{\rm if }}k\leq n.
\end{array}
\right.
\]

Next, let $f\eql g$ iff there are $k,l\in\N$ such that
$F^kf=F^lg$ in \Lc. It is easy to check that $\eql$ is an
equivalence relation on the arrows of \Lc, which satisfies moreover
\[
({\mbox{\it congr}}\;\ast)\quad {\mbox{\rm if }}f_1\eql f_2\;
{\mbox{\rm and }} g_1\eql g_2, \;
{\mbox{\rm then }}g_1\ast f_1\eql g_2\ast f_2.
\]

For every arrow $f$ of \Lc, let $[f]$ be $\{g\mid f\eql g\}$, and let $\Lc^{\ast}$ be
$\{[f]\mid f\; {\mbox{\rm is an arrow of \Lc}}\}$. With
\[
\begin{array}{l}
\mj=_{\mbox{\scriptsize\it def}}[\mj_0],
\\[.1cm]
[g][f]=_{\mbox{\scriptsize\it def}}[g\ast f],
\end{array}
\]
we can check that $\Lc^{\ast}$ is a monoid. We will show that this
monoid is isomorphic to the monoid \Lo.

Consider the map $\psi$ from the arrow-terms of \Lc\ to
the terms of \Lo\ defined inductively by

\[
\begin{array}{lcl}
\psi(\mj_n) & {\mbox{\rm is}} & \mj,
\\
\psi(\varphi_n) & {\mbox{\rm is}} & \dk{n+1},
\\
\psi(\gamma_n) & {\mbox{\rm is}} & \gk{n+1},
\\
\psi(Ff) & {\mbox{\rm is}} & \psi(f),
\\
\psi(g\cirk f) & {\mbox{\rm is}} & \psi(g)\psi(f).
\end{array}
\]

Let $m,n\in\N$. If in a frieze we have for every $k\in\N^+$
that $(m+k,0)$ and $(n+k,a)$ are the end points of a transversal thread,
i.e., $[-(m+k),n+k]$ identifies this thread, and there are no circular components
in the regions that correspond to $[-(m+k+1),n+k+1]$, i.e., the
ordinal 0 of circular components is assigned to this thread, then we say that this frieze
{\it is of type} $(n,m)$.
Note that a frieze of type $(n,m)$ is also of type $(n+k,m+k)$.
The $n$-friezes of
Section 9 are the friezes of type $(n,n)$.

We can easily establish the following by induction on the length of $f$.

\vspace{.2cm}

\noindent {\sc Remark I.}\quad {\it For every arrow-term} $f:n\str m$ {\it of} \Lc,
{\it the frieze} $\delta(\psi(f))$ {\it is of type} $(n,m)$.

\vspace{.2cm}

We also have the following.

\vspace{.2cm}

\noindent {\sc Remark II.}\quad {\it If the frieze} $\delta(t)$
{\it is of type} $(n,m)$, {\it then}
$\delta(t\dk{n+1})\cong_{\cal L}\delta(\dk{m+1}t)$ {\it and}
$\delta(t\gk{n+1})\cong_{\cal L}\delta(\gk{m+1}t)$.

\vspace{.2cm}

Then we can prove the following lemma.

\vspace{.2cm}

\noindent $\psi$ {\sc Lemma.}\quad {\it If} $f=g$ {\it in} \Lc, {\it then}
$\psi(f)=\psi(g)$ {\it in} \Lo.

\vspace{.1cm}

\noindent {\it Proof.}\quad
We proceed by induction on the length of the derivation of $f=g$ in \Lc. All the
cases are quite straightforward except when $f=g$ is an instance of
({\it nat} $\varphi$) or ({\it nat} $\gamma$), where we use Remarks I and II.
In case $f=g$ is an instance of ($\varphi\gamma$) we use ({\it cup-cap} 3).
\qed

\vspace{.2cm}

As an immediate corollary we have that if $f\eql g$, then $\psi(f)=\psi(g)$ in
\Lo. Hence we have a map from $\Lc^{\ast}$ to \Lo, which we also
call $\psi$, defined by $\psi([f])=\psi(f)$. Since
\[
\begin{array}{l}
\psi([\mj_0])=\psi(\mj_0)=\mj,
\\
\psi([g\ast f])=\psi(g)\psi(f),
\end{array}
\]
this map is a monoid homomorphism.

Consider next the map $\chi$ from the terms of \Lo\
to the arrow-terms of \Lc\ defined
inductively by

\[
\begin{array}{lcl}
\chi(\mj) & {\mbox{\rm is}} & \mj_0,
\\
\chi(\dk{k}) & {\mbox{\rm is}} & \varphi_{k-1},
\\
\chi(\gk{k}) & {\mbox{\rm is}} & \gamma_{k-1},
\\
\chi(tu) & {\mbox{\rm is}} & \chi(t)\ast\chi(u).
\end{array}
\]
Then we establish the following lemmata.

\vspace{.2cm}

\noindent $\chi$ {\sc Lemma.}\quad {\it If} $t=u$ {\it in} \Lo, {\it then}
$\chi(t)\eql\chi(u)$.

\vspace{.1cm}

\noindent {\it Proof.}\quad We proceed by induction on the length of the derivation of
$t=u$ in \Lo. The cases where $t=u$ is an instance of (1) and (2) are quite
straightforward. If $t=u$ is an instance of ({\it cup}), we have
\[
\varphi_{k-1}\cirk F^{k-l+2}\varphi_{l-1}=F^{k-l}\varphi_{l-1}\cirk\varphi_{k+1},\quad
{\mbox{\rm by ({\it nat} $\varphi$)}}.
\]
If $t=u$ is an instance of ({\it cup-cap} 1), we have
\[
F^{k-l+2}\varphi_{l-1}\cirk \gamma_{k+1}=
\gamma_{k-1}\cirk F^{k-l}\varphi_{l-1},
\quad {\mbox{\rm by ({\it nat} $\gamma$)}}.
\]

We proceed analogously for ({\it cap}), using ({\it nat} $\gamma$),
and for ({\it cup-cap} 2), using ({\it nat} $\varphi$). Finally, if
$t=u$ is an instance of ({\it cup-cap} 3), we have by ($\varphi\gamma$)
\[
\begin{array}{ll}
F\varphi_{k-1}\cirk\gamma_k=\mj_k,
\\
\varphi_{k-1}\cirk F\gamma_{k-2}=\mj_{k-1}.
\end{array}
\]

We have already established that $\eql$ is an equivalence relation that
satisfies ({\it congr} $\ast$). So the lemma follows.
\qed

\vspace{.2cm}

\noindent $\chi\psi$ {\sc Lemma.}\quad {\it For every arrow-term} $f$ {\it of} \Lc\
{\it we have} $\chi(\psi(f))\eql f$.

\vspace{.1cm}

\noindent {\it Proof.}\quad We proceed by induction on the length of $f$. We have
\[
\begin{array}{lcl}
\chi(\psi(\mj_n)) & {\mbox{\rm is}} & \mj_0\eql\mj_n,
\\
\chi(\psi(\varphi_n)) & {\mbox{\rm is}} & \varphi_n,
\\
\chi(\psi(\gamma_n)) & {\mbox{\rm is}} & \gamma_n,
\\[.1cm]
\chi(\psi(Ff)) & {\mbox{\rm is}} & \chi(\psi(f))
\\
               & \eql            & f,\;{\mbox{\rm by the induction hypothesis}}
\\
               & \eql            & Ff,
\\[.1cm]
\chi(\psi(g\cirk f)) & {\mbox{\rm is}} & \chi(\psi(g))\ast\chi(\psi(f))
\\
                     & \eql            & g\cirk f,\;{\mbox{\rm by the induction hypothesis, }} ({\mbox{\it congr }} \ast )\; {\mbox{\rm  and}}
\\
                     &                 & {\mbox{\rm the definition of }} \ast.
\end{array}
\]
\qed

\vspace{.2cm}

By a straightforward induction we can prove also the following lemma.

\vspace{.2cm}

\noindent $\psi\chi$ {\sc Lemma.}\quad {\it For every term} $t$ {\it of} \Lo\
{\it we have that} $\psi(\chi(t))$ {\it is the term} $t$.

\vspace{.2cm}

This establishes that $\Lc^{\ast}$ and \Lo\ are isomorphic monoids.

Let $\Kc^{\ast}$ and $\Mc^{\ast}$ be monoids defined analogously to
$\Lc^{\ast}$, by replacing everywhere $\cal L$ by $\cal K$ and $\cal J$
respectively. Then we can easily extend the foregoing results to establish that
$\Kc^{\ast}$ is isomorphic to the monoid \Ko\ and that $\Mc^{\ast}$ is
isomorphic to the monoid \Mo.

We also have the following lemma.

\vspace{.2cm}

\noindent $\cal L$ {\sc Cancellation Lemma.}\quad {\it In every} $\cal L$-{\it
adjunction, for} $f,g:a\str Fb$ {\it or} $f,g:Fa\str b$, {\it if}
$Ff=Fg$, {\it then} $f=g$.

\vspace{.1cm}

\noindent {\it Proof.}\quad Suppose $Ff=Fg$ for $f,g:a\str Fb$. Then
\[
\begin{array}{rcl}
F\varphi_b\cirk FFf\cirk\gamma_a & = & F\varphi_b\cirk FFg\cirk\gamma_a
\\
F\varphi_b\cirk\gamma_{Fb}\cirk f & = & F\varphi_b\cirk\gamma_{Fb}\cirk g,\;
{\mbox{\rm by ({\it nat} $\gamma$)}}
\\
f & = & g,\;{\mbox{\rm by ($\varphi\gamma$) and ({\it cat} 1)}}.
\end{array}
\]
Suppose $Ff=Fg$ for $f,g:Fa\str b$. Then
\[
\begin{array}{rcl}
\varphi_b\cirk FFf\cirk F\gamma_a & = & \varphi_b\cirk FFg\cirk F\gamma_a
\\
f\cirk\varphi_{Fa}\cirk F \gamma_{a} & = & g\cirk\varphi_{Fa}\cirk F\gamma_{a},\;
{\mbox{\rm by ({\it nat} $\varphi$)}}
\\
f & = & g,\;{\mbox{\rm by ($\varphi\gamma$) and ({\it cat} 1)}}.
\end{array}
\]
\qed

\vspace{.2cm}

As an instance of this lemma we obtain that $Ff=Fg$ implies $f=g$ in \Lc\
provided that for $f,g:m\str n$ we have $m+n>0$. As a matter of fact, this implication
holds for $f,g:0\str 0$ too, but the proofs we know of that fact are
rather involved, and are pretty lengthy. We know two proofs, which
are both based on reducing arrow-terms of \Lc\ to a unique normal form that
corresponds to $c_1^{\alpha}$ of Section 4 when $m+n=0$.
One of these normal forms is based on the normal form for terms of \Lo\ in
Section 4, but with a number of complications brought
in by the types of arrow terms. The other normal form is a composition-free
normal form in a particular language, and reduction to it
(achieved in the style of Gentzen's famous proof-theoretical
cut elimination theorem; see \cite{G69}, Paper 3) is at least as
complicated as reduction to the other normal form. We will omit these proofs,
since the importance of the fact in question, on which we will not rely in the
sequel, does not warrant spending too much on establishing it.

As a corollary of the $\cal L$ Cancellation Lemma, and of previously
established results, we have that for
$f,g:m\str n$ with $m+n>0$
\[
\begin{array}{rcl}
f=g\;{\mbox{\rm in \Lc}} & {\mbox{\rm iff}} & f\eql g
\\
                         & {\mbox{\rm iff}} & \psi(f)=\psi(g)\; {\mbox{\rm in \Lo}}
\\
                         & {\mbox{\rm iff}} & \delta(\psi(f))\cong_{\cal L} \delta(\psi(g)).
\end{array}
\]

Since the $\cal L$ Cancellation Lemma applies also to \Kc, we have exactly analogous
equivalences when $\cal L$ is replaced by $\cal K$. However, for
this replacement we can rather easily lift the restriction $m+n>0$.

\vspace{.2cm}

\noindent \Kc\ {\sc Cancellation Lemma.}\quad
{\it In} \Kc, {\it if} $Ff=Fg$, {\it then} $f=g$.

\vspace{.1cm}

\noindent {\it Proof.}\quad If in $f,g:m\str n$ we have $m+n>0$, we apply the
$\cal L$ Cancellation Lemma. If $m=n=0$, then $f$ is equal either
to $\mj_0$ or
to $\varphi_0\cirk f'$, and $g$ is equal either to $\mj_0$ or
to $\varphi_0\cirk g'$.
Here $f'$ must be $f'' \cirk\gamma_0$ and $g'$ must be $g''\cirk\gamma_0$.
(So we could alternatively consider $f$ and $g$ being equal to
$f'''\cirk\gamma_0$ or $g'''\cirk\gamma_0$, and reason analogously below.)

If both $f=g=\mj_0$, we are done. It is excluded that $f=\mj_0$ while
$g=\varphi_0\cirk g'$. We have that $\psi(f)$ is $\psi(Ff)$ and
$\psi(g)$ is $\psi(Fg)$. Since from $Ff=Fg$ in \Kc\ it follows that
$\delta(\psi(Ff))\cong_{\cal K}\delta(\psi(Fg))$, we have
$\delta(\psi(f))\cong_{\cal K}\delta(\psi(g))$.
But $\delta(\psi(\mj_0))$ is not $\cal K$-equivalent
to $\delta(\psi(\varphi_0\cirk g''\cirk\gamma_0))$,
because, by Remark I, the frieze $\delta(\psi(g''))$ must
be a 2-frieze, from which we obtain that
there is at least one circular
component in $\delta(\psi(\varphi_0\cirk g''\cirk\gamma_0))$.
It is excluded in the same manner that $g=\mj_0$ while $f=\varphi_0\cirk f'$.

If $f=\varphi_0\cirk f'$ and $g=\varphi_0\cirk g'$, then from
$\delta(\psi(\varphi_0\cirk f'))\cong_{\cal K}\delta(\psi(\varphi_0\cirk g'))$,
we conclude $\delta(\psi(f'))\cong_{\cal K} \delta(\psi(g'))$.
This is because $\delta(\psi(\varphi_0\cirk f'))$ and
$\delta(\psi(\varphi_0\cirk g'))$ are both $\cal K$-equivalent to $\delta(c^k)$, for
$k\geq 1$, while $\delta(\psi(f'))$ and $\delta(\psi(g'))$ must both be
$\cal K$-equivalent to $\delta(c^{k-1}\gk{1})$. But $f'$ and $g'$ are of type $0\str 2$,
and hence, by the $\cal L$ Cancellation Lemma, $f'=g'$ in \Kc, from which
we obtain that $\varphi_0\cirk f'=\varphi_0\cirk g'$ in \Kc.
\qed

\vspace{.2cm}

This proof wouldn't go through for \Lc, because $\delta(\psi(f'))$ need not be
$\cal L$-equivalent to $\delta(\psi(g'))$. For example, with $h$ being
$F(\varphi_3\cirk\gamma_3)\cirk\gamma_2\cirk\gamma_0$, we have in \Lc
\[
\varphi_0\cirk\varphi_2\cirk h=\varphi_0\cirk F^2\varphi_0\cirk h,
\]
but $\delta(\psi(\varphi_2\cirk h))$ is not $\cal L$-equivalent to
$\delta(\psi(F^2\varphi_0\cirk h))$, and
$\varphi_2\cirk h=F^2\varphi_0\cirk h$ doesn't hold in \Lc. It holds in \Kc.

The \Kc\ Cancellation Lemma implies that $f\eqk g$ for $f:m\str n$
and $g:k\str l$ could be defined by $F^{k-n}f=g$ in \Kc\ when $n\leq k$, and by
$f=F^{n-k} g$ in \Kc\ when $k\leq n$. So for arbitrary $f,g:m\str n$
we have established that $f=g$ in \Kc\ iff $f\eqk g$.

Since, by Remark I, for $f,g:n\str n$, where $n\in\N$, the friezes $\delta(\psi(f))$
and $\delta(\psi(g))$ are $n$-friezes, we can conclude that \Kn\ is isomorphic to the
monoid of endomorphisms $f:n\str n$ of \Kc. We have this isomorphism for every
$n\in\N$, but the monoids \Kn\ are interesting only when $n\geq 2$.

We could conclude analogously that \Ln\ is isomorphic to the monoid of
endomorphisms $f:n\str n$ of \Lc, relying on the proof of the isomorphism of
\Ln\ with $[\Fn]_{\cal L}$, which we have only indicated, and not given
in Section 9. The $\cal L$ Cancellation Lemma
guarantees this isomorphism for \Ln\ if $n>0$, though,
as we mentioned above,
at the cost of additional arguments this restriction can be lifted.

We have the following lemma for $\cal J$-adjunctions.

\vspace{.2cm}

\noindent $\cal J$ {\sc Cancellation Lemma.}\quad {\it In every}
$\cal J$-{\it adjunction, for f and g arrows of the same type,
if} $Ff=Fg$, {\it then} $f=g$.

\vspace{.1cm}

\noindent {\it Proof.}\quad Take $f,g:a\str b$, for $a$ and $b$ arbitrary objects.
From $Ff=Fg$ we infer
\[
\varphi_b\cirk FFf\cirk\gamma_a=\varphi_b\cirk FFg\cirk\gamma_a,
\]
from which by ({\it nat} $\varphi$) or ({\it nat} $\gamma$), followed by
($\varphi\gamma{\cal J}$) and ({\it cat} 1), we obtain $f=g$.
\qed

\vspace{.2cm}

\noindent By proceeding as in this proof in an arbitrary $\cal L$-adjunction
we can conclude only that if $Ff=Fg$, then $f\cirk \kappa_a=g\cirk \kappa_a$.

The $\cal J$ Cancellation Lemma enables us to show that
for arbitrary $f,g:m\str n$
we have $f=g$ in \Mc\ iff $\psi(f)=\psi(g)$ in \Mo, which means
that we can check
equations of \Mc\ through $\cal K$-equivalence of $\cal J$-friezes.
The monoids \Mn\ are isomorphic to the monoids of endomorphisms of \Mc.

We may now confirm what we stated at the end of the preceding section,
namely, that if
$f=g$ in \Mc\ but not in \Kc, then for some $k,l\geq 0$ such that
$k\neq l$ we have $f\cirk\kappa_0^k=g\cirk\kappa_0^l$ in \Kc. If
$f=g$ doesn't hold \Kc, but it holds in \Mc, then
the friezes $\delta(\psi(f))$ and $\delta(\psi(g))$ differ only
with respect to the number of circular components.

\section{Self-adjunction in \Mat}

Let \Mat\ be the skeleton of the category of finite-dimensional vector
spaces over a number field $\cal F$ with linear transformations as arrows.
A number field is any subfield of the field of complex numbers {\bf C},
and hence it is an extension of the field of rational numbers {\bf Q}.
A skeleton of a category $\cal C$ is any full subcategory $\cal C'$ of
$\cal C$
such that each object of $\cal C$ is isomorphic in $\cal C$ to exactly
one object
of $\cal C'$. Any two skeletons of $\cal C$ are isomorphic categories,
so that,
up to isomorphism, we may speak of {\it the} skeleton of $\cal C$.

More precisely, the objects of the category
\Mat\ are natural numbers (the
dimensions of our vector spaces), an arrow $A:m\str n$
is an $n\times m$ matrix,
composition of arrows $\cirk$ is matrix multiplication, and the identity
arrow $\mj_n:n\str n$ is the $n\times n$ matrix with 1
on the diagonal and 0 elsewhere. (The number 0 is a null object in the
category \Mat, which, as far as we are here interested in this category,
we could as well exclude.)

For much of what we say at the beginning concerning self-adjunction in
\Mat\ it would be enough to assume that the scalars in $\cal F$
are just elements of the commutative
monoid $\langle \N,+,0\rangle$. However,
then we would not have vector spaces, but
something more general, which has no standard name.
Later (see Section 21) we will indeed need that
the scalars make {\bf Q} or an extension of it.

Let $p\in\N^+$, and consider the functor $p\otimes$
from \Mat\ to \Mat\ defined
as follows: for the object $m$ of \Mat\ we have that $p\otimes m$ is $pm$, and
for the arrow $B:m\str n$ of \Mat, i.e. an $n\times m$ matrix $B$, let
$p\otimes B:pm\str pn$ be the Kronecker product $\mj_p\otimes B$ of the
matrices
$\mj_p$ and $B$ (see \cite{J53}, Chapter VII.5, pp. 211-213). It is not difficult
to check that $p\otimes$ is indeed a functor. The essential
properties of the Kronecker product $\otimes$ we will need below are
that $\otimes$ is associative and that
\[
\alpha(A\otimes B)=\alpha A\otimes B= A\otimes\alpha B.
\]

The functor $1\otimes$ is just the identity functor on \Mat.
The interesting functors $p\otimes$ on \Mat\ will have $p\geq 2$.

Let $E_p$ be the
$1\times p^2$ matrix that for $1\leq i,j\leq p$ has the entries
\[E_p(1,(i-1)p+j)=\delta(i,j),\]
where $\delta$ is the Kronecker delta. For example, $E_2$ is $[1\;0\;0\;1]$
and $E_3$ is $[1\;0\;0\;0\;1\;0\;0\;0\;1]$. Let $E_p'$ be the
transpose of $E_p$.
Then $\varphi_m$ is $E_p \otimes \mj_m$, and
$\gamma_m$ is its transpose, i.e. $E_p' \otimes \mj_m$.
We can check that $\varphi$ and $\gamma$ are natural transformations,
which satisfy moreover ($\varphi\gamma$) and ($\varphi\gamma{\cal K}$).
Namely, we can check that
$\langle\Mat,\cirk,\mj,p\otimes,\varphi,\gamma\rangle$ is a
$\cal K$-adjunction.

This self-adjunction is based on the fact that $2\otimes A=A\oplus A$,
where $A\oplus A$ is the sum of matrices

\begin{center}
\begin{picture}(92,70)
{\thicklines
\put(0,7){\line(1,0){3}}
\put(0,7){\line(0,1){66}}
\put(0,73){\line(1,0){3}}
\put(82,7){\line(-1,0){3}}
\put(82,7){\line(0,1){66}}
\put(82,73){\line(-1,0){3}}}

\put(23,25){\makebox(0,0){$0$}}
\put(59,25){\makebox(0,0){$A$}}
\put(23,55){\makebox(0,0){$A$}}
\put(59,55){\makebox(0,0){$0$}}
\linethickness{0.05pt}
\put(5,10){\line(0,1){60}}
\put(41,10){\line(0,1){60}}
\put(77,10){\line(0,1){60}}
\put(5,10){\line(1,0){72}}
\put(5,40){\line(1,0){72}}
\put(5,70){\line(1,0){72}}

\end{picture}
\end{center}

\noindent and behind this sum we have a bifunctor that is both a product and a
coproduct.
The category \Mat\ is a linear category in the sense of \cite{LS}; namely,
in it finite products and coproducts are isomorphic---actually, they
coincide.
Finite products and coproducts coincide in the category of commutative monoids
with monoid homomorphisms, of which the category of vector spaces
over $\cal F$ is a subcategory.
Since we always have that the product bifunctor is right adjoint to the diagonal
functor into the product category, and the coproduct bifunctor is left adjoint
to this diagonal functor, by composing the product bifunctor, which in
\Mat\ coincides with the coproduct bifunctor, with the diagonal
functor we obtain
in \Mat\ a self-adjoint functor.

The category \Mat\ is a strict monoidal category with the bifunctor $\otimes$,
whose unit object is 1. This category is also symmetric monoidal (see
\cite{ML71}, Chapter VII.1,7).

The self-adjunction of $p\otimes$ in \Mat\ is not a $\cal J$-adjunction
for $p\geq 2$, because $\varphi_m\cirk\gamma_m=p\mj_m$ (here the
natural number
$p$ is a scalar, and in $p\mj_m$ the matrix $\mj_m$ is multiplied
by this scalar). However, there is still a possibility to interpret \Mo,
which
is derived from \Mc, in \Mat, as we will see in the next section.

\section{Representing \Mo\ in \Mat}

For $p\in\N^+$, consider the operation $\ast$
on the arrows $A:p^m\str p^n$ and $B:p^k\str p^l$ of \Mat, which is analogous
to the operation $\ast$ of \Lc\ in Section 16:
\[
B\ast A=_{\mbox{\scriptsize\it def}}
\left\{
\begin{array}{ll}
B\cirk(\mj_{p^{k-n}}\otimes A) & {\mbox{\rm if }}n\leq k
\\
(\mj_{p^{n-k}}\otimes B)\cirk A & {\mbox{\rm if }}k\leq n.
\end{array}
\right.
\]

Consider next the map $\eta_p$ from the terms of \Lo\ to the arrows
of \Mat\
defined inductively as follows:
\[
\begin{array}{lcl}
\eta_p(\mj) & {\mbox{\rm is}} & \mj_1:p^0\str p^0,
\\[.05cm]
\eta_p(\dk{k}) & {\mbox{\rm is}} & \varphi_{p^{k-1}}:p^{k+1}\str p^{k-1},
\\[.05cm]
\eta_p(\gk{k}) & {\mbox{\rm is}} & \gamma_{p^{k-1}}:p^{k-1}\str p^{k+1},
\\[.05cm]
\eta_p(tu) & {\mbox{\rm is}} & \eta_p(t)\ast\eta_p(u).
\end{array}
\]

Next, let $A\eqm B$ in \Mat\ iff there are numbers $k,l,m\in\N$ such that
$p^m(\mj_{p^k}\otimes A)=\mj_{p^l}\otimes B$ or
$\mj_{p^k}\otimes A=p^m(\mj_{p^l}\otimes B)$
in \Mat. The relation $\eqm$ is an equivalence relation on the arrows of
\Mat, congruent with respect to the operations $\ast$ and
$\mj_p\otimes$. We can then prove the following lemma.

\vspace{.2cm}

\noindent $\eta_p$ {\sc Lemma.}\quad {\it For} $p\geq 2$ {\it we have}
$t=u$ {\it in} \Mo\ {\it iff} $\eta_p(t)\eqm\eta_p(u)$ {\it in} \Mat.

\vspace{.1cm}

\noindent {\it Proof.} From left to right we proceed, in principle, by
induction on the length of derivation of $t=u$ in \Mo. However,
most cases in this induction are already covered by the $\chi$ Lemma
of Section 16, and by our having established in the
preceding section that \Mat\ is
a $\cal K$-adjunction. The only case specific for \Mo, namely,
when $t=u$ is an
instance of $\dk{k}\gk{k}=\mj$, is covered by the fact that
$\varphi_{p^{k-1}}\cirk\gamma_{p^{k-1}}=p\mj_{p^{k-1}}.$

To prove the lemma from right to left, suppose that we don't have $t=u$
in \Mo, but $\eta_p(t)\eqm\eta_p(u)$ in \Mat. Then, by the
left-to-right direction of the lemma, which we have just
established, we should have in \Mat\
an $\eta_p$ image of \Mo\ extended with $t=u$. By the results of
Section 12, we should have that
$\eta_p(\gk{i}\dk{i})\eqm \eta_p(\mj)$ in \Mat.
We have, however, that $\eta_p(\gk{i}\dk{i})$ is
$\gamma_{p^{i-1}}\cirk\varphi_{p^{i-1}}$,
and for no $k,l,m\in\N$
we can have
$p^m(\mj_{p^k}\otimes (\gamma_{p^{i-1}}\cirk\varphi_{p^{i-1}}))=\mj_{p^l}$ or
$\mj_{p^k}\otimes (\gamma_{p^{i-1}}\cirk\varphi_{p^{i-1}})=p^m\mj_{p^l}$,
provided $p \geq 2$.
From this the lemma follows.
\qed

\vspace{.2cm}

An alternative proof of this lemma could be obtained by relying on
the results of \cite{DP02}.

\section{Representing \Kc\ in \Mat}

Let $p\in\N^+$, and consider
the category \Kc\ of the free $\cal K$-adjunction generated by the object $0$
(see Section 15). We define inductively a functor $H_p$ from \Kc\
to \Mat\ in the following manner:
\[
\begin{array}{lcl}
H_p(0) & {\mbox{\rm is}} & 1=p^0,
\\[.05cm]
H_p(m+1) & {\mbox{\rm is}} & pH_p(m)=p^{m+1},
\\[.2cm]
H_p(\mj_m) & {\mbox{\rm is}} & \mj_{p^m}:p^m\str p^m,
\\[.05cm]
H_p(\varphi_m) & {\mbox{\rm is}} & \varphi_{p^m}:p^{m+2}\str p^m,
\\[.05cm]
H_p(\gamma_m) & {\mbox{\rm is}} & \gamma_{p^m}:p^m\str p^{m+2},
\\[.05cm]
H_p(Ff) & {\mbox{\rm is}} & \mj_p\otimes H_p(f):p^{m+1}\str p^{n+1},
\quad {\mbox{\rm for }}f:m\str n,
\\[.05cm]
H_p(g\cirk f) & {\mbox{\rm is}} & H_p(g)\cirk H_p(f).
\end{array}
\]
That this defines a functor indeed follows from the fact that
$\langle\Mat,\cirk,\mj,p\otimes,\varphi,\gamma\rangle$ is a
$\cal K$-adjunction, as established in Section 17.

The function $H_p:\Lc\str\Mat$ on arrow-terms of \Lc, which we have
above, is not obtained by
composing $\psi:\Lc\str\Lo$ of Section 16 and
$\eta_p:\Lo\str\Mat$ of the previous section, but we can check by induction
on the length of $f$ that $H_p(f)\eqm\eta_p(\psi(f))$ in \Mat.

The functor $H_1$ is not faithful, since for every arrow $f$ of \Kc\
we have $H_1(f)=\mj_1$, but for $p\geq 2$ the functors $H_p$ are faithful.
As a
matter of fact, these functors, which are one-one on objects,
are one-one on arrows.
This is shown by the following lemma.

\vspace{.2cm}

\noindent {\sc Faithfulness of} $H_p$.\quad {\it For f and g arrow-terms of} \Lc\
{\it of the same type and} $p \geq 2$,
{\it if} $H_p(f)=H_p(g)$
{\it in} \Mat, {\it then} $f=g$ {\it in} \Kc.

\vspace{.1cm}

\noindent {\it Proof.}\quad Suppose $H_p(f)=H_p(g)$ in \Mat, but not $f=g$ in \Kc.
If $f=g$ in \Mc, then, as we have seen
at the very end of Section 16,
for some $k,l\geq 0$ such that $k\neq l$ we have
$f\cirk\kappa_0^k=g\cirk\kappa_0^l$ in \Kc. But then in \Mat\ we have
\[
p^kH_p(f)=p^lH_p(f),
\]
which is impossible, since $H_p(f)$ is never a zero matrix.

So we don't have $f=g$ in \Mc. Hence,
according to what we have established before the end of Section 16,
we don't have $\psi(f)=\psi(g)$ in \Mo.
Then, by the $\eta_p$ Lemma of the previous section, we don't have
$\eta_p(\psi(f))\eqm\eta_p(\psi(g))$ in \Mat. However, from
$H_p(f)=H_p(g)$ in \Mat\ it follows that
$\eta_p(\psi(f))\eqm\eta_p(\psi(g))$ in \Mat, which yields a contradiction.
\qed

\vspace{.2cm}

So in \Mat\ we have an isomorphic representation of \Kc. We also have
for every $n\in\N$ isomorphic
representations of the monoids \Kn\ as monoids of endomorphisms of $p^n$,
provided $p\geq 2$.

Our proof of the faithfulness of these representations of \Kn\
relies on the maximality result of Section 12.
An alternative proof is obtained either as in \cite{DP02}, or
by relying on the faithfulness result of \cite{J94} (Section 3) and
\cite{DKP}, as mentioned in the Introduction.

\section{The algebras $\End(p^n)$}

Let $\End(p^n)$ be the set of all
endomorphisms $A:p^n\str p^n$ in \Mat, i.e. of all $p^n\times p^n$
matrices in \Mat.
Let us first consider $\End(p^n)$ when $p$ is 2.
We have remarked at the end of the preceding section that we have in
$\End(2^n)$ an isomorphic representation of the monoid \Kn. Let
us denote by $h_k^n$ the representation of the diapsis $h_k=\gk{k}\dk{k}$
of \Kn\ in $\End(2^n)$. The matrix $h_k^n$ is $\mj_{2^{n-k-1}}
\otimes(\gamma_{2^{k-1}}\cirk\varphi_{2^{k-1}})$.
To define $\gamma_{2^{k-1}}$ and $\varphi_{2^{k-1}}$ we need
the matrix $E_2$, namely $[1\; 0\; 0\; 1]$, and its transpose $E_2'$.
The matrix $E_2'\cirk E_2$ is
\[
\left [
\begin{array}{cccc}
1 & 0 & 0 & 1 \\ 0 & 0 & 0 & 0 \\ 0 & 0 & 0 & 0 \\ 1 & 0 & 0 & 1
\end{array}
\right ]
\]
So $h_k^n=\mj_{2^{n-k-1}}\otimes(E_2'\cirk E_2)\otimes \mj_{2^{k-1}}$.

For example, in $\End(2^2)$ the diapsis $h_1$ is represented by the matrix
$h_1^2$, which is $\gamma_1\cirk\varphi_1=E_2'\cirk E_2$.
In $\End(2^3)$ the diapsides $h_1$ and
$h_2$ are represented by $h_1^3=\mj_2\otimes h_1^2=\mj_2\otimes
(E_2'\cirk E_2)$ and
$h_2^3=\gamma_2\cirk\varphi_2=(E_2'\cirk E_2)\otimes\mj_2$.

Every $n\times m$ matrix $A$ whose entries are only 0 and 1 may be
identified with
a binary relation $R_A\subseteq n\times m$ such that $A(i,j)=1$
iff $(i,j)\in R_A$.
Every binary relation may of course be drawn as a bipartite graph. Here are a
few examples of such graphs for matrices we have introduced up to now, with
$p=2$:

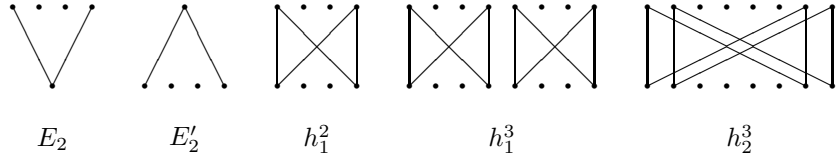
\begin{center}
\begin{picture}(310,90)
\put(0,70){\circle*{2}}
\put(10,70){\circle*{2}}
\put(20,70){\circle*{2}}
\put(30,70){\circle*{2}}
\put(65,70){\circle*{2}}
\put(100,70){\circle*{2}}
\put(110,70){\circle*{2}}
\put(120,70){\circle*{2}}
\put(130,70){\circle*{2}}
\put(150,70){\circle*{2}}
\put(160,70){\circle*{2}}
\put(170,70){\circle*{2}}
\put(180,70){\circle*{2}}
\put(190,70){\circle*{2}}
\put(200,70){\circle*{2}}
\put(210,70){\circle*{2}}
\put(220,70){\circle*{2}}
\put(240,70){\circle*{2}}
\put(250,70){\circle*{2}}
\put(260,70){\circle*{2}}
\put(270,70){\circle*{2}}
\put(280,70){\circle*{2}}
\put(290,70){\circle*{2}}
\put(300,70){\circle*{2}}
\put(310,70){\circle*{2}}
\put(15,40){\circle*{2}}
\put(50,40){\circle*{2}}
\put(60,40){\circle*{2}}
\put(70,40){\circle*{2}}
\put(80,40){\circle*{2}}
\put(100,40){\circle*{2}}
\put(110,40){\circle*{2}}
\put(120,40){\circle*{2}}
\put(130,40){\circle*{2}}
\put(150,40){\circle*{2}}
\put(160,40){\circle*{2}}
\put(170,40){\circle*{2}}
\put(180,40){\circle*{2}}
\put(190,40){\circle*{2}}
\put(200,40){\circle*{2}}
\put(210,40){\circle*{2}}
\put(220,40){\circle*{2}}
\put(240,40){\circle*{2}}
\put(250,40){\circle*{2}}
\put(260,40){\circle*{2}}
\put(270,40){\circle*{2}}
\put(280,40){\circle*{2}}
\put(290,40){\circle*{2}}
\put(300,40){\circle*{2}}
\put(310,40){\circle*{2}}

\put(0,70){\line(1,-2){15}}
\put(30,70){\line(-1,-2){15}}

\put(50,40){\line(1,2){15}}
\put(80,40){\line(-1,2){15}}

\put(100,40){\line(0,1){30}}
\put(100,40){\line(1,1){30}}
\put(130,40){\line(-1,1){30}}
\put(130,40){\line(0,1){30}}

\put(150,40){\line(0,1){30}}
\put(150,40){\line(1,1){30}}
\put(180,40){\line(-1,1){30}}
\put(180,40){\line(0,1){30}}
\put(190,40){\line(0,1){30}}
\put(190,40){\line(1,1){30}}
\put(220,40){\line(-1,1){30}}
\put(220,40){\line(0,1){30}}

\put(240,40){\line(0,1){30}}
\put(240,40){\line(2,1){60}}
\put(300,40){\line(-2,1){60}}
\put(300,40){\line(0,1){30}}
\put(250,40){\line(0,1){30}}
\put(250,40){\line(2,1){60}}
\put(310,40){\line(-2,1){60}}
\put(310,40){\line(0,1){30}}

\put(15,20){\makebox(0,0){$E_2$}}
\put(65,20){\makebox(0,0){$E_2'$}}
\put(115,20){\makebox(0,0){$h^2_1$}}
\put(185,20){\makebox(0,0){$h^3_1$}}
\put(275,20){\makebox(0,0){$h^3_2$}}

\end{picture}
\end{center}

In $\End(2^4)$ we have $h_1^4=\mj_2\otimes h_1^3$,
$h_2^4=\mj_2\otimes h_2^3$ and
$h_3^4=\gamma_4\cirk\varphi_4=(E_2'\cirk E_2)\otimes\mj_{2^2}$, etc. for
$\End(2^5)$, $\End(2^6), \ldots$ In $\End(2^n)$ the unit $\mj$ of \Kn\
is represented by the $2^n \times 2^n$ identity matrix
$\mj_{2^n}$, whose entries are
$\mj_{2^n}(i,j)=\delta (i,j)$, where $\delta$ is Kronecker's delta.
As usual, we denote this matrix also by $I$.
The circle $c$ is
represented by $\mj_{2^n}\otimes(E_2\cirk E_2')=\mj_{2^n}\otimes [2]=
2\mj_{2^n}=2I$.

We proceed analogously when $p>2$ in $\End(p^n)$. Then
$h_k^n=\mj_{p^{n-k-1}}\otimes(\gamma_{p^{k-1}}\cirk \varphi_{p^{k-1}})=
\mj_{p^{n-k-1}}\otimes(E'_p\cirk E_p)\otimes \mj_{p^{k-1}}$, the
unit matrix $I$ is the identity matrix $\mj_{p^n}$, and the circle is
represented by $pI$.

If $\mj^n, h^n_1, \ldots,h^n_{n-1}, c$ denote the 0-1 matrices
we have assigned to
these expressions, then these matrices
satisfy the equations of \Kn, with multiplication
being matrix multiplication.
If $\mj^n,h^n_1,\ldots,h^n_{n-1}$ denote the corresponding binary relations,
then for multiplication being composition of binary relations the
equations of \Mn\ are satisfied.

Composition of binary relations is easy to read
from bipartite graphs. Here is an
example:

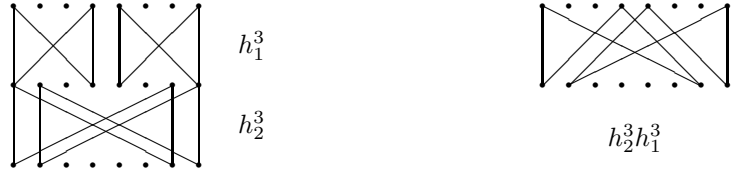
\begin{center}
\begin{picture}(270,100)
\put(0,20){\circle*{2}}
\put(10,20){\circle*{2}}
\put(20,20){\circle*{2}}
\put(30,20){\circle*{2}}
\put(40,20){\circle*{2}}
\put(50,20){\circle*{2}}
\put(60,20){\circle*{2}}
\put(70,20){\circle*{2}}

\put(0,50){\circle*{2}}
\put(10,50){\circle*{2}}
\put(20,50){\circle*{2}}
\put(30,50){\circle*{2}}
\put(40,50){\circle*{2}}
\put(50,50){\circle*{2}}
\put(60,50){\circle*{2}}
\put(70,50){\circle*{2}}

\put(0,80){\circle*{2}}
\put(10,80){\circle*{2}}
\put(20,80){\circle*{2}}
\put(30,80){\circle*{2}}
\put(40,80){\circle*{2}}
\put(50,80){\circle*{2}}
\put(60,80){\circle*{2}}
\put(70,80){\circle*{2}}

\put(200,50){\circle*{2}}
\put(210,50){\circle*{2}}
\put(220,50){\circle*{2}}
\put(230,50){\circle*{2}}
\put(240,50){\circle*{2}}
\put(250,50){\circle*{2}}
\put(260,50){\circle*{2}}
\put(270,50){\circle*{2}}

\put(200,80){\circle*{2}}
\put(210,80){\circle*{2}}
\put(220,80){\circle*{2}}
\put(230,80){\circle*{2}}
\put(240,80){\circle*{2}}
\put(250,80){\circle*{2}}
\put(260,80){\circle*{2}}
\put(270,80){\circle*{2}}

\put(0,50){\line(0,1){30}}
\put(0,50){\line(1,1){30}}
\put(30,50){\line(-1,1){30}}
\put(30,50){\line(0,1){30}}
\put(40,50){\line(0,1){30}}
\put(40,50){\line(1,1){30}}
\put(70,50){\line(-1,1){30}}
\put(70,50){\line(0,1){30}}

\put(0,20){\line(0,1){30}}
\put(0,20){\line(2,1){60}}
\put(60,20){\line(-2,1){60}}
\put(60,20){\line(0,1){30}}
\put(10,20){\line(0,1){30}}
\put(10,20){\line(2,1){60}}
\put(70,20){\line(-2,1){60}}
\put(70,20){\line(0,1){30}}

\put(200,50){\line(0,1){30}}
\put(200,50){\line(1,1){30}}
\put(210,50){\line(1,1){30}}
\put(210,50){\line(2,1){60}}
\put(260,50){\line(-2,1){60}}
\put(260,50){\line(-1,1){30}}
\put(270,50){\line(-1,1){30}}
\put(270,50){\line(0,1){30}}

\put(90,35){\makebox(0,0){$h^3_2$}}
\put(90,65){\makebox(0,0){$h^3_1$}}
\put(235,30){\makebox(0,0){$h^3_2h^3_1$}}

\end{picture}
\end{center}

\noindent By so composing binary relations we
can assign to every element of \Mn\ a binary
relation, and then from this binary relation we can recover
the 0-1 matrix assigned to our element of \Kn.

However, $\End(p^n)$ is a richer structure than \Kn. It is an
associative $\cal F$ algebra under matrix addition +, the multiplication
of a matrix $A$ by a scalar $\alpha$ (which is written $\alpha A$)
and matrix
multiplication, which we continue to write as composition $\cirk$.
We will consider in the next section representations of braid groups in the
algebras $\End(p^n)$.

The representation of \Kn\ in $\End(p^n)$ we dealt with above is
obtained by restricting to \Kn\ the orthogonal group case of
Brauer's representation of Brauer algebras from \cite{B37} (see also
\cite{W88}, Section 3, and \cite{J94}, Section 3).

\section{Representing braid groups in $\End(p^n)$}

The {\it braid group} \Bn\ has for every $k\in \{1,\ldots,n-1\}$ a generator $\sigma_k$.
The number $n$ here could in principle be any natural number, but, as
for \Ln, \Kn\ and \Mn, the interesting groups \Bn\ have $n\geq 2$. When
$n$ is 0 or 1, we have no generators $\sigma_k$. The terms of \Bn\ are
obtained from these generators and \mj\ by closing under inverse $^{-1}$
and multiplication. The following equations are assumed for \Bn:
\[
\begin{array}{ll}
{\makebox[1cm][l]{$(1)$}} & {\makebox[6cm][l]{$\mj t =t\mj=t,$}}
\\[.05cm]
(2) & t(uv)=(tu)v,
\\[.05cm]
(3) & tt^{-1}=t^{-1}t=\mj,
\\[.05cm]
(\sigma 0) & \sigma_i\sigma_j=\sigma_j\sigma_i,\quad {\mbox{\rm for }}|i-j|\geq 2,
\\[.05cm]
(\sigma 3) & \sigma_i\sigma_{i+1}\sigma_i=\sigma_{i+1}\sigma_i\sigma_{i+1}.
\end{array}
\]
We can replace $(3)$ by
\[
\begin{array}{ll}
{\makebox[1cm][l]{$(3.1)$}} & {\makebox[6cm][l]{$\mj^{-1}=\mj,$}}
\\[.05cm]
(3.2) & (tu)^{-1}=u^{-1}t^{-1},
\\[.05cm]
(\sigma 2) & \sigma_i\sigma_i^{-1}=\sigma_i^{-1}\sigma_i=\mj.
\end{array}
\]
(In naming the equations $(\sigma 3)$ and $(\sigma 2)$ we paid attention
to the fact that $(\sigma 3)$ corresponds to the {\it third}
Reidemeister move, and $(\sigma 2)$ to the {\it second}.)

Inspired by the bracket equations (see \cite{KL}, pp. 11, 15, and
references therein), we define inductively as follows a map $\rho$
from the terms of \Bn\ to $\End(p^n)$:

\[
\begin{array}{lcl}
\rho(\mj)\;{\mbox{\rm and }}\rho(\mj^{-1}) & {\mbox{\rm are}} & I,
\\[.05cm]
\rho(\sigma_i) & {\mbox{\rm is}} & \alpha_i h_i^n+\beta_iI,
\\[.05cm]
\rho(\sigma_i^{-1}) & {\mbox{\rm is}} & \alpha_i'h_i^n+\beta_i'I,
\\[.05cm]
\rho(tu)  & {\mbox{\rm is}} & \rho(t)\cirk\rho(u),
\\[.05cm]
\rho((tu)^{-1}) & {\mbox{\rm is}} & \rho(u^{-1})\cirk\rho(t^{-1}).
\end{array}
\]
We will now find conditions for $\alpha_i$, $\alpha_i'$,
$\beta_i$ and $\beta_i'$
sufficient to make $\rho$ a group homomorphism from \Bn\ to $\End(p^n)$.

The equations $(1)$, $(2)$, $(3.1)$ and $(3.2)$ are always satisfied.
The equation $(\sigma 0)$ will also be satisfied always, because
$h_i^n\cirk h_j^n=h_j^n\cirk h_i^n$ when $|i-j|\geq 2$. This is the equation
$(h1)$ of \Kn.

Consider now the equation $(\sigma 3)$. For
$\rho(\sigma_i\sigma_{i+1}\sigma_i)=\rho(\sigma_{i+1}\sigma_i\sigma_{i+1})$
to hold in $\End(p^n)$, we compute that it is sufficient if we have
$\beta_i=\beta_{i+1}$, $\alpha_i=\alpha_{i+1}$, $\beta_i'=\beta_{i+1}'$
and $\alpha_i'=\alpha_{i+1}'$, so that
\[
\begin{array}{lcl}
\rho(\sigma_i) & {\mbox{\rm is}} & \alpha h_i^n+\beta I,
\\[.05cm]
\rho(\sigma_i^{-1}) & {\mbox{\rm is}} & \alpha' h_i^n+\beta'I,
\end{array}
\]
together with
\[
p=-\alpha\beta^{-1}-\alpha^{-1}\beta.
\]

Consider next the equation $(\sigma 2)$. We will have always that
$\rho(\sigma_i\sigma_i^{-1})=\rho(\sigma_i^{-1}\sigma_i)$
in $\End(p^n)$, while in order that $\rho(\sigma_i\sigma_i^{-1})=I$ it is
sufficient that
\[
(p\alpha\alpha'+\alpha\beta'+\alpha'\beta)h_i^n+\beta\beta'I=I,
\]
which yields
\[
\begin{array}{l}
\beta'=\beta^{-1},
\\[.05cm]
p=-(\alpha')^{-1}\beta'-\alpha^{-1}\beta.
\end{array}
\]
Since for $(\sigma 3)$ we required
$p=-\alpha\beta^{-1}-\alpha^{-1}\beta$,
we obtain $\alpha'=\alpha^{-1}$. So with the clauses
\[
\begin{array}{l}
\rho(\sigma_i)=\alpha h_i^n +\beta I,
\\[.05cm]
\rho(\sigma_i^{-1})=\alpha^{-1} h_i^n +\beta^{-1}I,
\\[.05cm]
p=-\alpha\beta^{-1}-\alpha^{-1}\beta,\;{\mbox{\rm i.e. }}
\alpha\beta^{-1}=(-p\pm\sqrt{p^2-4})/2,
\end{array}
\]
which amount to the clauses:
\[
\begin{array}{l}
\rho(\sigma_i)=\alpha(h_i^n+(2/(-p\pm\sqrt{p^2-4}))I),
\\[.05cm]
\rho(\sigma_i^{-1})=\alpha^{-1}(h_i^n+((-p\pm\sqrt{p^2-4})/2)I),
\end{array}
\]
we obtain that $\rho$ is a homomorphism from \Bn\ to $\End(p^n)$, i.e.
a representation of \Bn\ in $\End(p^n)$.

The conditions we have found sufficient to make $\rho$ a representation
of \Bn\ in $\End(p^n)$ are also necessary, since the 0-1 matrices in our
representation of \Kn\ in $\End(p^n)$ are linearly independent (see
\cite{J94}, Section 3, and also \cite{DKP} for an elementary self-contained
proof). Actually,
for the necessity of our conditions we have to prove linear independence
just for the 0-1 matrices in the
representation of ${\cal K}_3$. Then we have only five of these
matrices, whose linear independence one can rather easily check
in case $p$ is equal to 2 or 3 by listing them all, and by finding
for each an entry with 1 where all the others have 0. In \cite{J94} and
\cite{DKP}
linear independence is established for every $n\geq 2$ in \Kn\ and every
$p\geq 2$.

The clauses of the bracket equations
\[
\begin{array}{l}
\rho(\sigma_i)=\alpha h_i^n+\alpha^{-1}I,
\\[.05cm]
\rho(\sigma_i^{-1})=\alpha^{-1} h_i^n+\alpha I,
\\[.05cm]
p=-\alpha^2-\alpha^{-2},
\end{array}
\]
are obtained from ours by requiring that $\beta=\alpha^{-1}$,
and by not requiring as we do that $p$ be a natural number. So
our representation is in a certain sense
more general, but it requires that $p$ be a natural number.

If $p$ is 2 in our representation, then we obtain that
$\beta=-\alpha$. In this case, however, $\rho(\sigma_i\sigma_i)=
\rho(\sigma_{i+1}\sigma_{i+1})=\alpha^2I$, and the representation
is not faithful.
Is this representation faithful for $p>2$?
(The question whether the representation of braid groups in
Temperley-Lieb algebras based on the bracket equations
is faithful is raised in \cite{J00}.)

\end{document}